\documentclass[12pt]{article}
\usepackage{amssymb}
\usepackage[french]{babel}
\usepackage[latin1]{inputenc}
\usepackage[all]{xy}
\newtheorem{thm}{Th\'eor\`eme}[section]

\newtheorem{q}[thm]{Question}

\newtheorem{cor}[thm]{Corollaire}
\newtheorem{lem}[thm]{Lemme}
\newtheorem{prop}[thm]{Proposition}
\newtheorem{defn}[thm]{D\'efinition}
\newtheorem{rem}[thm]{Remarque}

\def\Ltens{{\mathop\otimes\limits^{L}}}

\def\limproj{{\displaystyle\lim_{\longleftarrow}}}
\def\limind{{\displaystyle\lim_{\longrightarrow}}}
\def\hfl#1#2{\smash{\mathop{\hbox to
10mm{\rightarrowfill}}\limits^{\scriptstyle#1}_{\scriptstyle#2}}}
\def\hflrev#1#2{\smash{\mathop{\hbox to
10mm{\leftarrowfill}}\limits^{\scriptstyle#1}_{\scriptstyle#2}}}
\def\hflcourte#1#2{\smash{\mathop{\hbox to
3mm{\rightarrowfill}}\limits^{\scriptstyle#1}_{\scriptstyle#2}}}
\def\hflrevcourte#1#2{\smash{\mathop{\hbox to
3mm{\leftarrowfill}}\limits^{\scriptstyle#1}_{\scriptstyle#2}}}
\def\injfl#1#2{\smash{\mathop{\hookrightarrow}\limits^{\scriptstyle#1}_{\scriptstyle#2}}}
\def\surjfl#1#2{\smash{\mathop{\twoheadrightarrow}\limits^{\scriptstyle#1}_{\scriptstyle#2}}}

\newbox\frogdown
\newcommand{\hookdownarrow}{\setbox\frogdown=\hbox to 0pt{\hss $\displaystyle \downarrow $\hss }\vrule height0pt width0pt depth1.5\ht\frogdown
\setlength{\unitlength}{0.4pt}
\begin{picture}(10,5)
\put(5,0){\oval(10,10)[t]}
\end{picture}
\lower\ht\frogdown\box\frogdown}

\def\nekovar{Nekov\'a\v{r}}

\def\G{\Gamma}

\def\Z{\textbf{$\mathbb Z$}}

\def\limproj{{\displaystyle\lim_{\longleftarrow}}}

\def\X{\mathcal X}

\def\C{\mathcal C}

\def\Q{\mathbb Q}
\def\N{{\mathbb N}}
\def\Zp{{\mathbb Z_p}}
\def\Qp{{\mathbb Q_p}}
\def\L{\Lambda}

\def\Xi{\X^{(-i)}}

\def\Z{\mathbb Z}

\def\vfl#1#2{\llap{$\scriptstyle #1$}\left\downarrow\vbox to
6mm{}\right.\rlap{$\scriptstyle #2$}}

\def\vflcourte#1#2{\llap{$\scriptstyle #1$}\left\downarrow\vbox to
2mm{}\right.\rlap{$\scriptstyle #2$}}
\def\vflupcourte#1#2{\llap{$\scriptstyle #1$}\left\uparrow\vbox to
2mm{}\right.\rlap{$\scriptstyle #2$}}

\def\diagram#1{\def\normalbaselines{\baselineskip=0pt\lineskip=10pt\lineskiplimit=1pt}
  \matrix{#1}}
\def\build#1_#2^#3{\mathrel{
\mathop{\kern 0pt#1}\limits_{#2}^{#3}}} \font\totodix=eufm10 scaled 1000 \font\totosept=eufm7 scaled 1000 \font\totocinq=eufm5 scaled 1000
\newfam\totofam \textfont\totofam=\totodix
\scriptfont\totofam=\totosept \scriptscriptfont\totofam=\totocinq

\begin{document}

\title{Sur la dualité et la descente d'Iwasawa}

%

\author{par David Vauclair\vspace{0.2cm}\\  \hspace{4cm}{\small\textit{dédié à T. Nguyen Quang Do.}}}

\date{2 octobre 2008}

\maketitle

\begin{abstract}
Guided by the concrete examples of cyclotomic units and the ideal class group in cyclotomic Iwasawa theory, we develop a general tool for studying
descent and codescent, with a special interest in relating the two of them.

Given any ``normic system'' $A=(A_n)$ (that is a collection of Galois modules plus additional data), attached to a fixed $p$-adic Lie extension with
Iwasawa algebra $\L$, we mainly show that there is a natural morphism $$R\limproj A_n\rightarrow RHom_\L(RHom_\Zp(\limind A_n,\Zp),\L)$$ which can be
given a functorial cone measuring the defect of descent as well as the defect of codescent (for the $A_n$'s). Thanks to a sharpening of the usual
Poincaré duality, this results in an enlightening relation between these two.

We show in great detail how known results in the cyclotomic situation fit into this setting, and give a generalization to multiple $\Zp$-extensions.
\end{abstract} \vspace{0.5cm}

\noindent {\footnotesize AMS Classification 2000: 11R34, 13C05, 11R23} 

\newpage

\tableofcontents

\section{Introduction}

Soit $F$ un corps de nombres, $p$ un nombre premier, et $F_\infty=\cup F_n$ une extension galoisienne dont le groupe de Galois, noté $\G$, est un
groupe de Lie $p$-adique de dimension cohomologique finie. Supposons donnée une collection $A=(A_n)$, où $A_n$ est un $\Zp[Gal(F_n/F)]$-module (en
pratique d'origine arithmétique, par exemple groupe de classes, unités, groupe de Selmer...), munie à la fois d'une structure de système inductif et
de système projectif, ces deux structures vérifiant une condition naturelle de compatibilité normique. A une telle famille de données (appelée
système normique dans le texte, cf \ref{defIw}), la théorie d'Iwasawa suggère d'associer deux $\L=\Zp[[\G]]$-modules: $A_\infty:=\limind A_n$ et
$X_\infty:=\limproj A_n$. Idéalement, on aimerait pouvoir reconstruire la collection des $A_n$, au moins asymptotiquement (ie. pour $n$ suffisament
grand) à partir de la seule donnée de $X_\infty$ et/ou $A_\infty$. Dans cette optique, il est naturel d'étudier les applications de descente
$j_n:A_n\rightarrow (A_\infty)^{\G_n}$ et de codescente $k_n:(X_\infty)_{\G_n}\rightarrow A_n$.

\begin{q}\label{q1} Que peut-on dire des noyaux et conoyaux de $j_n$ (resp.  $k_n$) pour $n$ suffisamment grand?
En reste-t-il une trace sur $X_\infty$ (resp.  $A_\infty$)?
\end{q}

Malheureusement, la seule donnée de $A_\infty$ et $X_\infty$ ne suffit pas à caractériser un système normique, même asymptotiquement. Aussi on peut
se demander de quelles structures supplémentaires sont munis les couples $(X_\infty,A_\infty)$ provenant d'un système normique. Dans cette direction:

\begin{q}\label{q2} Quel lien existe-t-il entre $X_\infty$ et $A_\infty$ (ou plutôt son dual)?
\end{q}

Malgré l'intérêt pratique évident d'une réponse précise et générale à ces deux questions, seuls peu d'efforts semblent leur avoir été consacrés dans
cette généralité. On peut citer \cite{Gr} pour une étude du noyau et conoyau de $k_n$ (du point de vue des modules de normes universelles) dans le
cas $\G\simeq\Zp$, lorsque $A$ vérifie un certain nombre d'hypothèses restrictives, dont celle que les $j_n$ sont des isomorphismes. Dans le même
ordre d'idées, citons ausi \cite{KY}. Concernant la seconde question, le cas où les $k_n$ sont des isomorphismes est étudié de manière systématique
dans \cite{J3}.

En revanche, l'étude de ces deux questions pour un système normique donné est une activité fréquente en théorie d'Iwawasa, si bien qu'on dispose de
nombreux exemples qui peuvent servir de guide à une théorie générale. Inspiré par les méthodes de \cite{J3} et \cite{Ne}, on propose dans cet article
un cadre pour une telle théorie, et on établit quelques résultats pratiques qui permettent l'unification, et la généralisation de nombreux résultats
connus.

Pour avoir une première idée de ce qu'on peut espérer, on s'inspire essentiellement de deux études parallèles (\cite{LFMN} pour le groupe de classes
et \cite{NL} pour les unités cyclotomiques) dans le cas où $F_\infty=\cup F_n$ est la $\Zp$-extension cyclotomique d'un corps de nombres $F$
($\G\simeq \Zp$). On peut reformuler quelques résultats de la manière suivante (cf \ref{resAB} et \ref{resEC} pour des résultats plus exhaustifs):

\begin{prop} \label{fait1} (\cite{I1}, \cite{LFMN}) Si $A_n=\C l(\mathcal O_{F_n}[{1\over p}])\otimes\Zp$, alors il existe une flèche naturelle
$$X_\infty\mathop\rightarrow\limits^{\alpha_j} Ext^1_\L(Hom_\Zp(A_\infty,\Qp/\Zp),\L)$$ dont le noyau et le conoyau sont pseudo-isomorphes.
Si ceux-ci sont finis, alors $Ker\,\alpha_j\simeq\limproj Ker\, j_n$ et $Coker\alpha_j\simeq \limproj Coker\, j_n\simeq \limind Ker\, k_n$.
\end{prop}

\begin{prop} \label{fait2} (\cite{NL}, \cite{KN}, \cite{Be}, \cite{BN}) Supposons $F$ réel, abélien sur $\Q$ et soit $A_n=C_n$ le
$\Zp$-module des unités cyclotomiques à la Sinnott de $F_n$ (cf. \cite{Si}). Alors il existe une flèche naturelle, injective et pseudo-surjective:
$$X_\infty\mathop\rightarrow\limits^{\alpha_j} Hom_\L(Hom_\Zp(A_\infty,\Zp),\L)$$ De plus $Coker\,\alpha_j\simeq\limproj Coker\, j_n\simeq \limind
Ker\, k_n$. \end{prop}

Ainsi reformulés, ces résultats suggèrent que pour $A$ quelconque, on devrait pouvoir construire une flèche
$$\alpha_j:X_\infty\rightarrow RHom_\L(X_\infty^*,\Zp))$$ avec $X_\infty^*:=RHom_\Zp(A_\infty,\Zp)$, dont ``le'' cône soit relié à la fois au défaut
de descente et de codescente. La première partie de l'article est consacrée à la construction d'une telle flèche. \vspace{0.1cm}

Expliquons brièvement les idées de la construction. D'abord, il convient de préciser la notion de défaut de (co)descente: si $A=(A_n)$ est un système
normique et $A_\infty=\limind A_n$, il est possible de construire un objet de $R\underline \G(A_\infty)\in D^b(_{\underline \G}\C)$, qui relève la
collection des $R\G(\G_n,A_\infty)$ usuels. Celui-ci est naturellement muni d'une flèche ``de descente'' $A\rightarrow R\underline\G(A_\infty)$, dont
on construira un cône fonctoriel: $C(j_A)\in D^b(_{\underline\G}\C)$. C'est cet objet qui encode le défaut de descente. Concrètement la cohomologie
de $C(j_A)$ est nulle en degrés $<-1$, puis vaut ensuite $(Ker\, j_n)$, $(Coker\, j_n)$, $(H^1(\G_n,A_\infty))$, $(H^2(\G_n,A_\infty))$ ... De
manière ``duale'', on définit le défaut de codescente $C(k_A)\in D^b({_{\underline\G}\C})$ (ici il faut prendre certaines précautions, cf texte). On
procède alors approximativement comme suit:

\noindent - Si $A_n$ est de la forme $R\G(\G_n,A_\infty)$ (ie. si $C(j_A)=0$), on construit un isomorphisme $\alpha_j$ (théorème \ref{dualite}).

\noindent - Dans le cas cas général, la définition de $C(j_A)$ permet de se ramener au cas précédent pour contruire un triangle distingué fonctoriel
$$\alpha_j(A)=\left(X_\infty\mathop\rightarrow\limits^{\alpha_j} RHom_\L(X_\infty^*,\L)\rightarrow \Delta_j(A)\rightarrow X_\infty[1]\right)$$ dans lequel
$\Delta_j(A)$ est relié explicitement à $C(j_A)$.

\noindent - On procède de manière analogue pour construire un triangle distingué fonctoriel $$\alpha_k(A)=\left(
X_\infty^*\mathop\rightarrow\limits^{\alpha_k} RHom_\L(X_\infty,\L)\rightarrow \Delta_k(A)\rightarrow X_\infty^*[1]\right)$$ dans lequel
$\Delta_k(A)$ est relié explicitement à $C(k_A)$.

\noindent - Sous l'hypothèse de finitude adéquate, on utilise la condition normique de compatibilité entre la structure inductive et projective pour
construire un isomorphisme de triangles: $\alpha_j(A)\simeq RHom_\L(\alpha_k(A),\L)$, fonctoriel en $A$.

En pratique, ces trois dernières étapes sont menées de front et aboutissent avec l'énoncé et la preuve du théorème principal (\ref{theoprincipal}).
Dans le texte, tout ce qui précède est énoncé dans le cadre plus général où $\G$ est profini de dimension cohomologie finie, tel que $\L$ soit
noethérien. Lorsque $\G$ est un groupe de Lie $p$-adique de dimension cohomologique finie, alors la dualité de Poincaré permet d'exprimer de manière
plus éclairante le rapport entre descente et codescente (cf \ref{poincare}).

Dans ce qui précède, on a négligé les questions de finitude. On y apporte une certaine attention dans le texte, notamment pour établir des critères
de finitude utiles en pratique. \vspace{.2cm}

A ce point, s'impose la question du contrôle de $\Delta_j(A)$. Ne sachant comment aborder le cas général, on s'intéresse au cas particulier où
$A_n=M_{\G_n}$ pour $M$ un $\L$-module de type fini et $\G\simeq \Zp^d$, $d\ge 1$. Dans cette situation, on s'attache notamment à établir un critère
maniable pour la pseudo-nullité des modules de cohomologie de $\Delta_j(A)$, ce qui sera utile dans les applications. Il serait sans doute utile de
pousser plus loin l'étude embryonnaire commencée ici.

Dans le cas particulier où $\G\simeq\Zp$, apparaît un phénomène nouveau: lorsque $A$ est un système normique raisonnable, les systèmes inductifs
(resp. projectif) sous-jacents à la cohomologie de $C(k_A)$ (resp. $C(j_A)$) ont tendance à se stabiliser. Lorsque cela se produit (notamment pour
les systèmes normiques usuels provenant de l'arithmétique) il est intéressant de noter que l'on peut reconstruire le système normique $(A_n)$
asymptotiquement à partir de la seule donnée de la flèche $\alpha_j:X_\infty\rightarrow RHom_\L(RHom_\Zp(A_\infty,\Zp),\L)$. En pratique, le
phénomène de stabilisation - lorsqu'il se produit - permet surtout de dresser un tableau général des relations entre descente, codescente et strcture
de $X_\infty$ et $A_\infty$. A cet égard, la proposition \ref{exhaust} devrait être considérée comme un formulaire. \vspace{.2cm}

Pour illustrer nos méthodes, on donne deux applications: \vspace{0.2cm}

\noindent 1.  à l'étude du groupe des $(p)$-classes - noté $\C l'$ dans une $\Zp^d$-extension. Outre une preuve unifiée de nombreux résultats connus,
on obtient notamment (cf \ref{propcl}):

- l'existence d'un morphisme $\limproj \C l'_n \rightarrow Hom_\Zp(\limind \C l'_n,\Qp/\Zp)$ dont le noyau et le conoyau sont pseudo-isomorphes. Ce
résultat fait suite aux tentatives précédentes de \cite{MC}, \cite{LN}, \cite{NV} (cf. \ref{remrem}), et complète \cite{Ne} 9.4.1. Sa preuve dépend
de l'étude de $\Delta_j(A)$ mentionnée plus haut ainsi que d'un résultat de \cite{Ne}.

- une nouvelle description du sous-module pseudo-nul maximal de $\limproj \C l'_n$, sous $(Dec_2)$ (cf texte). L'intérêt de ce résultat tient à sa
relation à la conjecture de Greenberg généralisée (cf \cite{V2} ou \cite{NV}). \vspace{0.2cm}

\noindent 2.  à la théorie d'Iwasawa cyclotomique d'un corps de nombres abélien. Dans ce cadre on retrouve \ref{fait1} et \ref{fait2} comme cas
particuliers de la situation générale \ref{exhaust}. Cette approche permet notamment de répondre à une question de \cite{NL}. Cette partie est
rédigée de manière à ce que la lecture des énoncés soit indépendante du reste du texte.  \vspace{.2cm}

L'un des objectifs recherchés dans ces deux applications est une meilleure compréhension de la nature (arithmétique ou algébrique) des résultats
connus, ainsi que les hypothèses minimales nécessaires à leur validité. Pour cette raison, mais aussi pour faire apparaître plus clairement la ligne
conductrice de notre étude, on s'abstiendra d'invoquer les résultats préexistants dont la preuve relève des méthodes de ce travail.

En appendice, on donne une construction directe du complexe dualisant, laquelle permet notamment son utilisation dans le contexte des systèmes
normiques.

Cet article puise visiblement sa motivation et son inspiration dans les méthodes développées par T. Nguyen Quand Do. Ce travail lui est dédié. Le
paragraphe 6.2, reprend et améliore les résultats d'une version antérieure, grâce aux suggestions de J. \nekovar. En particulier, le théorème
\ref{propcl} et ses corollaires lui sont partiellement dûs, et je le remercie de m'avoir permis de reproduire ici ses arguments. Je remercie aussi B.
Angles, J-R. Belliard et R. Sharifi pour leurs encouragements à rédiger ce travail, ainsi que F. Nuccio pour quelques discussions intéressantes.
Enfin, je remercie le rapporteur anonyme pour les nombreuses corrections qu'il m'a suggérées.

\section{Préliminaires}

La présente étude requiert un cadre dans lequel on puisse à la fois utiliser le formalisme des catégories dérivées, et effectuer facilement les
passages à la limite habituels de la théorie d'Iwasawa. A cet effet, on introduit la catégorie des ``systèmes normiques'' le long de la tour
d'Iwasawa. On fixe d'abord quelques notations, puis on explique rapidement comment fonctionnent l'algèbre homologique et les passages à la limite
dans cette catégorie.

\subsection{$\L$-modules, systèmes normiques} \label{not1}

Fixons un groupe profini $\G$, et indexons la famille des sous-groupes ouverts distingués de $\G$ par un ensemble $N$. Pour $n\in N$, on note $\G_n$
le sous-groupe de $\G$ correspondant, et $G_n:=\G/\G_n$. L'anti-inclusion des sous-groupes confère à $N$ une structure d'ensemble ordonné filtrant
($n\le m\Leftrightarrow \G_n\supset \G_m$), suivant laquelle $\G=\limproj G_n$.

On fixe une fois pour toutes un nombre premier $p$, et $\L:=\limproj \Zp[G_n]$ désigne l'algèbre d'Iwasawa de $\G$. Notons ${_\L\C}$ (resp. $_\Zp\C$)
la catégorie des $\L$-modules à gauche (resp. $\Zp$-modules). On considère aussi ${_\G\C}$ (resp. $_{G_n}\C$) la sous-catégorie pleine de ${_\L\C}$,
formée par les modules sur lesquels $\G$ agit discrètement (resp. $\G_n$ agit trivialement), ie. $M\in {{_\G\C}}\Leftrightarrow M=\cup M^{\G_n}$. On
adopte des notations analogues pour les catégories de modules à droite correspondantes. Si $*_1$ et $*_2$ sont deux listes de symboles choisis parmi
$\L$, $\Zp$, $\G$, $G_n$, on note ${_{*_1}\C_{*_2}}$ la catégorie de multi-modules correspondante (les différentes structures doivent être
compatibles entre elles). Ainsi, ${{_\G\C}_{\L}}$ désigne par exemple la catégorie des $\L$-bimodules sur lesquels l'action de $\G$ à gauche est
discrète.

\begin{defn}\label{defIw} $ $ \\
\noindent Soit ${_{*_1}\C_{*_2}}$ une catégorie de multi-modules comme ci-dessus. La catégorie des ``systèmes normiques'', notée
${_{\underline\G,*_1}\C_{*_2}}$ est définie comme suit:

- un objet $A$ de ${_{\underline\G,*_1}\C_{*_2}}$ consiste en la donnée, pour tout couple $n\le m$, de flèches de ${_{\G,*_1}\C_{*_2}}$: $j_{n,m}:
A_n\rightarrow A_m$ et $k_{m,n}:A_m\rightarrow A_n$ vérifiant:

\noindent $(Norm1)$: $A_n\in {_{G_n,*_1}\C_{*_2}}$.

\noindent $(Norm2)$: Si $n_1\le n_2\le n_3$, $j_{n_2,n_3}\circ j_{n_1,n_2}=j_{n_1,n_3}$ et $k_{n_2,n_1}\circ k_{n_3,n_2}=k_{n_3,n_1}$.

\noindent $(Norm3)$: Pour $n\le m$, $j_{n,m}\circ k_{m,n}:A_m\rightarrow A_m$ est la multiplication à gauche par l'élément $\sum_{g\in \G_n/\G_m}g\in
\Zp[G_m]$.

\noindent $(Norm4)$: Pour $n\le m$, $k_{m,n}\circ j_{n,m}:A_n\rightarrow A_n$ est la multiplication par $(\G_n:\G_m)$.

- une flèche $A\rightarrow B$ est la donnée pour chaque $n$ d'une flèche de ${_{*_1}\C_{*_2}}$: $A_n\rightarrow B_n$, qui commute avec les $j_{n,m}$
et les $k_{m,n}$.
\\ On définit une catégorie ${_{*_1}\C_{*_2,\underline \G}}$ de façon analogue.
\end{defn}

L'exemple fondamental, noté $\underline\L$, est défini comme suit:  $j_{n,m}:\Zp[G_n]\rightarrow \Zp[G_m]$, $x\mapsto
\sum_{g\in\G_n/\G_m}gx=x\sum_{g\in\G_n/\G_m}g$ ($\G_n/\G_m$ est distingué dans $G_m$) et $k_{m,n}:\Zp[G_m]\rightarrow\Zp[G_n]$ la projection
naturelle. $\underline\L$ est au choix un objet de $_{\underline{\G}}\C_\L$, ${_\L\C}_{\underline{\G}}$ (ou encore, par oubli de structure, de
${_{\underline\G}\C}$ ou $\C_{\underline\G}$).

Lorsque les flèches de transition d'un objet $A\in {_{\underline{\G}}\C}$ sont évidentes, on désignera simplement l'objet $A\in
{_{\underline{\G}}\C}$ par la collection $(A_n)$ (eg. $\underline\L=(\Zp[G_n])$). On remarque que ${_{\underline\G,*_1}\C_{*_2}}$ est une catégorie
abélienne, munie pour chaque $n$ d'un foncteur de projection exact: $\pi_n:{_{\underline\G,*_1}\C_{*_2}}\rightarrow _{G_n,*_1}\C_{*_2}$, $A\mapsto
A_n$.

\begin{defn}\label{homtens} $(i)$ Si $(A,B)\in {_{\L,*_1}\C_{*_2}}\times {_{\L,*_3}\C_{*_4}}$ ou ${_{*_1}\C_{*_2,\L}}\times _{*_3}\C_{*_4,\L}$,
on note $\underline{Hom}(A,B)$ le système normique formé de la collection des $G_n$-modules $Hom_{\Zp[[\G_n]]}(A,B)$ (action par conjugaison), munie
des applications de transition évidentes. C'est un objet de ${_{\underline\G,*_2,*_3}\C_{*_1,*_4}}$ ou $_{*_2,*_3}\C_{*_1,*_4,\underline\G}$, suivant
qu'on considère l'action de $G_n$ par conjugaison à gauche ou à droite.

$(i)$ Si $(A,B)\in{_{*_1}\C_{*_2,\L}}\times {_{\L,*_3}\C_{*_4}}$, on note $A\underline\otimes B$ le système normique formé de la collection des
$G_n$-modules $A\otimes_{\Zp[[\G_n]]} B$ (action par conjugaison), munie de ses applications de transition naturelles. C'est un objet de
$_{\underline\G,*_1,*_3}\C_{*_2,*_4}$ ou $_{*_1,*_3}\C_{*_2,*_4,\underline\G}$, suivant qu'on considère l'action de $G_n$ par conjugaison à gauche ou
à droite.
\end{defn}

Si $A\in {_{\underline\G,*_1}\C_{*_2}}$, on note simplement $\limind A\in {_{\G,*_1}\C_{*_2}}$ (resp. $\limproj A\in {_{\L,*_1}\C_{*_2}}$) la limite
inductive (resp. projective) prise suivant les $j_{n,m}$ (resp. $k_{n,m}$). De même pour les modules à droite.

L'axiome $(Norm1)$ donne lieu pour chaque $n\in N$ à une flèche de descente $A_n\rightarrow (\limind A)^{\G_n}$, ainsi qu'à une flèche de codescente
$(\limproj A)_{\G_n}\rightarrow A_n$. L'axiome $(Norm3)$ assure que ces deux collections de flèches commutent aux applications de transition
$j_{n,m}$ et $k_{m,n}$. On a donc, dans ${_{\underline\G,*_1}\C_{*_2}}$ une flèche de descente $A\rightarrow \underline{Hom}(\Zp,\limind A)$, ainsi
qu'une flèche de codescente $\Zp\underline\otimes\limproj A\rightarrow A$. Pour des raisons d'algèbre homologique évidentes, on préfèrera noter
$\underline\G(\limind A)$ au lieu de $\underline{Hom}(\Zp,\limind A)$ (voir prochain paragraphe). On notera que le noyau et le conoyau de
l'application de descente $A\rightarrow \underline\G(\limind A)$ sont de $\Zp$-torsion (au sens où chaque composante est de $\Zp$-torsion), grâce à
l'axiome $(Norm4)$. C'est la seule utilité de cette axiome, dont la totalité du présent texte est en fait indépendante.

D'après la discussion précédente, le foncteur $\Zp\underline\otimes-:{_{\L,*_1}\C_{*_2}}\rightarrow {_{\underline\G,*_1}\C_{*_2}}$ (resp.
$\underline\G:{_{\G,*_1}\C_{*_2}}\rightarrow {_{\underline\G,*_1}\C_{*_2}}$) est adjoint à gauche au foncteur
$\limproj:{_{\underline\G,*_1}\C_{*_2}}\rightarrow {_{\L,*_1}\C_{*_2}}$ (resp.  adjoint à droite au foncteur
$\limind:{_{\underline\G,*_1}\C_{*_2}}\rightarrow {_{\G,*_1}\C_{*_2}}$). En particulier, $\underline\G(-)$ conserve les injectifs (car $\limind$ est
exact à gauche). Pour $\Zp\underline\otimes -$ et les projectifs, voir \ref{jensen}.

Un foncteur additif $F:{_{\G,*_1}\C_{*_2}}\rightarrow {_{\G,*_3}\C_{*_4}}$ est dit $\G$-linéaire s'il vérifie $F(_gm)={_gm}$ pour tout $g\in \G$,
$_gm$ désignant la multiplication à gauche par $g$. Un tel foncteur possède alors une extension évidente, $F:{_{\underline\G,*_1}\C_{*_2}}\rightarrow
{_{\underline\G,*_3}\C_{*_4}}$ définie par $F(j_{n,m})=j_{n,m}$ et $F(k_{m,n})=k_{m,n}$. De même, un foncteur additif contravariant $\G$-linéaire
$F:{_{\G,*_1}\C_{*_2}}\rightarrow {_{*_3}\C_{*_4,\G}}$ possède une extension $F:{_{\underline\G,*_1}\C_{*_2}}\rightarrow
{_{*_3}\C_{*_4,\underline\G}}$, définie par $F(j_{n,m})=k_{m,n}$ et $F(k_{m,n})=j_{n,m}$. Il y a bien sûr une discussion similaire pour les foncteurs
covariants ${_{*_1}\C_{*_2,\G}}\rightarrow {_{*_3}\C_{*_4,\G}}$ et les foncteurs contravariants ${_{*_1}\C_{*_2,\G}}\rightarrow {_{\G,*_3}\C_{*_4}}$.

De cette manière, on obtient par exemple un foncteur $Hom_\Zp(-,\Zp):{_{\underline\G}\C_\L}\rightarrow {_\L\C_{\underline\G}}$. On relève
l'isomorphisme important (cas particulier de \ref{indvscoind}):
$$Hom_\Zp(\underline\L,\Zp)\simeq\underline\L$$ où le premier $\underline\L$ est un objet de ${_{\underline\G}\C_\L}$ et le second un objet de
${_\L\C}_{\underline\G}$. \vspace{.2cm}

\noindent \emph{\textbf{Induction normique.}} \\ \label{induction}

On définit quatre foncteurs (induction et co-induction):
$$\underline{Ind},\underline{Coind}: {_\L\C}\rightarrow {_\L\C}_{\underline\G} \hspace{1cm}
\underline{Ind},\underline{Coind}: \C_\L\rightarrow {_{\underline \G}\C_\L}$$ de la manière suivante (cf def. \ref{homtens})

\noindent - $\underline{Ind}A:=\L\underline\otimes A$ (resp. $A\underline\otimes \L$), muni de la conjugaison à droite (resp. gauche) si
$A\in{_\L\C}$ (resp. $A\in \C_\L$). C'est un $\L$-module à gauche (resp. à droite) via le facteur $\L$.

\noindent - $\underline{Coind}A:=\underline{Hom}(\L,A)$, muni de la conjugaison à droite (resp. gauche) si $A\in{_\L\C}$ (resp. $A\in \C_\L$). C'est
un $\L$-module à gauche (resp. droite) via la structure à droite (resp. gauche) de $\L$.

Ont alors lieu les isomorphismes bifonctoriels suivants:

\noindent - Si $(A,B)\in {_\L\C}\times{_\L\C}$, dans ${_{\underline\G}\C}$ (resp. $\C_{\underline\G}$):
$$\underline{Hom}(A,B)\simeq Hom_\L(\underline{Ind}A,B)    \hspace{1cm}
\hbox{(resp.}\,\,   \underline{Hom}(A,B)\simeq Hom_\L(A,\underline{Coind}B)\hbox{)}$$

\noindent - Si $(A,B)\in \C_\L\times\C_\L$, le même isomorphisme  a lieu, mais cette fois dans $\C_{\underline\G}$ (resp. ${_{\underline\G}\C}$).

\noindent - Si $(A,B)\in \C_\L\times{_\L\C}$, on a dans ${_{\underline\G}\C}$ (resp. $\C_{\underline\G}$):
$$\underline{Ind}(A)\otimes_\L B\simeq A\underline\otimes B\hspace{1cm}
\hbox{(resp.}\,\, A\otimes_\L \underline{Ind}B\simeq A\underline\otimes B\hbox{)}$$

Regarder $\Zp$ comme un $\L$-module à gauche provoque dans ${_\L\C}_{\underline\G}$ les isomorphismes suivants
$$\underline{Ind}(\Zp)\simeq \underline \L\simeq Hom_\Zp(\underline\L,\Zp)\simeq\underline{Coind}(\Zp)$$

Plus généralement:

\begin{prop}\label{indvscoind} Les foncteurs $\underline{Ind},\underline{Coind}:{_\L\C}\rightarrow{_\L\C_{\underline\G}}$ sont naturellement
isomorphes. De même pour les foncteurs  $\C_\L\rightarrow {_{\underline \G}\C_\L}$. \end{prop} Preuve: On laisse au lecteur le soin de vérifier que
la collection d'applications naturelles $\phi_n:Hom_{\Zp[[\G_n]]}(\L,A)\rightarrow \L\otimes_{\Zp[[\G_n]]} B$, $f\mapsto \sum_{g\in
\G_n\backslash\G}g^{-1}\otimes f(g)$ respecte la structure de système normique à droite, ainsi que la structure de $\L$-module à gauche, et que c'est
un isomorphisme.
\begin{flushright}$\square$\end{flushright}

\noindent \emph{\textbf{Induction discrète.}} \, (inutile, sauf pour \ref{compdu}) \\ \label{inductiondisc}

On voit facilement que $\underline{Ind}:{_{\L}\C}\rightarrow {_\L\C_{\underline\G}}$ transforme objets de ${_\G\C}$ en objets de
${_\G\C_{\underline\G}}$. De même pour les trois variantes. A partir de là, on peut définir des foncteurs d'induction discrète à partir des
précédents, par limite inductive:
$$Ind_\G,Coind_\G: {_\G\C}\rightarrow {_\G\C}_{\G} \hspace{1cm}
Ind_\G,Coind_\G: \C_\G\rightarrow {_{\G}\C_\G}$$ $Ind_\G(A):=\limind
\underline{Ind}(A)$ et $Coind_\G(A):=\limind\underline{Coind}(A)$.

Un processus de double limite permet d'écrire des isomorphismes bifonctoriels à partir des précédents (se ramener au cas où $A$ est finiment
engendré). Par exemple:

- Si $(A,B)\in {_\G\C}\times {_\G\C}$, on a dans ${_\L\C}$ (resp. ${\C_\L}$): {\small$$Hom_\Zp(A,B)\simeq Hom_\L(Ind_\G(A),B)\hspace{.3cm}
\hbox{(resp.}\, \,  Hom_\Zp(A,B)\simeq Hom_\L(A,Coind_\G(A))\hbox{)}$$} La structure de $\L$-module de $Hom_\Zp(A,B)$ étant ici induite par l'action
de $\G$ par conjugaison. On prendra garde au fait que celle-ci n'est pas discrète en général. Cette structure de $\L$-module peu naturelle (et
contraire aux conventions générales de cet article, expliquées plus haut), ne sera utilisée en pratique que dans la situation où l'action de $\G$ sur
$B$ est triviale (voir notamment \ref{inddu}).

\subsection{(Bi-)foncteurs dérivés}

Soit $F$ un bifoncteur additif exact à gauche. Si l'une des deux catégories de départ possède suffisamment d'objets $F$-acycliques (ou déployants),
on peut définir le foncteur dérivé droit $RF$. Nous appliquons ce principe pour dériver certains des foncteurs définis au paragraphe précédent. On
indique ensuite comment s'écrivent quelques unes des compatibilités habituelles entre $Hom$ et $\otimes$ dans ce cadre. Pour une catégorie abélienne
$\C$, $Kom^a(\C)$, $a\ge 1$ (resp. $Kom^a_{naifs}(\C)$, $Kom^*(\C)$, $K^*(\C)$, $D^*(\C)$, $*=+,-,b$) désigne la catégorie des complexes $a$-uples
(resp. celle des complexes $a$-uples naifs, celle des complexes simples, des complexes simples à homotopie près, la catégorie dérivée). Un objet de
$Kom(\C)=Kom^1(\C)$ sera le plus souvent symbolisé par $A$, et par $A^\bullet$ seulement lorsque la clarté l'exige.

Comme notre étude est restreinte au cas où $\G$ est de dimension cohomologique finie et qu'en pratique ou travaillera toujours dans $D^b$, on ne
cherche pas à étendre au maximum le domaine de définition des foncteurs dérivés.

\begin{prop} \label{bifder}
1. Soient $*_1,*_2,*_3,*_4$ quatre listes de symboles choisis parmi $\Zp$ et $\L$. Si $*_1=*_2=\emptyset$ ou $*_3=*_4=\emptyset$, alors:

$(i)$ Le bifoncteur $\underline{Hom}:{_{\L,*_1}\C_{*_2}}\times{_{\L,*_3}\C_{*_4}}\rightarrow {_{\underline\G,*_2,*_3}\C_{*_1,*_4}}$ est dérivable à
droite, ainsi que ses trois variantes obtenues en changeant $\L$ et/ou ${\underline\G}$ de côté. On note $R\underline{Hom}:D^-\times D^+\rightarrow
D^+$ le foncteur ainsi obtenu.

$(ii)$ Le bifoncteur $\underline\otimes:{_{*_1}\C_{*_2,\L}}\times{_{\L,*_3}\C_{*_4}}\rightarrow {_{\underline\G,*_1,*_2}\C_{*_3,*_4}}$ (ou
${_{*_1,*_2}\C_{*_3,*_4,\underline\G}}$) est dérivable à gauche.  On note $\underline\Ltens:D^-\times D^-\rightarrow D^-$.

\noindent 2. Soit $*=\L$ ou $\Zp$, et $*_1,*_2,*_3,*_4$ quatre listes de symboles choisis parmi $\Zp,\L,\underline\G$ de façon à ce que
$\underline\G$ apparaisse au plus dans une seule liste. Si $*_1=*_2=\emptyset$ ou $*_3=*_4=\emptyset$, alors:

$(i)$ Le bifoncteur $Hom_*:{_{*,*_1}\C_{*_2}}\times {_{*,*_3}\C_{*_4}}\rightarrow {_{*_2,*_3}\C_{*_1,*_4}}$ est dérivable à droite. On note
$RHom_*:D^-\times D^+\rightarrow D^+$. Idem si l'on change $*$ de coté.

$(ii)$ Le bifoncteur $\otimes_*:{_{*_1}\C_{*_2,*}}\times {_{*,*_3}\C_{*_4}}\rightarrow {_{*_1,*_3}\C_{*_2,*_4}}$ est dérivable à gauche. On note
$\Ltens_*:D^-\times D^-\rightarrow D^-$. \end{prop}

Preuve: 1. $(i)$ Supposons $*_3=*_4=\emptyset$. La catégorie ${_{\L,*_3}\C_{*_4}}={_\L\C}$ possède alors suffisamment d'injectifs. Comme la catégorie
des injectifs de ${_\L\C}$ vérifie les deux conditions de \cite{KS} 13.4.5, cela permet de conclure. Bien sûr, la flèche naturelle de
$D^+({_{\underline\G,*_2}\C_{*_1}})$: $\underline{Hom}(A,I)\rightarrow R\underline{Hom}(A,I)$ est un isomorphisme dès que $I$ est injectif. Un objet
$I\in Kom(_\L\C)$ possédant cette propriété pour tout $A\in Kom(_{\L,*_1}\C_{*_2})$ sera dit déployant (et les couples $(A,I)$ déployés) pour le
foncteur $R\underline{Hom}$. On emploie une terminologie analogue pour tous les bifoncteurs dérivés.

Les autres points se traitent de façon analogue.
\begin{flushright}$\square$\end{flushright}

\begin{rem} \label{rempin} Les foncteurs de projection $\pi_n:{_{\underline\G,*_2,*_3}\C_{*_1,*_4}}\rightarrow _{G_n,*_2,*_3}\C_{*_1,*_4}$
étant exacts, ils passent aux catégories dérivées. En composant $\pi_n$ avec $R\underline{Hom}$ défini ci-dessus, on retrouve les bifoncteurs usuels
$RHom_{\Zp[[\G_n]]}$. De même $\pi_n\circ\underline\Ltens=\Ltens_{\Zp[[\G_n]]}$. Il est utile de noter que la famille $(\pi_n)$ est conservative.
\end{rem}

Les compatibilités habituelles entre $Hom$ et $\otimes$ s'étendent naturellement à ce cadre (adjonction, évaluation). Détaillons quelques cas utiles.

\begin{prop} \label{comphomtens} Soit $*=\L$ ou $\Zp$.

\noindent 1.  Pour $(A,B,C)\in D^-(\C_{\L})\times D^-(_{\L}\C_{*})\times D^+(\C_{*})$, il y a dans $D^+({_{\underline\G}\C})$ et
$D^+(\C_{\underline\G})$ un isomorphisme d'adjonction tri-fonctoriel:
$$adj:RHom_*(A\underline\Ltens B,C)\simeq R\underline{Hom}(A,RHom_*(B,C))$$

\noindent 2. $(i)$ Pour $(A,B,C)\in D^b(\C_\L)\times D^b(_*\C_\L)\times D^b(_*\C)$, il y a dans $D^b({_{\underline\G}\C})$ et
$D^b(\C_{\underline\G})$ un morphisme d'évaluation trifonctoriel:
$$ev:A\underline\Ltens RHom_*(B,C)\rightarrow RHom_*(R\underline{Hom}(A,B),C)$$

$(ii)$ Pour $(A,B,C)\in D^b({_\L\C})\times D^b({_\L\C}_*)\times D^b(\C_*)$, il y a dans $D^b({_{\underline\G}\C})$ et $D^b(\C_{\underline\G})$ un
morphisme d'évaluation trifonctoriel:
$$ev:RHom_*(B,C)\underline\Ltens A\rightarrow RHom_*(R\underline{Hom}(A,B),C)$$

$(iii)$ Les morphismes $(i)$ et $(ii)$ ci-dessus sont des isomorphismes lorsque $A$ est isomorphe (dans $D^b$) à un complexe parfait.
\end{prop}
Preuve: 1. $(i)$ Pour chaque $n\in N$, il y a un isomorphisme naturel entre les trifoncteurs $\C_\L\times {_\L\C}_*\times\C_*\rightarrow _{G_n}\C$
$(A,B,C)\mapsto Hom_*(A\otimes_{\Zp[[\G_n]]} B,C)$ et $(A,B,C)\mapsto Hom_{\Zp[[\G_n]]}(A,Hom_*(B,C))$. Il est facile de voir que cet isomorphisme
est compatible à la structure de système normique, au sens où il identifie $Hom_*(k_{m,n},C)$ avec $j_{n,m}$ et $Hom_*(j_{n,m},C)$ avec $k_{m,n}$. On
obtient ainsi un isomorphisme entre les trifoncteurs $\C_\L\times {_\L\C}_*\times\C_*\rightarrow {_{\underline\G}\C}$, $(A,B,C)\mapsto
Hom_*(A\underline\otimes B,C)$ et $(A,B,C)\mapsto\underline{Hom}(A,Hom_*(B,C))$, d'où immédiatement un isomorphisme de trifoncteurs
$K^-(\C_{\L})\times K^-(_{\L}\C_{*})\times K^+(\C_*)\rightarrow K^+({_{\underline\G}\C})$ par extension des trifoncteurs aux complexes. Mais alors,
l'extension au complexes du trifoncteur $(A,B,C)\mapsto Hom_*(A\underline\otimes B,C)$ est isomorphe (avec les conventions de signes adéquates, cf
\cite{De}, 1.1.8):

- d'une part au trifoncteur obtenu en composant l'extension aux complexes du bifoncteur $(D,C)\mapsto Hom_*(D,C)$ avec l'extension aux complexes du
bifoncteur $(A,B)\mapsto A\underline\otimes B$,

- d'autre part au trifoncteur obtenu en composant l'extension aux complexes du bifoncteur $(A,D)\mapsto \underline{Hom}(A,D)$ avec l'extension aux
complexes du bifoncteur $(B,C)\mapsto Hom_*(B,C)$.

Aussi les deux derniers trifoncteurs sont-ils isomorphes. Il suffit de choisir $A$ à objets projectifs et $C$ à objets injectifs pour obtenir
l'isomorphisme de l'énoncé.

2. $(i)$ Se traite de façon analogue, à partir du morphisme fonctoriel d'évalu\-ation entre les trifoncteurs
${_\L\C}\times{_\L\C}_*\times\C_*\rightarrow {_{\underline\G}\C}$. Idem pour $(ii)$.

$(iii)$ est une conséquence directe du fait suivant: si $A\in \C_\L$ (resp. ${_\L\C}$) est projectif de type fini, alors la flèche d'évaluation
$$A\otimes_{\Zp[[\G_n]]}Hom_*(B,C)\rightarrow Hom_*(Hom_{\Zp[[\G_n]]}(A,B),C)$$ $$\hbox{(resp.}\hspace{0.2cm}
Hom_*(B,C)\otimes_{\Zp[[\G_n]]} A\rightarrow Hom_*(Hom_{\Zp[[\G_n]]}(A,B),C)\hspace{0.2cm}\hbox{)}$$ est un isomorphisme pour tout $(B,C)\in
{_*\C_\L}\times {_*\C}$ (resp. ${_\L\C_*}\times \C_*$).
\begin{flushright}$\square$\end{flushright}

Soit $\C'$ une sous-catégorie d'une catégorie $\C$. Si  $F:D^+(\C)\rightarrow D^+(\C'')$ est un foncteur, on notera encore $F:D^+(\C')\rightarrow
D^+(\C'')$ le foncteur obtenu en composant $F$ avec le foncteur évident $D^+(\C')\rightarrow D^+(\C)$, et on parlera de la ``restriction de $F$ à
$D^+(\C')$''. Il s'agit d'un abus de langage puisqu'en général le foncteur évident n'est pas fidèle.

\begin{prop} \label{compcohgal} Notons $\underline\G:{{_\G\C}}\rightarrow {_{\underline\G}\C}$ (resp. ${{\C_\G}}\rightarrow {\C_{\underline\G}}$) la restriction à
${{_\G\C}}\subset {_\L\C}$ (resp. ${\C_\G}\subset {\C_\L}$) du foncteur $\underline{Hom}(\Zp,-)$. Alors:

$(i)$ Le foncteur $\underline\G$ est dérivable à droite, et donne donc $R\underline\G:D^+({_\G\C})\rightarrow D^+({_{\underline\G}\C})$. On a de même
$R\underline\G:D^+(\C_\G)\rightarrow D^+(\C_{\underline\G})$.

$(ii)$ Il y a une flèche naturelle de $R\underline\G$ vers la restriction à $D^+({_\G\C})$ du foncteur
$R\underline{Hom}(\Zp,-):D^+({_\L\C})\rightarrow D^+({_{\underline\G}\C})$. Idem en remplaçant ${_\G\C}$ par $\C_\G$.

$(iii)$ Si $\L$ est noethérien et si les objets de cohomologie de $A\in D^+({_\G\C})$ (resp. $D^+(\C_\G)$) sont tous des $\Zp$-modules de torsion,
alors la flèche précédente induit un isomorphisme $R\underline\G(A)=R\underline{Hom}(\Zp,A)$.
\end{prop}
Preuve: $(i)$ ${_\G\C}$ et $\C_\G$ possèdent suffisamment d'injectifs.

$(ii)$ résulte de la propriété universelle du foncteur dérivé droit $R\underline\G$.

$(iii)$ Par troncature, on se ramène tout de suite au cas où $A$ est concentré en degré $0$. Dans ce cas, $A$ possède une résolution $A\rightarrow I$
dans laquelle les objets de $I$ sont de $\Zp$-torsion, injectifs dans ${_\G\C}$. Le résultat suit, d'après le lemme \ref{discdep} (appendice).
\begin{flushright}$\square$\end{flushright}

\begin{defn} \label{Bid} Soit $*=\L$ ou $\Zp$, et $*_1,*_2$ deux listes de symboles choisis parmi $\Zp,\L,\underline\G$, dont l'une au moins contient $*$.

On définit $Bid_*:D^b({_{*_1}\C_{*_2}})\rightarrow D^b({_{*_1}\C_{*_2}})$ par $$Bid_*(A):=RHom_*(RHom_*(A,*),*)$$ et l'on note $bid_*:A\rightarrow
Bid_*(A)$ le morphisme naturel de bidualité.
\end{defn}

Si $*$ est noethérien de dimension homologique finie, $bid_*$ est un isomorphisme dès $A$ est d'amplitude bornée, à objets de cohomologie de type
fini sur $*$  (dualité de Grothendieck). Pour $*=\Zp$, on décrira un peu plus loin l'effet de $Bid_\Zp$ sur une classe plus large d'objets; le cas
significatif est celui où $A$ un module divisible de cotype fini: $Bid_\Zp(A)$ est alors le module de Tate de $A$ placé en degré $-1$ (cf
\ref{bidettate}).

\subsection{Limites}

Soient $\C_1,\C_2$ (resp. $\C'$) deux catégories abéliennes dans lesquelles les produits infinis existent et sont exacts (resp. une cat. ab. avec
sommes infinies exactes). On pourra par exemple prendre $\C_1={_{*}\C}$ avec $*=\Zp,\L$ ou $\underline\G$ (resp. $\C'={_*\C}$ avec $*=\Zp,\L,\G$ ou
$\underline\G$). Soit $I$ un ensemble ordonné filtrant et considérons $\C_1^{I^\circ}$ la catégorie des systèmes projectifs de $\C_1$ indexés sur
$I$, et ${\C'}^I$ celle des systèmes inductifs de $\C'$ indexés sur $I$. On dispose alors de foncteurs $\limind:{\C'}^I\rightarrow \C'$ et
$\limproj:\C_1^{I^\circ}\rightarrow \C_1$. Le premier est exact et passe donc aux catégories dérivées; le second est seulement exact à gauche, et son
usage nécessite les résultats suivants:

\begin{prop} \label{proplimproj} Soient $I$, $\C_1,\C_2,\C'$ comme ci-dessus. Alors:

\noindent  1. $\limproj_I$ possède un dérivé droit $R\limproj_I:D^+(\C_1^{I^\circ})\rightarrow D^+(\C_1)$.

\noindent 2. Supposons $\C_1$ munie d'un foncteur $\Phi$ vers la catégorie $Ab$ des groupes abéliens. Si $\Phi$ est exact et commute aux limites
projectives filtrantes, alors $R\limproj_I$ commute aussi à $\Phi$. Si de plus $\Phi$ est fidèle, alors

$(i)$ $R\limproj_I$ est de dimension cohomologique $\le 1$ dès que $I$ possède un ensemble cofinal dénombrable.

$(ii)$ $R^q\limproj_I A$ s'annule pour $q\ge 1$, dès que $\Phi(A)\in Ab^{I^\circ}$ est dans l'image essentielle du foncteur d'oubli
$Comp^{I^\circ}\rightarrow Ab^{I^\circ}$, $Comp$ désignant la catégorie des groupes topologiques compacts.

\noindent 3. Soit $F:\C_1\times \C'\rightarrow \C_2$ un bifoncteur covariant exact à gauche dérivable à droite (resp. covariant exact à droite et
dérivable à gauche de dimension homologique finie) et notons $*=+$ (resp. $*=b$). On suppose que:

 - $\C'$ possède suffisamment d'objets déployants (cf preuve de \ref{bifder}) pour $RF$ (resp. $LF$).

 - via la première variable, $F$ commute aux produits infinis.

\noindent Alors $RF$ s'étend naturellement en un bifoncteur $D^*(\C_1^{I^\circ})\times D^*(\C')\rightarrow D^*(\C_2^{I^\circ})$ et l'on a pour
$(A,B)\in D^*(\C_1^{I^\circ})\times D^*(\C)$ un isomorphisme bifonctoriel $$RF(R\limproj_I A,B)\simeq R\limproj_I RF(A,B)\hspace{0.2cm} \hbox{(resp.
$LF(R\limproj_I A,B)\simeq R\limproj_I LF(A,B)$)}$$

\noindent 4. Soit $G:\C'\rightarrow \C_1$ un foncteur contravariant exact. On suppose que $G$ transforme sommes infinies en produits infinis. Il y a
alors un isomorphisme $G(\limind_I A)\simeq R\limproj_I G(A)$, fonctoriel en $A\in D^-({\C'}^I)$.

\end{prop}
Preuve: 1. résulte de l'existence dans $\C_1^{I^\circ}$ de suffisamment de systèmes projectifs flasques au sens de \cite{Ro}. Comme nous en aurons
l'utilité, indiquons brièvement de quoi il s'agit. On dit qu'un système projectif $A=(A_i)\in \C_1^{I^\circ}$ est flasque si pour pour tout
$I''\subset I'\subset I$, $\limproj_{I'}A\twoheadrightarrow \limproj_{I''}A$. Pour $A\in\C_1^{I^\circ}$, on définit un ``effacement flasque''
$A\rightarrow Fl(A)$, où $Fl(A)_i=\prod_{j\le i}A_i$ est le système projectif évident et où les flèches $A_i\rightarrow Fl(A)_i$ sont induites par
celles de $A$. Il est utile de noter que l'on dispose ici d'un foncteur résolvant ``à la Godement'': si
$Fl^\bullet:Kom^a_{naifs}(\C_1^{I^{\circ}})\rightarrow Kom_{naifs}^{a+1}(\C_1^{I^{\circ}})$ est induit par $Fl$, on note
$Fl^*:Kom(\C_1^{I^{\circ}})\rightarrow Kom(\C_1^{I^{\circ}}), A^\bullet\mapsto Tot_\Pi Fl^\bullet(A^\bullet)$, où $Tot_\Pi$ est le foncteur
``complexe total'' avec des produits (et non des sommes). On a donc $R\limproj A^\bullet=\limproj Fl^*(A^\bullet)$.

2. La construction ci-dessus montre que $R\limproj_I$ commute bien à $\Phi$.

$(i)$ On peut toujours supposer $\C_1=Ab$, et invoquer \cite{Ro}, corollaire à la proposition 5. D'après \cite{Mi}, il faut prendre garde à la
proposition 5 elle-même (reportée aussi dans \cite{J1}), qui est fausse en toute généralité. Lorsque $\C_1=Ab$, la version usuelle du critère de
Mittag-Leffler remplace \cite{Ro} prop. 5 et permet cependant de conclure.

$(ii)$ Il est bien connu que $\limproj_I:Comp^{I^\circ}\rightarrow Comp$ est exact. Le résultat s'en déduit aisément. On renvoie à \cite{Je} pour une
étude plus complète des foncteurs $R^q\limproj_I$.

3. Traitons le cas où $F$ est exact à gauche dérivable à droite.  Si $B^\bullet\in Kom^*(\C)$ est un complexe à objets déployants pour $RF$ et
$A^\bullet\in Kom^*(\C_1^{I^\circ})$, $RF(A^\bullet,B^\bullet):=F(A^\bullet,B^\bullet)\in D^+(\C_2^{I^\circ})$ définit l'extension de $RF$ (resp.
$LF$) souhaitée. Les foncteurs $F(-,B^q)$ ($q\ge 0$) sont exacts, et commutent aux produits infinis. Il y a donc un isomorphisme naturel de complexes
triples naifs
$$F(\limproj Fl^\bullet(A^\bullet),B^\bullet)\simeq \limproj Fl^\bullet(F(A^\bullet,B^\bullet))$$ fonctoriel en $(A^\bullet,B^\bullet)$.
En groupant de deux manières différentes, on obtient un isomorphisme de complexes simples, duquel se déduit celui de l'énoncé.

Le cas $F$ exact à droite dérivable à gauche se traite de façon analogue. Expliquons la condition de dimension homologique finie: l'extension de $F$
au complexes fait intervenir des sommes directes, alors que ci-dessus, le passage aux complexes simples fait intervenir des produits. La condition de
dimension homologique finie permet de faire uniquement des produits finis.

4. Dualisons la contruction du 1. $(i)$: pour $A\in \C^I$, posons $Cofl(A)\rightarrow A$, où $Cofl(A)_i=\oplus_{j\le i}A_i$ avec les flèches
naturelles. De même, on définit $Cofl^{\bullet}: Kom^a_{naifs}(\C')\rightarrow Kom^{a+1}_{naifs}(\C')$ et $Cofl^*(A^\bullet):=Tot_\oplus
Cofl^\bullet(A^\bullet)$. Comme $G$ est exact et transforme sommes en produits, on a pour $A^\bullet\in Kom^-(\C')$ un isomorphisme de doubles
complexes naifs
$$G(\limind Cofl^\bullet(A^\bullet))\simeq \limproj Fl^\bullet(G(A^\bullet))$$
dont se déduit celui de l'énoncé en prenant les complexes totaux (noter que $Tot_\Pi\circ G\simeq G\circ Tot_\oplus$).
\begin{flushright}$\square$\end{flushright}

\begin{rem} \label{remlimproj} Dans le point 3., on peut jouer sur la variance en les deux variables: soit en remplaçant $\C'$ par sa catégorie opposée,
soit en remplaçant $\C_1$ par sa catégorie opposée (dans ce cas, on utilise le point 4.). Ainsi, soit par exemple $F:\C_1'\times \C'\rightarrow \C_2$
contravariant en la première variable et covariant en la seconde. On suppose que $F$ est dérivable à droite, transforme sommes infinies en la
première variable en produits infinis et que $\C'$ possède suffisamment d'objets déployants. Soit $\C_1$ la catégorie opposée à $\C_1'$,
$F':\C_1\times\C'\rightarrow \C_2$ le bifoncteur covariant déduit de $F$ et $G:\C_1\rightarrow \C_1'$ le foncteur contravariant exact induit par
l'identité.  Combinant 2. et 3. on obtient:
$$RF(\limind_I A,B)=RF'(G(\limind_I A),B)\simeq RF'(R\limproj_I G(A),B)$$ $$\hspace{3cm} \simeq R\limproj_I RF'(G(A),B)=R\limproj_I RF(A,B)$$ \end{rem}

Dans ce travail, les passages à la limites apparaissent dans deux contextes différents:

- Le plus souvent, pour former les modules d'Iwasawa: $I=N$, la catégorie des sous-groupes ouverts de $\G$ et $\C_1={_{\L}\C}$. Par abus de langage,
on parlera de la limite d'un système normique; aussi $R\limproj$ désignera-t-il souvent le foncteur composé $D^+({_{\underline\G}\C})\rightarrow
D^+({_{\L}\C}^I)\mathop\rightarrow\limits^{R\limproj}D^+(_{\L}\C)$.

On dit d'un système normique $A\in {_{\underline\G}\C}$ qu'il est de type fini sur $\Zp$ si chaque $\pi_n(A)=A_n$ l'est. On note
${_{\underline\G}\C}_{tf}\subset {_{\underline\G}\C}$ la sous-catégorie des systèmes normiques de type fini, et $Kom^*({_{\underline\G}\C})_{tf}$,
(resp. $D^*({_{\underline\G}\C})_{tf}$) la sous-catégorie pleine de $Kom^*({_{\underline\G}\C})$ (resp.  $D^*({_{\underline\G}\C})$) des complexes
dont les objets de cohomologie sont dans dans ${_{\underline\G}\C}_{tf}$. Pour les catégories de modules l'indice $tf$ prend sa signification
habituelle. Le fait suivant se déduit de \ref{proplimproj} 2. $(ii)$:

\begin{prop} \label{jensen} Si $A\in Kom^+({_{\underline\G}\C})_{tf}$, alors $\limproj A\rightarrow R\limproj A$
est un isomorphisme (dans $D^+({_\L\C})$). \end{prop}
\begin{flushright}$\square$\end{flushright}

Comme la restriction de $\limproj:{_{\underline\G}\C}\rightarrow {_\L\C}$ à ${_{\underline\G}\C}_{tf}$ est exacte, on voit que $\Zp\underline\otimes
-:{_\L\C}\rightarrow {_{\underline\G}\C}$ (qui est l'adjoint à gauche de $\limproj$) vérifie la propriété suivante: Soit $P$ un objet projectif de
${_\L\C}$, si $\Zp\underline\otimes P$ est de type fini alors c'est un objet projectif dans ${_{\underline\G}\C}_{tf}$.

- Occasionnellement, pour former les modules de Tate, ou la cohomologie continue: $I=\mathbb N$, la catégorie des entiers naturels ordonnés, et
$\C_1={{_{\underline\G}\C}}$. On conviendra alors d'indexer les systèmes projectifs par la lettre $k$, et on notera $R\limproj_k$ au lieu de
$R\limproj_I$. Dans ce cadre, indiquons le rapport entre module de Tate et bidualité:

\begin{lem} \label{bidettate} Soit $A\in D^b({_{\underline\G}\C})$, $A_n:=\pi_nA$. Alors:

$(i)$ $A\Ltens_\Zp \Z/p^k\in D^b({_{\underline\G}\C})_{tf}\, \Leftrightarrow Bid_\Zp(A)\in D^b({_{\underline\G}\C})_{tf}$. Cela se produit si et
seulement si chaque $\Zp$-module $H^q(A_n)$ possède une suite de composition dont chaque quotient est: soit de torsion et de $\Zp$-cotype fini, soit
uniquement $p$-divisible, soit de $\Zp$-type fini.

$(ii)$ Il y a un isomorphisme et une flèche naturelle $$RHom_\Zp(\Qp/\Zp,A)[1]\simeq R\limproj_k A\Ltens_\Zp\Z/p^k\rightarrow Bid_\Zp(A)$$ c'est
isomorphisme si la condition $(i)$ est satisfaite.

$(iii)$ Supposons $(i)$ satisfaite. Si de plus chaque $\Zp$-module $H^q(A_n)$ est de torsion, alors il y a un isomorphisme naturel $A\simeq
Bid_\Zp(A)\Ltens_\Zp\Qp/\Zp[-1]$
\end{lem}
Preuve: $(i)$ Commençons par montrer trois faits:

$(a)$ $Bid_\Zp(A\Ltens_\Zp \Z/p)\simeq Bid_\Zp(A)\Ltens_\Zp \Z/p$. \\ En effet, si $I$ (resp. $P$) désigne le complexe de $\Zp$-modules
$[\Qp\rightarrow \Qp/\Zp]$ (resp. $[\Zp\mathop\rightarrow\limits^p\Zp]$) placé en degrés $[0,1]$ (resp. $[-1,0]$), alors
$$Hom_\Zp(Hom_\Zp(A\otimes_\Zp P,I),I)\simeq Hom_\Zp(Hom_\Zp(A,I),I)\otimes_\Zp P$$ d'où le résultat annoncé.

$(b)$ $A\Ltens_\Zp\Z/p\in D^b({_{\underline\G}\C})_{tf}\Leftrightarrow Bid_\Zp(A\Ltens_\Zp \Z/p) \in D^b({_{\underline\G}\C})_{tf}$. \\ En effet, si
$M$ est un $\Zp$-module tué par $p$, alors on a simplement $$\begin{array}{rcl}Bid_\Zp(M) & \simeq & Hom_\Zp(Hom_\Zp(M,\Qp/\Zp)[-1],\Qp/\Zp)[-1]\\
& \simeq & Hom_\Zp(Hom_\Zp(M,\Qp/\Zp),\Qp/\Zp)\\ & \simeq & Hom_\Zp(Hom_\Zp(M,\Z/p),\Z/p)\end{array}$$ et l'on voit que $M$ est de type fini si et
seulement si $Bid_\Zp(M)$ l'est; l'assertion annoncée s'en déduit par troncature de $A\Ltens_\Zp \Z/p$.

$(c)$ Si $M\in {_{\Zp}\C}$, alors $M\Ltens_\Zp\Z/p\in D^b({_{\Zp}\C})_{tf}$ $\Leftrightarrow$ $M$ possède une suite de composition comme dans
l'énoncé. \\ L'implication $\Leftarrow$ est évidente. Montrons $\Rightarrow$. Soient $M_1\subset M_2\subset M$ définis de la manière suivante: $M_1$
est le sous-module de torsion de $M$, et l'image de $M_2$ dans $M/M_1$ est le sous-groupe constitué des éléments infiniment $p$-divisibles (noter que
$M_2/M_1$ est $p$-divisible, car $M/M_1$ est sans torsion). Alors $M_1,M_2/M_1,M/M_2$ est la suite de composition annoncée. En effet:
$Tor_1^\Zp(M_1,\Z/p)=Tor_1^\Zp(M,\Z/p)$ est fini, donc $M_1$ est de cotype fini; $M_2/M_1$ est uniquement $p$-divisible; $M_2/M_1=\cap_k p^k(M/M_1)$,
si bien que $M/M_2$ est $p$-adiquement séparé; comme de plus $Tor_0^\Zp(M/M_2,\Z/p)$ est fini, on voit tout de suite que $\limproj (M/M_2)/p^k$ est
de type fini, donc $M/M_2$ aussi.

Pour montrer $(i)$, on fait un raisonnement en cercle:

- Si $Bid_\Zp(A)\in D^b({_{\underline\G}\C})_{tf}$, alors il en est de même de $Bid_\Zp(A)\Ltens_\Zp\Z/p$, puis de $A\Ltens_\Zp\Z/p$, par $(a)$ et
$(b)$.

- Si $A\Ltens_\Zp\Z/p \in D^b({_{\underline\G}\C})_{tf}$, alors pour tout $q$, on a $H^q(A_n)\Ltens_\Zp\Z/p\in D^b({_\Zp\C})_{tf}$ (observer que la
suite suite spectrale d'hypercohomologie du foncteur $(-)\Ltens_\Zp\Z/p$ dégénère en suites exactes courtes). D'après $(c)$, $H^q(A_n)$ possède donc
une suite de composition comme dans l'énoncé.

- Si les $H^q(A_n)$ possèdent une suite de composition comme dans l'énoncé, alors $Bid_\Zp(H^q(A_n))\in D^b({_\Zp}\C)_{tf}$  (se ramener par
dévissage aux trois cas suivants: 1) si $M\in {_\Zp\C}_{tf}$ est de torsion de cotype fini, alors $RHom_\Zp(M,\Zp)=Hom_\Zp(M,\Qp/\Zp)[-1]$ est de
type fini, donc $Bid_\Zp(M)$ aussi; 2) si $M\in {_\Zp\C}$ est uniquement $p$-divisible, $M\Ltens_\Zp\Qp/\Zp=0$, puis
$RHom_\Zp(M,\Zp)=RHom_\Zp(M\Ltens_\Zp\Qp/\Zp,\Qp/\Zp)=0$ et $Bid_\Zp(A)=0$; 3) si $M$ est de type fini, alors $M$ possède un résolution parfaite, et
$M\simeq Bid_\Zp(M)$ est aussi de type fini), et donc $Bid_\Zp(A)\in D^b({_{\underline\G}\C})_{tf}$.

$(ii)$ Expliquons d'abord l'isomorphisme. Malheureusement, \ref{proplimproj} 3. ne suffit pas, et nous procéderons donc à la main. Soit $(P_k)\in
Kom^b({_\Zp\C}^{\N^\circ})$ (resp. $(Q_k)\in Kom^b({_\Zp\C}^\N)$) le système projectif (resp inductif) évident des complexes
$[\Zp\mathop\rightarrow\limits^{p^k}\Zp]$ placés en degrés $[-1,0]$ (resp. $[0,1]$), de sorte que $(P_k)=Hom_\Zp((Q_k),\Zp)$. Pour $A\in
Kom^b({_{\underline\G}\C})$, on a dans $Kom^b({_{\underline\G}\C}^{\N^\circ})$ des isomorphismes naturels
$$Hom_\Zp(Cofl^*((Q_k)),A)\simeq Fl^*(Hom_\Zp((Q_k),A))\simeq Fl^*(A\otimes_\Zp (P_k))$$ dont celui de l'énoncé se déduit en appliquant $\limproj_k$
(noter le quasi-isomorphis\-me $\limind Q_k\simeq \Qp/\Zp[-1]$).

Passons à l'étude de la seconde flèche de l'énoncé. Dans $D^b({_{\underline\G}\C})$ (resp. $D^b({_\Zp\C^{\N^\circ}})$) on a $A\rightarrow Bid_\Zp(A)$
(resp. $(\Z/p^k)\leftarrow \Zp$); d'où, dans $D^b({_{\underline\G}\C}^{\N^\circ})$: $$A\Ltens_\Zp(\Z/p^k)\mathop\rightarrow\limits^1
 Bid_\Zp(A)\Ltens_\Zp\Z/p^k \mathop\leftarrow\limits^2 Bid_\Zp(A)$$
(le système projectif de droite est constant). On laisse au lecteur le soin de montrer que la flèche ``2'' donne un isomorphisme lorsqu'on lui
applique $R\limproj_k$ (les foncteurs $\pi_n$ et d'oubli permettent de travailler dans ${_{\Zp}\C}$ et ${_\Zp\C}^{\N^\circ}$ au lieu de
${_{\underline\G}\C}$ et ${_{\underline\G}\C}^{\N^\circ}$; on peut alors invoquer un argument de compacité, en mimant le raisonnement de \ref{dupoin}
$(iii)$). La flèche de l'énoncé est donc construite.

Pour mémoire, notons que la flèche ainsi obtenue $RHom_\Zp(\Qp/\Zp,A)[1]\rightarrow Bid_\Zp(A)$ redonne la flèche de bidualité, lorsqu'on la compose
avec la flèche évidente $A=RHom_\Zp(\Zp,A)\rightarrow RHom_\Zp(\Qp/\Zp,A)[1]$.

Reste à vérifier que la flèche ``1'' est un isomorphisme lorsque la condition $(i)$ est vérifiée. Mais dans ce cas, l'analogue modulo $p^k$ de $(i)$
$(a)$ permet d'interpréter celle-ci comme la flèche de bidualité $A\Ltens_\Zp\Z/p^k\rightarrow Bid_\Zp(A\Ltens_\Zp\Z/p^k)$, laquelle est un
isomorphisme, puisque $A\Ltens_\Zp\Z/p^k\in D^b({_{\underline\G}\C})_{tf}$.

$(iii)$ D'après $(ii)$, la flèche de bidualité $A\rightarrow Bid_\Zp(A)$ s'identifie à $A\rightarrow RHom_\Zp(\Qp/\Zp,A)[1]$. De là un triangle
distingué $$A\rightarrow Bid_\Zp(A)\rightarrow RHom_\Zp(\Qp,A)[1]\rightarrow A[1]$$ duquel on tire $A\Ltens_\Zp \Qp/\Zp\simeq
Bid_\Zp(A)\Ltens_\Zp\Qp/\Zp$. Sous nos hypothèses on a de plus $A\Ltens_\Zp\Qp=0$, si bien que la flèche naturelle $A\Ltens_\Zp
\Qp/\Zp[-1]\rightarrow A\Ltens_\Zp \Zp=A$ est un isomorphisme. Le résultat suit. \begin{flushright}$\square$\end{flushright}

\section{Adjonction et dualité}

On présente une étude de la dualité $\L$-linéaire, selon une approche inspirée de \cite{J3}. L'utilisation des catégories dérivées, largement
influencée par \cite{Ne}, montre ici toute son efficacité; elle permet notamment la synthèse et la généralisation des résultats de \cite{V1}.

Comme on a choisi de travailler avec des $\L$-modules non topologisés, il convient de faire une hypothèse de noethérianité: dans cette section et
\textbf{dans tout le reste de l'article}, on fait sur $\G$ les hypothèses suivantes:

\noindent - $\L=\Zp[[\G]]$ est noethérien.

\noindent - La dimension homologique globale de $\L$ est finie (ie. la $p$-dimension cohomologique de $\G$ est finie, cf \cite{Br}). \vspace{0.1cm}

\noindent La seconde hypothèse nous est imposée par le domaine de définition des divers foncteurs dérivés qu'on utilise. Elle est superflue en de
nombreux endroits, si l'on agrandit précautionneusement les domaines de définition en question.  Il est bien connu que la première condition est
vérifiée notamment si $\G$ est un pro-$p$-groupe analytique (cf \cite{La} V.2.2.4). Lorsque de plus la seconde est vérifiée, alors $\G$ est
automatiquement un groupe de Poincaré (\textit{loc. cit} V.2.5.8 et ref.).

\subsection{Un résultat général}

Dans notre contexte, la dualité homologie/cohomologie prend la forme du lemme suivant:

\begin{lem} \label{lem1} Il y a dans $D^b(\C_{\underline\G})$ un isomorphisme fonctoriel en $A\in D^b({_\G\C})$: $$RHom_\Zp(R\underline\G(A),\Zp)\simeq
RHom_\Zp(A,\Zp)\underline\Ltens\Zp$$
\end{lem}
Preuve: D'après \ref{compcohgal} $(ii)$, il y a $D^b({_{\underline\G}\C})$ une flèche fonctorielle en $A$:
$$R\underline\G(A)\rightarrow \underline{RHom}(\Zp,A)$$ Le triangle distingué tautologique
$A\rightarrow A\Ltens_\Zp \Qp\rightarrow A\Ltens_\Zp \Qp/\Zp\rightarrow A[1]$ donne donc lieu donne lieu à un carré commutatif \vspace{-0.5cm}
$$\diagram{RHom_\Zp(R\underline\G(A),\Zp)&\hfl{}{}&RHom_\Zp(R\underline\G(A\Ltens_\Zp \Qp/\Zp),\Zp)[1]\cr \vflupcourte{}{}&&\vflupcourte{}{}\cr
 RHom_\Zp(\underline{RHom}(\Zp,A),\Zp)&\hfl{}{}
&RHom_\Zp(\underline{RHom}(\Zp,A\Ltens_\Zp \Qp/\Zp),\Zp)[1]}$$ \vspace{-0.5cm}

\noindent Les deux flèches horizontales sont des isomorphismes, car $RHom_{\Zp}(-,\Zp)$ s'annule sur tous les complexes dont la cohomologie est
uniquement divisible. La seconde flèche verticale est un isomorphisme par \ref{compcohgal} $(iii)$. La première flèche verticale est donc elle aussi
un isomorphisme, et l'on conclut en invoquant l'isomorphisme d'évaluation \ref{comphomtens} 2. $(ii)$-$(iii)$:
$$RHom_\Zp(A,\Zp)\underline{\Ltens}\Zp\simeq RHom_\Zp(\underline{RHom}(\Zp,A),\Zp)$$
\begin{flushright}$\square$\end{flushright}

\begin{thm} \label{dualite} $ $ \\
\noindent On rappelle que $\underline\L\in {_{\underline\G}\C_\L}$ désigne le système normique canonique $(\Zp[G_n])$.

$(i)$ Pour $M\in D^b(\C_\L)$, il y a dans ${_{\underline\G}\C}$ un isomorphisme fonctoriel:
$$RHom_\Zp(M\underline\Ltens\Zp,\Zp)\simeq RHom_\L(M,\underline\L)$$

$(ii)$ Pour $A\in D^b(_\G\C)$,  il y a dans ${_{\underline\G}\C}$ un isomorphisme fonctoriel:

$$Bid_\Zp(R\underline\G(A))\simeq
RHom_\L(RHom_\Zp(A,\Zp),\underline\L)$$
\end{thm}
Preuve: $(i)$ Par adjonction, puis induction, on a: $$RHom_\Zp(M\underline\Ltens\Zp,\Zp)\simeq R\underline{Hom}(M,RHom_\Zp(\Zp,\Zp))$$
$$\hspace{3cm}\simeq R\underline{Hom}(M,\Zp)$$ $$\hspace{3cm}\simeq RHom_\L(M,\underline{Coind}\Zp)$$
d'où le résultat, puisque dans ${_{\underline\G}\C_\L}$: $\underline\L\simeq\underline{Coind}\Zp$.

$(ii)$ Appliquant $RHom_\Zp(-,\Zp)$ au lemme \ref{lem1}, on obtient dans $D^b({_{\underline\G}\C})$: $$Bid_\Zp(R\underline\G(A))\simeq
RHom_\Zp(RHom_\Zp(A,\Zp)\underline\Ltens\Zp,\Zp)$$ On conclut par $(i)$, appliqué à $M=RHom_\Zp(A,\Zp)\in D^b(\C_\L)$.
\begin{flushright}$\square$\end{flushright}

\begin{cor} \label{cordualite}
Il y a dans $D^b({_\L\C})$ des isomorphismes fonctoriels:

$(i)$ Pour $M\in D^b(\C_\L)$, $RHom_\Zp(\limind M\underline\Ltens\Zp,\Zp)\simeq RHom_\L(M,\L)$.

$(ii)$ Pour $A\in D^b({_\G\C})$, $R\limproj Bid_\Zp(R\underline\G(A))\simeq RHom_\L(RHom_\Zp(A,\Zp),\L)$.
\end{cor}
Preuve: $(i)$ Les compatibilités entre homomorphismes et limites donnent ici (\ref{proplimproj} 3. et remarque):

$$RHom_\Zp(\limind M\underline\Ltens\Zp,\Zp)\simeq R\limproj RHom_\Zp(M\underline\Ltens\Zp,\Zp)$$ $$\hspace{4cm}\simeq R\limproj
RHom_\L(M,\underline\L)$$ $$\hspace{4cm} \simeq RHom_\L(M,R\limproj \underline\L)$$ $$\hspace{4cm}\simeq RHom_\L(M,\L)$$

De même, $(ii)$ découle de \ref{dualite} en appliquant $R\limproj$.
\begin{flushright}$\square$\end{flushright}

\subsection{Le cas des groupes de Poincaré}

\begin{prop}\label{dupoin} On suppose que $\G$ est un groupe de Poincaré de dimension $d$, ie. $D(\G)\simeq \Qp/\Zp[d]$ (cf. \ref{compdu}); alors:

$(i)$ Dans $D^b(\C{_{\underline\G}})$, fonctoriellement en $A\in D^b(_\G\C)$: $$RHom_\Zp(A,\Zp)\underline\Ltens\Zp\simeq RHom_\Zp(\Zp\underline\Ltens
A,\Zp)[d]$$

$(ii)$ Si $A\in D^b(_\G\C)_{tf}$, alors dans $D^b(\C{_{\underline\G}})$: $$RHom_\Zp(R\underline\G(A),\Zp)\simeq
Bid_\Zp(R\underline\G(RHom_\Zp(A,\Zp)))[d]$$ Aussi, si $A\in {_\G\C}$ est fini: $Hom_\Zp(R\underline\G(A),\Qp/\Zp)\simeq
R\underline\G(Hom_\Zp(A,\Qp/\Zp))[d]$.\vspace{0.1cm}

$(iii)$ Si $A\in D^b(_\G\C)$, alors dans $D^b(\C_\L)$: $RHom_\Zp(A,\Zp)\simeq RHom_\L(A,\L)[d]$.
\end{prop}
Preuve: $(i)$ Le théorème \ref{compdu}, ou plutôt \ref{compdurem} $(iii)$, donne ici $$RHom_\Zp(A,\Zp)\underline\Ltens\Zp\simeq
R\underline{Hom}(A,\Zp[d])$$ Le résultat s'en déduit via l'isomorphisme de \ref{dualite} $(i)$:
$$R\underline{Hom}(A,\Zp)\simeq RHom_\Zp(\Zp\underline\Ltens A,\Zp)$$

$(ii)$ Comme $A=Bid_\Zp(A)$, le lemme \ref{lem1} et sa version à droite donnent:

$$RHom_\Zp(A,\Zp)\underline\Ltens\Zp\simeq RHom_\Zp(R\underline\G(A),\Zp)$$
$$\Zp\underline\Ltens A\simeq RHom_\Zp(R\underline\G(RHom_\Zp(A,\Zp)),\Zp)$$ d'où le résultat de l'énoncé, en remplaçant dans $(i)$.

Si $A\in {_\G\C}$ est fini, alors les groupes de cohomologie de $R\underline\G(RHom_\Zp(A,\Zp))$ le sont aussi (car $\L$ est noethérien; utiliser par
exemple \ref{compcohgal}), et l'isomorphisme de l'énoncé suit.

$(iii)$ Par \ref{compdurem} $(iii)$ et induction, on a dans $D^b(\C_{\underline\G})$: $$RHom_\Zp(A,\Zp)\Ltens_\L \underline\L\simeq
RHom(A,\Zp)\underline\Ltens \Zp \simeq R\underline{Hom}(A,\Zp[d])\simeq RHom_\L(A,\underline\L)[d]$$ D'après \ref{proplimproj} 3.,
$RHom_\L(A,\L)\simeq R\limproj RHom_\L(A,\underline\L)$. Reste à montrer que la flèche naturelle de $D^b(\C_\L)$
$$(*)\hspace{2.5cm}RHom_\Zp(A,\Zp)\rightarrow R\limproj (RHom_\Zp(A,\Zp)\Ltens_\L \underline\L)\hspace{3cm}$$ est un isomorphisme
(ce qu'il suffit de vérifier dans $D^b(_\Zp\C)$, par oubli). Maintenant les objets de cohomologie de $RHom_\Zp(A,\Zp)$ sont des $\L$-modules
topologiques compacts, puisque d'après \ref{remlimproj} et \ref{proplimproj} 2. $(ii)$:
$$Ext^q(A,\Zp)\simeq R^q\limproj RHom_\Zp(A_i,\Zp)\simeq \limproj Ext^q_\Zp(A_i,\Zp)$$ si $A=\limind A_i$ avec $A_i$ des complexes à objets de
$\Zp$-type fini. Par troncature du complexe $RHom_\Zp(A,\Zp)$, il suffit donc de montrer que la flèche naturelle de $D^b(\C_\Zp)$
$$(**)\hspace{3.5cm}M\rightarrow R\limproj (M\Ltens_\L \Zp[G_n])\hspace{4cm}$$ est un isomorphisme  pour tout $M\in {\C_\L}$ provenant par oubli
de la catégorie ${\C_{\L-comp}}$ des $\L$-modules (à droite) topologiques compacts.

D'après \cite{Br}, $M$ possède une résolution $P\rightarrow M$, où $P\in Kom^b(\C_{\L-comp})$ est à objets projectifs (dans $\C_{\L-comp}$). Fixons
encore  $L_n\rightarrow \Zp[G_n]$ un système projectif de résolutions parfaites dans ${_\L\C}$. Ont alors lieu dans $Kom^b({_\Zp\C}^{N^\circ})$ des
quasi-isomorphismes (pour le premier, cf \cite{Br} 2.1):
$$P\hat\otimes \Zp[G_n]\leftarrow P\hat\otimes_\L L_n=P\otimes_\L L_n \rightarrow M\otimes_\L L_n$$
Ici $\hat\otimes_\L:{\C_{\L-comp}}\times {_{\L-comp}\C}\rightarrow {_\Zp\C}$ désigne le produit tensoriel complété (il coïncide avec $\otimes_\L$
lorsque le second argument est dans ${_\L\C}_{tf}$, vue comme sous-catégorie pleine de ${_{\L-comp}\C}$). Aussi $P\hat\otimes_\L \Zp[G_n]\simeq
M\Ltens_{\L}\Zp[G_n]$ dans $D^b({_\Zp\C}^{N^\circ})$. Comme les objets du complexe $P\hat\otimes_\L \Zp[G_n]$ sont des systèmes projectifs de
$\Zp$-modules topologiques compacts, donc déployants pour $R\limproj$ (\ref{proplimproj} 2. $(ii)$), la flèche $(**)$ est représentée par
$P\rightarrow \limproj P\otimes_\L \Zp[G_n]$. C'est bien un isomorphisme, puisque les objets de $P$ sont compacts.
\begin{flushright}$\square$\end{flushright}

\begin{cor}\label{dupoin2} Si $\G$ est un groupe de Poincaré de dimension $d$, alors, fonctoriellement en $A\in D^b({_\G\C})$:
$$R\limproj Bid_\Zp R\underline\G(A)\simeq Bid_\L(A)[-d]$$\end{cor}
Preuve: On réécrit \ref{cordualite} $(ii)$ en tenant compte de \ref{dupoin} $(iii)$.
\begin{flushright}$\square$\end{flushright}

\section{Descente et codescente} \label{descetcod}

Si $A\in D^b({_{\underline\G}\C})$ est de la forme $R\underline\G(A_\infty)$ pour un certain $A_\infty\in D^b({_\G\C})$, la connaissance de $A$
équivaut à celle de $\limind A\simeq A_\infty$. Dans ces conditions, il y a un isomorphisme naturel, donné par le corollaire \ref{cordualite}:
$$X_\infty\simeq RHom_\L(RHom_\Zp(A_\infty),\Zp),\L)$$
et cela répond à la question \ref{q2}.

Le but de cette section est la construction, pour $A\in D^b({_{\underline\G}\C})$ quelconque, d'un triangle distingué qui généralise l'isomorphisme
précédent. Comme expliqué dans l'introduction, on en donne en fait deux contructions duales, dont la comparaison aboutit à une relation non triviale
entre le ``défaut de descente'' du sytème $A$ et son ``défaut de codescente''. \vspace{0.2cm}

\subsection{Triangles de (co)-descente: $j(A)$ et $k(A)$} \label{parcod}

A chaque objet $A\in D^b(_{\underline\G}\C)$, sont associés fonctoriellement les objets suivants: \vspace{0.1cm}

\noindent - $A_n:=\pi_n(A)\in D^b(_{G_n}\C)$, pour chaque $n\in N$ (cf \ref{rempin}).

\noindent - $A^*:=RHom_\Zp(A,\Zp)\in D^b(\C_{\underline\G})$, \hbox{$A^*_n:=\pi_n(A^*)\simeq RHom_\Zp(A_n,\Zp)\in D^b(\C_{G_n})$.}

\noindent - $A_\infty:=\limind A\in D^b({_\G\C})$ et $A_\infty^*:=\limind A^*\in D^b(\C_\G)$.

\noindent - $X_\infty:=R\limproj A\in D^b({_\L\C})$ et $X_\infty^*:=R\limproj A^*\simeq RHom_\Zp(A_\infty,\Zp)\in D^b(\C_\L)$. \vspace{0.1cm}

On rappelle que pour $A\in Kom^b({_{\underline\G}\C})_{tf}$, la flèche canonique $\limproj A\rightarrow X_\infty$ est un isomorphisme dans
$D^b({_\L\C})$ (\ref{jensen}), ie. $H^q(X_\infty)\simeq \limproj H^q(A_n)$. Commençons par un critère de noethérianité pour $X_\infty$ et
$X_\infty^*$.

\begin{lem}\label{nakayama} Soit $\G_n$ un sous-groupe ouvert de $\G$.

\noindent 1. Soit $A\in {{_{\underline\G}\C}_{tf}}$, alors $X_\infty=\limproj A\in {{_\L\C}_{tf}}$ si et seulement si $\Zp\otimes_{\Zp[[\G_n]]}
X_\infty$ est de $\Zp$-type fini. Si c'est le cas, alors $\Zp\Ltens_{\Zp[[\G_n]]}X_\infty\in D^b(_\Zp\C)_{tf}$.

\noindent 2. Soit $A\in D^b({_{\underline\G}\C})_{tf}$, alors:

$(i)$ $X_\infty\in D^b({_\L\C})_{tf}$ si et seulement si $\Zp\Ltens_{\Zp[[\G_n]]}X_\infty\in D^b(_{G_n}\C)_{tf}$.

$(ii)$ $X_\infty^*\in D^b({_\L\C})_{tf}$ est si et seulement $RHom_\Zp(R\G(\G_n,A_\infty),\Zp)\in D^b(\C_{G_n})_{tf}$. Cela se produit si et
seulement si $Ext^p_\Zp(R^q\G(\G_n,A_\infty),\Zp)$ est de type fini sur $\Zp$ pour $p=0,1, \, q\in \Z$.
\end{lem} Preuve:  1.  $X_\infty$ étant profini, on peut le munir de sa topologie compacte, et celle-ci est compatible avec l'action de $\L$. On
applique alors la version topologique du lemme de Nakayama.

\noindent 2. $(i)$  On va appliquer 1. aux $\L$-modules $H^q(X_\infty)=\limproj H^q(A_n)$. Observons la suite spectrale d'hypercohomologie:
$$E_2^{p,q}=Tor^{\Zp[[\G_n]]}_{-p}(\Zp,H^q(X_\infty))\Rightarrow Tor^{\Zp[[\G_n]]}_{p+q}(\Zp,X_\infty)=E^{p+q}$$

Si $X_\infty\in D^b({_\L\C})_{tf}$, alors les $E_2^{p,q}$ sont de type fini sur $\L$, donc les $E^{p+q}$ aussi, ie.
$\Zp\Ltens_{\Zp[[\G_n]]}X_\infty\in D^b(_\Zp\C)_{tf}$.

Réciproquement, supposons les $E^{p+q}$ de type fini sur $\L$, et montrons que les $H^q(X_\infty)$ le sont aussi. Ayant remarqué les équivalences
suivantes, données par le point 1: $$H^q(X_\infty)\in ({_\L\C})_{tf}\Leftrightarrow E_2^{0,q}\in ({_\L\C})_{tf}\Leftrightarrow (\forall p,\,
E_2^{p,q}\in ({_\L\C})_{tf})$$ on procède par récurrence descendante sur $q$, pour montrer que $E_2^{0,q}\in ({_\L\C})_{tf}$ (noter que $E_2^{0,q}=0$
pour $q>>0$): le terme initial $E_2^{0,q}=\Zp\otimes_{\Zp[[\G_n]]} H^q(X_\infty)$ est une extension d'un sous-module de
$E^{q}=Tor^{\Zp[[\G_n]]}_{q}(\Zp,X_\infty)$ (de type fini par hypothèse) par un $\L$-module qui possède une filtration dont les graduations sont des
quotients des  $E_2^{-r-1,q+r}=Tor^{\Zp[[\G_n]]}_{r+1}(\Zp,H^{q+r}(X_\infty))$, $r\ge 1$, lesquels sont noethériens par hypothèse de récurrence.

$(ii)$ Le point précédent possède bien sûr une variante à droite; appliquons celle-ci à $X_\infty^*=R\limproj A^*$. Comme
$X_\infty^*\Ltens_{\Zp[[\G_n]]}\Zp\simeq RHom_\Zp(R\G(\G_n,A_\infty),\Zp)$ (cf. \ref{lem1}), le résultat se lit sur la suite spectrale
$$Ext^p_\Zp(R^{-q}\G(\G_n,A_\infty),\Zp)\Rightarrow Ext^{p+q}_\Zp(R\G(\G_n,A_\infty),\Zp)$$
dans laquelle $E_2^{p,q}=0$ pour $p\ne 0,1$.
\begin{flushright}$\square$\end{flushright}

Passons à la construction des ``triangles de descente'' (resp.  de codescente). Pour chaque $A_n$ il existe dans $D^b(_{G_n}\C)$ un morphisme naturel
de descente $j_{A,n}:A_n\rightarrow R\G(\G_n,A_\infty)$ (resp. de codescente $k_{A,n}:\Zp\Ltens_{\Zp[[\G_n]]}X_\infty\rightarrow A_n$). Il s'agit de
relever cette collection de morphismes en une flèche de $D^b({_{\underline\G}\C})$, et de munir cette dernière d'un cône fonctoriel.

\begin{prop}\label{deftr} Soit $Tr(D^b({_{\underline\G}\C}))$ la catégorie des triangles distingués de $D^b({_{\underline\G}\C})$.

$(i)$ Il existe un foncteur $j:D^b({_{\underline\G}\C})\rightarrow Tr(D^b({_{\underline\G}\C}))$: $$A\hspace{0.3cm}\mapsto
\hspace{0.3cm}j(A):=\left(A\mathop\rightarrow\limits^{j_A} R\underline \G(A_\infty)\rightarrow C(j_A)\rightarrow A[1]\right)$$ tel que
$\pi_n(j_A)=j_{A,n}$.

$(ii)$ Il existe un foncteur $k:D^b({_{\underline\G}\C})_{tf}\rightarrow Tr(D^b({_{\underline\G}\C}))$: $$A\hspace{0.3cm}\mapsto \hspace{0.3cm}
k(A):=\left(\Zp\underline\Ltens X_\infty\mathop\rightarrow\limits^{k_A} A\rightarrow C(k_A)\rightarrow \Zp\underline \Ltens X_\infty[1]\right)$$ tel
que $\pi_n(k_A)=k_{A,n}$.
\end{prop}
Preuve: $(i)$ Il suffit de construire un foncteur $\tilde{j}:Kom^b({_{\underline\G}\C})\rightarrow Tr(Kom^b({_{\underline\G}\C}))$ préservant les
quasi-isomorphismes, et qui vérifie les conditions suivantes:

- les premier et second termes du triangle $\tilde{j}(A)$ ont respectivement pour image $A$ et $R\underline \G(A_\infty)$ dans
$D^b({_{\underline\G}\C})$

- la première flèche de $\tilde j(A)$ a pour image $j_{A,n}$ dans $D(_{G_n}\C)$.

\noindent Une fois $\tilde j$ construit, on définit $j(A)$ comme l'image de $\tilde j(A)$ dans $D^b({_{\underline\G}\C})$. Pour construire $\tilde
j$, considérons $K^*:(-):Kom^b({_\G\C})\rightarrow Kom^b({_\G\C})$, le foncteur résolvant construit par Verdier dans \cite{Se2}, annexe au chapitre
I. Celui-ci est muni d'un quasi-isomorphisme fonctoriel en $B$: $B\rightarrow K^*(B)$, et les objets du complexe $K^*(B)$ sont $\G$-cohomologiquement
triviaux, de sorte que $\G(\G_n,K^*(B))$ représente $R\G(\G_n,B)$. Par définition, le morphisme de descente $j_{A,n}$ est alors représenté par la
flèche composée $A_n\rightarrow \G(\G_n,\limind A)\rightarrow \G(\G_n,K^*(\limind A))$. Comme remarqué dans les préliminaires, cette flèche respecte
la structure de système normique, et donne une flèche $\tilde{j}_A:A\rightarrow \underline\G(K^*(\limind A))$ dans $Kom^b({_{\underline\G}\C})$. A
donc lieu dans $Kom^b({_{\underline\G}\C})$ un triangle distingué tautologique, fonctoriel en $A\in Kom^b({_{\underline\G}\C})$:
$$\tilde{j}(A):=\hspace{0.3cm}\left(A\mathop\rightarrow\limits^{\tilde{j}_A} \G(\G_n,K^*(\limind A))\rightarrow C(\tilde{j}_A)\rightarrow A[1]\right)$$
Celui-ci vérifie bien les deux propriétés annoncées, et cela termine la preuve.

$(ii)$ Fixons une fois pour toutes une résolution projective $P\rightarrow \Zp$ dans $D^b({_\L\C})$. Si $A\in Kom^b({_{\underline\G}\C})_{tf}$, alors
$X_\infty\in D^b({_\L\C})$ est représenté par $\limproj A$ (\ref{jensen}), et la flèche de codescente $k_{A,n}$ de $D^b(_{G_n}\C)$, par la flèche
composée: $P\otimes_{\Zp[[\G_n]]}\limproj A\rightarrow \Zp\otimes_{\Zp[[\G_n]]}\limproj A\rightarrow A_n$. Comme en $(i)$, on note que cette dernière
flèche a en fait lieu dans $Kom^b({_{\underline\G}\C})$, et donne un triangle distingué tautologique fonctoriel dont le foncteur $k$ est déduit par
localisation.
\begin{flushright}$\square$\end{flushright}

\begin{rem} \label{remdeftr} 1. Soit $A\in Kom^b(_{\underline\G}\C)_{tf}$. Alors:

$(i)$ Les flèches $j_A$ et $k_A$ ne dépendent pas des résolutions qu'on a choisies. Par ailleurs, on peut aussi expliciter $k_A$ à l'aide d'une
résolution projective de $\limproj A$: soit $P\rightarrow \limproj A$ une telle résolution, alors $k_A$ est représentée par la flèche composée
$\Zp\underline\otimes P\rightarrow \Zp\underline\otimes \limproj A\rightarrow A$.

$(ii)$  Remplacer $K^*(A)$ par une autre résolution $\G$-acyclique de $A$ dans la preuve ci-dessus remplace $j(A)$ par un autre triangle distingué,
qui lui est (non canoniquement) isomorphe. Remarque analogue pour $P$ et $k_A$.

2. Soit $A\in D^b({_{\underline\G}\C})$, non dans $D^b({_{\underline\G}\C})_{tf}$. Si $A$ est représenté par un complexe de ${_{\underline\G}\C}$
dont les objets sont acycliques pour le foncteur $R\limproj$, alors la construction de $k(A)$ tient encore. J'ignore si l'on peut trouver un tel
complexe pour tout $A\in D^b({_{\underline\G}\C})$. \end{rem}

\begin{rem}\label{remnoethjk} Soit $A\in D^b({_{\underline\G}\C})_{tf}$, alors \ref{nakayama} montre que:

$(i)$  $X_\infty\in D^b(_\L\C)_{tf}\Leftrightarrow C(k_A)\in D^b({_{\underline\G}\C})_{tf}$.

$(ii)$ $X_\infty^*\in D^b(\C_\L)_{tf}\Leftrightarrow RHom_\Zp(C(j_A),\Zp)\in D^b(\C_{\underline\G})_{tf}$
\end{rem}

\begin{prop} \label{compjk} Notons encore $k:D^b(\C_{\underline\G})_{tf}\rightarrow Tr(D^b(\C_{\underline\G}))$ l'analogue à droite du
foncteur $k$ de la proposition précédente. Si $A^*\in D^b({_{\underline\G}\C})_{tf}$, il y a dans $D^b({_{\underline\G}\C})$ un diagramme commutatif,
fonctoriel en $A\in D^b({_{\underline\G}\C})$:
$$\diagram{(X_\infty^* \underline\Ltens\Zp)&\hfl{k_{A*}}{}&A^* \cr \vfl{ev}{}&&\vfl{id}{}\cr
RHom_\Zp(R\underline\G(A_\infty),\Zp)&\hfl{RHom_\Zp(j_A,\Zp)}{}&RHom_\Zp(A,\Zp)}$$
\end{prop}
Preuve: Soit $A\in Kom^b({_{\underline\G}\C})$, tel que $A^*\in D^b(\C_{\underline\G})_{tf}$. Si $I\in Kom^b(_\Zp\C)$ est une résolution injective de
$\Zp$, alors $X_\infty^*\in D^b(\C_\L)$ est représenté par le complexe $\limproj Hom_\Zp(A,I)=Hom_\Zp(\limind A,I)$. Soient maintenant $P\rightarrow
\Zp$ et $\limind A\rightarrow K^*(\limind A)$ les résolutions de la preuve précédente, et considérons le diagramme commutatif de
$Kom^b(\C_{\underline\G})$:

{\small
$$\diagram{\limproj Hom_\Zp(A,I)\underline\otimes P&\hflcourte{}{}&\limproj Hom_\Zp(A,I)\underline\otimes \Zp&\hflcourte{}{}&Hom_\Zp(A,I)\cr
\vflupcourte{}{}&&\vflcourte{\wr}{}&&\vflcourte{id}{}\cr
Hom_\Zp(\limind K^*(A),I)\underline\otimes P&\hflcourte{}{}&Hom_\Zp(\limind A,I)\underline\otimes \Zp&\hflcourte{}{}&Hom_\Zp(A,I)\cr
\vflcourte{ev}{}&&\vflcourte{ev}{}&&\vflcourte{id}{}\cr  Hom_\Zp(\underline{Hom}(P,K^*(\limind A)),I) &\hflcourte{}{}&Hom_\Zp(\underline{Hom}(\Zp,
\limind
A),I)&\hflcourte{}{}&Hom_\Zp(A,I)\cr
\vflcourte{}{}&&\vflcourte{\wr}{}&&\vflcourte{id}{}\cr  Hom_\Zp(\underline\G(K^*(\limind A)),I) &\hflcourte{}{}&Hom_\Zp(\underline\G(\limind
A),I)&\hflcourte{}{}&Hom_\Zp(A,I)}$$}

\noindent Dans celui-ci, la ligne supérieure (resp. inférieure) représente $k_{A^*}$ (resp. \break $RHom_{\Zp}(j_A,\Zp)$) et toutes les flèches
verticales sont des quasi-isomorphismes, ainsi que les flèches horizontales de gauche. On obtient le carré commutatif de l'énoncé en identifiant,
dans $D^b({_{\underline\G}\C})$, les deux premières colonnes, ainsi que les deux premières et les deux dernières lignes.
\begin{flushright}$\square$\end{flushright}

\subsection{Triangles d'adjonction: $\alpha_j(A)$ et $\alpha_k(A)$}

Nous sommes maintenant en mesure d'établir le

\begin{thm} \label{theoprincipal} Soit $A\in D^b({_{\underline\G}\C})_{tf}$, alors

$(i)$ Dans $D^b({_{\underline\G}\C})$, fonctoriellement en $A$: $C(j_A)\simeq C(j_{C(k_A)})$. Si de plus $X_\infty^*\in D^b(\C_\L)_{tf}$, alors
$Bid_\Zp C(j_A)\in D^b({_{\underline\G}\C})_{tf}$ et $C(k_{A})\simeq C(k_{Bid_\Zp C(j_A)})[-1]$.

$(ii)$ Le foncteur $j$ donne lieu à un triangle distingué, fonctoriel en $A$: $$\alpha_j(A):=\left( X_\infty\mathop\rightarrow\limits^{\alpha_j}
RHom_\L(X^*_\infty,\L)\rightarrow \Delta_j(A)\rightarrow X_\infty[1]\right)$$ dans lequel $\Delta_j(A):=R\limproj Bid_\Zp(C(j_A))$.

$(iii)$ Le foncteur $k$ donne lieu à un triangle distingué, fonctoriel en $A$: $$\alpha_k(A):=\left( X_\infty^*\mathop\rightarrow\limits^{\alpha_k}
RHom_\L(X_\infty,\L)\rightarrow\Delta_k(A)\rightarrow  X_\infty^*[1]\right)$$ dans lequel $\Delta_k(A):=RHom_\Zp(\limind C(k_A),\Zp)[1]$.

$(iv)$ Supposons de plus $X_\infty\in D^b({_\L\C})_{tf}$, alors il y a un isomorphisme de triangles $\alpha_j(A)\simeq RHom_\L(\alpha_k(A),\L)$
fonctoriel en $A$, ie: {\small
$$\diagram{X_\infty&\hflcourte{\alpha_j}{}& RHom_\L(X^*_\infty,\L)&\hflcourte{}{} &\Delta_j(A)&\hflcourte{}{}& X_\infty[1]\cr
\vfl{}{}&&\vfl{id}{}&&\vfl{}{}&&\vfl{}{}\cr
Bid_\L(X_\infty)&\hflcourte{RHom_\L(\alpha_k,\L)}{}&RHom_\L(X_\infty^*,\L)&\hflcourte{}{}&RHom_\L(\Delta_k(A),\L)[1]&\hflcourte{}{}
&Bid_\L(X_\infty)[1]}$$} En particulier, $R\limproj Bid_\Zp(C(j_A))\simeq RHom_\L(RHom_\Zp(\limind C(k_A),\Zp),\L)$.
\end{thm}

\noindent Preuve: $(i)$ Considérons $B:=\Zp\underline\Ltens X_\infty\in D^b({_{\underline\G}\C})$, et appliquons le foncteur $j$ aux objets et
flèches du triangle $k(A)=\left(B\mathop\rightarrow\limits^{k_A} A\rightarrow C(k_A)\rightarrow B[1]\right)$. On obtient alors un diagramme
commutatif dont chaque ligne et colonne est un triangle distingué de $D^b({_{\underline\G}\C})$:
$$\diagram{B&\hfl{j_B}{}&R\underline\G(B_\infty)&\hfl{}{}&C(j_B)&\hfl{}{}&B[1]\cr
\vflcourte{}{}&&\vflcourte{R\underline\G(k_A)}{}&&\vflcourte{}{}&&\vflcourte{}{}\cr
A&\hfl{j_A}{}&R\underline\G(A_\infty)&\hfl{}{}&C(j_A)&\hfl{}{}&A[1]\cr
\vflcourte{}{}&&\vflcourte{}{}&&\vflcourte{}{}&&\vflcourte{}{}\cr
C(k_A)&\hfl{j_{C(k_A)}}{}&R\underline\G(C(k_A)_\infty)&\hfl{}{}&C(j_{C(k_A)})&\hfl{}{}&C(k_A)[1]\cr
\vflcourte{}{}&&\vflcourte{}{}&&\vflcourte{}{}\cr
B[1]&\hfl{j_B[1]}{}&R\underline\G(B_\infty)[1]&\hfl{}{}&C(j_B)[1]}$$

\noindent Pour obtenir le premier isomorphisme de l'énoncé, il suffit de montrer que $C(j_B)=0$, ie. que $j_B$ est un isomorphisme. Soit $P\in
Kom^b({_\L\C})$ un complexe à objets projectifs représentant $X_\infty$. Comme les objets du complexe $\limind \Zp\underline\otimes P$ sont
$\G$-cohomologiquement triviaux, la flèche $j_B$ est représentée par $\Zp\underline\otimes P\rightarrow \underline\G(\limind \Zp\underline\otimes
P)$. C'est un isomorphisme de $Kom^b(_{\underline \G}\C)$, puisque les objets de $P$ sont des facteurs directs de $\L$-modules libres, et que
$C(j_{\underline \L})=0$.

Passons au second isomorphisme. L'hypothèse de finitude sur $X_\infty^*$ nous assure que $D:=Bid_\Zp R\underline\G(A_\infty)$ est dans
$D^b({_{\underline\G}\C})_{tf}$. Aussi le triangle distingué
$$Bid_\Zp (j(A))=\left(Bid_\Zp A\mathop\rightarrow\limits^{Bid_\Zp(j_A)} D\rightarrow Bid_\Zp C(j_A)\rightarrow Bid_\Zp A[1]\right)$$
a-t-il lieu dans $D^b({_{\underline\G}\C})_{tf}$. On peut donc lui appliquer le foncteur $k$ et cela donne en particulier un nouveau triangle
distingué:
$$C(k_{Bid_\Zp A})\rightarrow C(k_D)\rightarrow C(k_{Bid_\Zp C(j_A)})\rightarrow C(k_{Bid_\Zp A})[1]$$ Comme $A\simeq Bid_\Zp(A)$,
il reste à montrer que  $C(k_D)=0$. Pour ce faire, nous utilisons le théorème \ref{dualite}, selon lequel $D\simeq RHom_\L(X_\infty^*,\underline\L)$.
Soit $P\in Kom^b(\C_\L)$ un complexe parfait représentant $X_\infty^*$, de sorte que le complexe $Hom_\L(P,\L)\in Kom^b({_\L\C})$ est parfait aussi.
Mais alors $k_D$ est représentée par $\Zp\underline\otimes Hom_\L(P,\L)\rightarrow Hom_\L(P,\underline\L)$; c'est visiblement un isomorphisme, et
l'on a bien $C(k_D)=0$.

$(ii)$ Appliquer le foncteur $R\limproj\circ Bid_\Zp(-)$ au triangle distingué $j(A)$ en donne un nouveau: {\small $$R\limproj Bid_\Zp(A)\rightarrow
R\limproj Bid_\Zp(R\underline\G(A_\infty))\rightarrow R\limproj Bid_\Zp C(j_A)\rightarrow R\limproj Bid_\Zp(A)[1]$$} D'où celui de l'énoncé, compte
tenu du corollaire \ref{cordualite} $(ii)$ appliqué à $A_\infty\in D^b(_\G\C)$, et de l'isomorphisme $A\simeq Bid_\Zp(A)\in
D^b({_{\underline\G}\C})_{tf}$.

$(iii)$ Appliquer le foncteur $RHom_\Zp(-,\Zp)\circ\limind$ au triangle distingué $k(A)$ en donne un nouveau:
$$X_\infty^*\rightarrow RHom_\Zp(\limind \Zp\underline\Ltens X_\infty,\Zp)\rightarrow RHom_\Zp(\limind C(k_A),\Zp)[1]\rightarrow
X_\infty^*[1]$$ D'où celui de l'énoncé, puisque d'après l'analogue à gauche de \ref{cordualite} $(i)$:
$$RHom_\Zp(\limind \Zp\underline\Ltens X_\infty,\Zp)\simeq RHom_\L(X_\infty,\L)$$

$(iv)$ Soit comme en $(i)$, $B=\Zp\underline\Ltens X_\infty$. L'hypothèse de finitude sur $X_\infty$ nous assure que $B,C(k_A)\in
D^b({_{\underline\G}\C})_{tf}$. Sachant que $A\simeq Bid_\Zp(A)$, $B\simeq Bid_\Zp(B)$, $C(k_A)\simeq Bid_\Zp(C(k_A))$ et $C(j_B)=0$, on obtient en
appliquant $R\limproj\circ Bid_\Zp(-)$ au premier diagramme de la preuve:

 {\small
$$\diagram{R\limproj B&\hflcourte{}{}&R\limproj Bid_\Zp(R\underline\G(B_\infty))&\hflcourte{}{}&
0&\hflcourte{}{}&R\limproj B[1]\cr \vflcourte{}{}&&\vflcourte{}{}&&\vflcourte{}{}&&\vflcourte{}{}\cr X_\infty&\hflcourte{}{}&R\limproj
Bid_\Zp(R\underline\G(A_\infty))&\hflcourte{}{}&R\limproj Bid_\Zp(C(j_A)) &\hflcourte{}{}&X_\infty[1]\cr
\vflcourte{}{}&&\vflcourte{}{}&&\vflcourte{}{}&&\vflcourte{}{}\cr R\limproj C(k_A)&\hflcourte{}{}&R\limproj
Bid_\Zp(R\underline\G(C(k_A)_\infty))&\hflcourte{}{}&R\limproj Bid_\Zp(C(j_{C(k_A)}))&\hflcourte{}{}&R\limproj C(k_A)[1]\cr
\vflcourte{}{}&&\vflcourte{}{}&&\vflcourte{}{}\cr R\limproj B[1]&\hflcourte{}{}&R\limproj Bid_\Zp(R\underline\G(B_\infty))[1]&\hflcourte{}{}&0}$$}

Dans ce diagramme commutatif, chaque ligne ou colonne est un triangle distingué de $D^b({_\L\C})$. De plus:

\noindent - La seconde ligne du diagramme s'identifie au triangle $\alpha_j(A)$, par construction de ce dernier.

\noindent - La seconde colonne s'identifie à $RHom_\L(\alpha_k(A),\L)$, comme cela résulte des isomorphismes de triangles distingués suivants
(\ref{cordualite}):
$$R\limproj Bid_\Zp(R\underline\G(\limind k(A)))\simeq RHom_\L(RHom_\Zp(\limind k(A),\Zp),\L)$$ $$\hspace{4cm}\simeq RHom_\L(\alpha_k(A),\L)$$

Pour obtenir l'isomorphisme souhaité entre les triangles $RHom_\L(\alpha_k(A),\L)$ et  $\alpha_j(A)$, il suffit donc de montrer que $R\limproj
C(k_A)=0$. Mais cela relève du même argument que la preuve de \ref{dupoin} $(iii)$.
\begin{flushright}$\square$\end{flushright}

\begin{cor} \label{corcjk} Si $A\in D^b({_{\underline\G}\C})_{tf}$ et $X_\infty^*\in D^b(\C_\L)_{tf}$, alors $$C(j_A)=0\Leftrightarrow C(k_A)=0$$
\end{cor}
Preuve: Cela résulte immédiatement de \ref{theoprincipal} $(i)$.
\begin{flushright}$\square$\end{flushright}

\begin{cor} \label{poincare} Si $\G$ est un groupe de Poincaré de dimension $d$ (par exemple $\G\simeq \Zp^d$),
alors $$R\limproj Bid_\Zp(C(j_A))\simeq Bid_\L(\limind C(k_A))[-d]$$ \end{cor} Preuve: On réécrit simplement \ref{theoprincipal} $(iv)$, à l'aide de
\ref{dupoin} $(iii)$. \begin{flushright}$\square$\end{flushright}

On relève un critère de noethérianité pour $X_\infty^*$, utile en pratique:
\begin{cor} \label{critnoethcod} On suppose $A\in D^b({_{\underline\G}\C})_{tf}$
et $X_\infty\in D^b({_\L\C})_{tf}$. Alors $X_\infty^*\in D^b(\C_\L)_{tf}$ si et seulement si $RHom_\Zp(\limind C(k_A),\Zp)\in D^b(\C_\L)_{tf}$.
\end{cor}
Preuve: L'équivalence se lit sur le triangle distingué $\alpha_k(A)$.
\begin{flushright}$\square$\end{flushright}

\begin{rem} \label{remtheoprincipal}
Pour $A\in D^b({_{\underline\G}\C})$ non nécessairement dans $D^b({_{\underline\G}\C})_{tf}$, le triangle distingué $\alpha_j(A)$ admet une
généralisation évidente: $$R\limproj Bid_\Zp(A)\rightarrow RHom_\L(X_\infty^*,\L)\rightarrow \Delta_j(A)\rightarrow R\limproj Bid_\Zp(A)[1]$$ De la
même façon, lorsqu'on est dans le cas de \ref{remdeftr} $(iii)$, il est possible de généraliser une partie des résultats.
\end{rem}

\begin{rem} \label{amplitude} Disons qu'un foncteur $F:D^b(\C)\rightarrow D^b(\C')$ covariant (resp. contravariant) est d'amplitude $[a,b]$ s'il possède
la propriété suivante: ``$A$ acyclique hors de $[\alpha,\beta]\Rightarrow F(A)$ acyclique hors de $[\alpha+a,\beta+b]$ (resp.
$[a-\beta,b-\alpha]$)''. On a jusqu'ici défini trois foncteurs exacts covariants: $C(j_{(-)}):D^b(_{\underline\G}C)\rightarrow
D^b({_{\underline\G}\C})$, $\Delta_j:D^b(_{\underline\G}C)_{tf}\rightarrow D^b({_\L\C})$, $C(k_{(-)}):D^b(_{\underline\G}C)_{tf}\rightarrow
D^b(_{\L}\C)$, et un contravariant: $\Delta_k:D^b(_{\underline\G}C)_{tf}\rightarrow D^b(\C_\L)$. Si $cd_p\G=d$, ils sont a priori d'amplitudes
respectives $[-1,d+1]$, $[-1,d]$, $[-d-1,0]$, $[-1,d+1]$ (que $\Delta_j(A)$ soit acyclique en degrés $>\beta+d$ si $A$ l'est en degrés $\ge \beta$ se
lit par exemple sur l'isomorphisme suivant: $\Delta_j(A)\simeq RHom_\Zp(\limind RHom_\Zp(C(j_A),\Zp),\Zp)$, compte-tenu de la nullité de
$Ext^1_\Zp(\limind Ext^1_\Zp(H^{d+1}(\G_n,H^\beta(A_\infty)),\Zp),\Zp)$ dûe à la $p$-divisibilité de $H^{d+1}(\G_n,-)$).
\end{rem}

\begin{q} Que peut-on dire en général concernant la ``taille'' de $\Delta_j(A)$ et $\Delta_k(A)$? En particulier, lorsqu'on dispose d'un foncteur $det$,
peut-on caractériser $det(\Delta_j)$ directement en termes de l'objet $A$? \end{q}

\noindent Si $\G\simeq\Zp^d$ et $A\in {_{\underline\G}\C}$ est de la forme $A=\Zp\underline\otimes X_\infty$, c'est l'objet de la section suivante
que de donner une condition suffisante pour avoir $\Delta_j(A)=0$ dans $D^b(_{\L}\C/\{pseudo-nuls\})$. Le cas où $A$ est quelconque reste assez
mystérieux, même pour $\G\simeq\Zp^d$.

Lorsque $\G\simeq\Zp$, $C(j_A)$ et $C(k_A)$ ont tendance à se stabiliser (en un sens adéquat, cf section \ref{sectionstab}) lorsque $A$ provient de
l'arithmétique (e.g. de la cohomologie d'une représentation $p$-adique); lorsque cela se produit, on parvient à une compréhension satisfaisante des
quatre foncteurs de la remarque précédente.

\section{Le cas commutatif: $\G\simeq\Zp^d$}

Dans ce paragraphe et \textbf{dans toute la suite de l'article}, on suppose que $\G\simeq\Zp^d$, si bien que $\L$ est un anneau commutatif. On
identifie donc ${_\L\C}$ et ${\C_\L}$. L'anneau $\L$ est régulier de dimension $d+1$ (\cite{Se1}); on note parfois $Q$ son corps des fractions. On
désigne les pseudo-isomorphismes par le symbole $\approx$, ou au besoin $\mathop\rightarrow\limits^{\approx}$; rappelons que $\approx$ est une
relation d'équivalence sur la classe des $\L$-modules de torsion. Si $M\in D^b({_\L\C})$ on utilisera les notations suivantes:
$$E^q(M):=Ext^q_\L(M,\L) \hspace{2cm} H_q(\G_n,M):=Tor_q^{\Zp[[\G_n]]}(\Zp,M)$$
\noindent  Pour un $\L$-module $M$, on note $t_\L M$ le sous-module de $\L$-torsion, et $f_\L M:=M/t_\L M$; de même avec $\Zp$ au lieu de $\L$. On
note aussi $M[p^k]=Tor_1^\Zp(M,\Z/p^k)$ (resp. $M/p^k=Tor_0^\Zp(M,\Z/p^k)$) le noyau (resp. conoyau) de la multiplication par $p^k$.

\subsection{$\Delta_k$ et la structure des $\L$-modules} \label{deltakstruc}

Rassemblons en un lemme quelques propriétés des $E^q$. Pour une étude systématique de ces foncteurs, on renvoie à \cite{J3}.

\begin{lem}\label{lemmod} Soit $M$ un $\L$-module de type fini. Alors:

$(i)$ $E^qM$ est de torsion (resp. pseudo-nul) si $q\ge 1$, (resp. $q\ge 2$).

$(ii)$ Si $M$ est de torsion, le bord de la suite spectrale de bidualité induit un pseudo-isomorphisme canonique
$M\mathop\rightarrow\limits^{\approx} E^1E^1M$.

$(iii)$ Si $M$ est de torsion, alors $E^1M$ ne possède aucun sous-module pseudo-nul.

$(iv)$ Il existe un pseudo-isomorphisme (non canonique) $t_\L M\approx E^1M$. En particulier, si $M$ est de torsion, on a les équivalences suivantes:
$$M\approx 0\Leftrightarrow E^1M\approx 0\Leftrightarrow E^1M=0$$
\end{lem}
Preuve (esquisse): $(i)$ se voit par localisation, puisque la dimension homologique globale d'un corps (resp. d'un anneau principal) vaut $0$ (resp.
$1$).

$(ii)$ Que le bord en question soit bien défini résulte de la nullité de $E^0E^0M$. Qu'il soit un pseudo-isomorphisme se vérifie en localisant.

$(iii)$ et $(iv)$ se déduisent facilement des deux faits suivants:

- Si $M$ est de torsion, alors $E^0M=0$ et $E^1M\simeq Hom_\L(M,Frac(\L)/\L)$.

- Si $M$ est sans torsion, alors $E^1M$ est pseudo-nul (localiser).
\begin{flushright}$\square$\end{flushright}

Lorsque $M\in D^b({_\L\C})_{tf}$ vérifie $M\Ltens_\L Q=0$ (ie. lorsque tous les modules de cohomologie de $M$ sont de torsion), on définit l'idéal
caractéristique de $M$ par
$$\chi(M):=\prod (\mathcal Char H^q(M))^{(-1)^q}$$
où $\mathcal Char$ désigne l'idéal caractéristique au sens habituel de la théorie des modules d'Iwasawa.

\begin{prop} \label{charjk} Soit $A\in D^b({_{\underline\G}\C})_{tf}$. On suppose  $X_\infty\in D^b({_\L\C})_{tf}$ et $X_\infty^*\in D^b(\C_\L)_{tf}$.
Alors:

$(i)$  $\Delta_j(A)\otimes_\L Q=0\Leftrightarrow \Delta_k(A)\otimes_\L Q=0$.

$(ii)$ Lorsque $\Delta_j(A)\otimes_\L Q=0$, on a $\chi(\Delta_j(A))=\chi(\Delta_k(A))$.
\end{prop} Preuve: d'après \ref{theoprincipal} $(iv)$, on a $\Delta_j(A)=RHom_\L(\Delta_k(A),\L)[1]$.
\begin{flushright}$\square$\end{flushright}

Dans ce paragraphe, on s'intéresse au cas où $A\in D^b({_{\underline\G}\C})_{tf}$ est concentré en degré $0$, ie. $A\in {{_{\underline\G}\C}}_{tf}$,
et vérifie la condition de codescente galoisienne, version non dérivée, ie. $H^{-1}(C(k_A))=H^0(C(k_A))=0$. Dans cette situation, on souhaite décrire
explicitement le complexe $\Delta_k(A)$.

\begin{lem} \label{deltak} Soit $M\in {{_\L\C}_{tf}}$, et $A:=\Zp\underline\otimes M\in {{_{\underline\G}\C}_{tf}}$, de sorte que $X_\infty\simeq M$.
On a alors:

\noindent 1. $X_\infty^*\in D^b(\C_\L)_{tf}$.

\noindent 2. $\Delta_k(A)\otimes_\L Q=0$, $\Delta_k(A)$ est acyclique en degrés négatifs, et vérifie:

$(i)$ $H^1(\Delta_k(A))=Hom_\Zp(\limind H_1(\G_n,X_\infty),\Zp)$.

$(ii)$ $H^q(\Delta_k(A))=E^q(X_\infty)$ pour $q\ge 2$.
\end{lem}
Preuve: 1. Le triangle distingué $\alpha_k(A)$ montre que $X_\infty^*\in D^b(\C_\L)_{tf}$ équivaut à $\Delta_k(A)\in D^b(\C_\L)_{tf}$, fait que le
calcul de la cohomologie de $\Delta_k(A)$ mettra en évidence ci-dessous.

2. On souhaite calculer $\Delta_k(A)=RHom_\Zp(\limind C(k_A),\Zp)[1]$. Examinons le triangle distingué $k(A)$:
$$\Zp\underline\Ltens X_\infty\mathop\rightarrow\limits^{k_A} A\rightarrow C(k_A)\rightarrow \Zp\underline\Ltens X_\infty[1]$$
Par hypothèse, $A$ est concentré en degré $0$ et $H^0(k_A)$ est un isomorphisme. On obtient donc $C(k_A)[-1]\simeq \tau_{\le -1}\Zp\underline\Ltens
X_\infty$, puis $$\Delta_k(A)\simeq RHom_\Zp(\limind \tau_{\le -1}\Zp\underline\Ltens M,\Zp)$$ Le résultat souhaité (y compris la noethérianité
annoncée en 1.) est alors un cas particulier du lemme suivant, dont nous aurons à nouveau l'utilité un peu plus loin.
\begin{flushright}$\square$\end{flushright}

\begin{lem}\label{calchypcoh} Soit $M\in D^b({{_\L\C}})$, et $a\in \Z$. Alors:

$(i)$ $Ext^q_\Zp(\limind \tau_{\le a}\Zp\underline\Ltens M,\Zp)=0$ pour $q\le -1-a$.

$(ii)$ $Ext^{-a}_\Zp(\limind \tau_{\le a}\Zp\underline\Ltens M,\Zp)\simeq Hom_\Zp(\limind H_{-a}(\G_n,M),\Zp)$.

$(iii)$ $Ext^q_\Zp(\limind \tau_{\le a}\Zp\underline\Ltens M,\Zp)=E^q(M)$ pour $q\ge 1-a$.

$(iv)$ Il y a pour tout $q\in \Z$ une suite exacte naturelle
$$Ext^1_\Zp(\limind H_{q-1}(\G_n,M),\Zp)\hookrightarrow E^q(M)\twoheadrightarrow Hom_\Zp(\limind H_q(\G_n,M),\Zp)$$
\end{lem}

Preuve: D'après \ref{cordualite} $(i)$, on a $Ext^q_\Zp(\limind \Zp\underline\Ltens X_\infty,\Zp)\simeq E^q(X_\infty)$.  Par ailleurs, la suite
spectrale d'hypercohomologie du foncteur $RHom_\Zp(-,\Zp):D^b({_{\G}\C})\rightarrow D^b(\C_\L)$, appliquée à la flèche $\limind\tau_{\le a}\Zp\Ltens
M\rightarrow \limind \Zp\underline\Ltens M$ produit un diagramme commutatif à lignes exactes:{\small
$$ \hspace{-1.5cm}\diagram{Ext^1_\Zp(\limind H^{1-q}(\tau_{\le a}\Zp\underline\Ltens M),\Zp)&\injfl{}{}&Ext_\Zp^{q}(\limind \tau_{\le a}\Zp\underline\Ltens M,\Zp)
&\surjfl{}{}&Hom_\Zp(\limind H^{-q}(\tau_{\le a}\Zp\underline\Ltens M),\Zp)\cr \vflupcourte{}{}&&\vflupcourte{}{}&&\vflupcourte{}{} \cr Ext^1_\Zp(\limind
H_{q-1}(\G_n,M),\Zp)&\injfl{}{}& Ext^q_\Zp(\limind \Zp\underline\Ltens M,\Zp)&\surjfl{}{}&Hom_\Zp(\limind H_q(\G_n,M),\Zp)}$$} Dans celui-ci, les flèches
verticales extrèmes sont, selon les valeurs de $q$, soit nulles soit des isomorphismes. Le résultat suit.
\begin{flushright}$\square$\end{flushright}

En pratique, l'expression de $H^1(\Delta_k(A))$ donnée par le lemme \ref{deltak} est peu maniable. On cherche dans ce qui suit à contrôler sa taille.

\begin{defn} Soit $M\in {_\L\C}_{tf}$ un module de torsion. On munit $M$ de deux filtrations décroissantes fonctorielles:

\noindent - $I^1M=\limproj t_\Zp(M_{\G_n})$, et $I^{q+1}M=I^1(I^q(M))$ pour $q\ge 1$.

\noindent - $F^1M$ est le noyau de l'application composée $$M\rightarrow E^1E^1M\rightarrow Hom_\Zp(\limind H_1(\G_n,E^1M),\Zp)$$ où la première
flèche est celle de \ref{lemmod} $(ii)$, et la seconde est la flèche de droite dans la suite exacte \ref{calchypcoh} $(iv)$, appliquée au module
$E^1M$. Puis $F^{q+1}M=F^1F^{q}M$ pour $q\ge 1$.
\end{defn}

\begin{rem} \label{FF}
La flèche composée qui définit $F^1M$ est pseudo-surjective, si bien que $M/F^1M\hookrightarrow Hom_\Zp(\limind H_1(\G_n,E^1M),\Zp)$ avec conoyau
pseudo-nul.
\end{rem}

\begin{prop} \label{IF} Soit $M\in {{_\L\C}_{tf}}$ un module de torsion, alors $I^qM\subset F^qM$.
\end{prop}
Preuve: On peut toujours supposer $q=1$. Il s'agit alors de montrer que l'image de $\limproj t_\Zp(M_{\G_n})$ par l'application $M\rightarrow
Hom_\Zp(\limind H_1(\G_n,E^1M),\Zp)$ est nulle. C'est évident si l'on écrit cette dernière comme la limite projective des applications
$M_{\G_n}\rightarrow Hom_\Zp(H_1(\G_n,E^1M),\Zp)$. \begin{flushright}$\square$\end{flushright}

\begin{q}\label{questionIF} Est-il vrai que $I^qM=F^qM$? \end{q}

Cette question possède une réponse affirmative si $d=1$, ie. si $\G\simeq \Zp$. En effet, on montre alors facilement que $M/I^1M$ et $M/F^1M$ sont
tous deux $\Zp$-libres de même rang, et le résultat suit.

\begin{lem} \label{pseudoIF}
Soient $M,M'\in {_\L\C_{tf}}$ de torsion. S'il y a un pseudo-isomorphisme $M\mathop\rightarrow\limits^\approx M'$, celui-ci induit
des pseudo-isomorphismes:

$(i)$ $I^qM\mathop\rightarrow\limits^\approx  I^q M'$ et $F^qM\mathop\rightarrow\limits^\approx  F^qM'$.

$(ii)$ $Hom_\Zp(\limind H_1(\G_n,M'),\Zp)\mathop\rightarrow\limits^\approx  Hom_\Zp(\limind H_1(\G_n,M),\Zp)$.
\end{lem} Preuve: $(i)$ Soit $J=I$ ou
$F$. Tout revient à montrer que la flèche $M/J^1M\rightarrow M'/J^1M'$ est un pseudo-isomorphisme. Or

- son conoyau est pseudo-nul, par hypothèse.

- il existe un pseudo-isomorphisme $M\mathop\leftarrow\limits^\approx M'$, lequel implique à son tour une pseudo-surjection $M/J^1M\leftarrow
M'/J^1M'$. En particulier $\mathcal Char(M/J^1M)$ divise $\mathcal Char(M'/J^1M')$.

- son noyau est donc nécessairement pseudo-nul lui aussi.

$(ii)$ $M\mathop\rightarrow\limits^\approx M'$ induit $E^1M'\mathop\rightarrow\limits^\approx E^1M$. Comme $Hom_\Zp(\limind H_1(\G_n,M),\Zp)$
s'identifie à un quotient de $E^1M$, on peut raisonner comme en $(i)$.
\begin{flushright}$\square$\end{flushright}

\begin{cor} \label{critDeltak} Soit $M\in {_\L\C_{tf}}$ un module de torsion.
Alors il existe une pseudo-surjection non canonique $$\limproj f_{\Zp}(M_{\G_n})\rightarrow Hom_\Zp(\limind H_1(\G_n,M),\Zp)$$
\end{cor}
Preuve: Par \ref{IF}, on sait déjà que $\limproj f_{\Zp}(M_{\G_n})=M/I^1M\twoheadrightarrow M/F^1M$. Mais par \ref{FF}: $$M/F^1M\approx
Hom_\Zp(\limind H_1(\G_n,E^1M),\Zp)$$ Par \ref{lemmod}, on peut trouver un pseudo-isomorphisme non canonique $M\approx E^1M$, et celui-ci induit à
son tour \ref{pseudoIF}:
$$Hom_\Zp(\limind H_1(\G_n,E^1M),\Zp)\approx Hom_\Zp(\limind H_1(\G_n,M),\Zp)$$ ce qui termine la preuve.
\begin{flushright}$\square$\end{flushright}

\subsection{Stabilisation (cas $\G\simeq \Zp$)} \label{sectionstab}

Nous étudions le phénomène de stabilisation pour le défaut de (co)-descente. Celui-ci semble particulier au cas $\G\simeq \Zp$, où les complexes
$C(j_A)$ et $C(k_A)$ sont particulièrement courts. Lorsqu'elle a lieu, la stabilisation rend particulièrement agréable la description du rapport
entre descente et codescente. On suppose $\G\simeq \Zp$ dans tout \ref{sectionstab}.

\begin{defn}\label{defstab} $(i)$ On dit d'un système inductif (resp. projectif) $(A_n)$ qu'il se stabilise s'il existe un rang $n_0$ tel que les
les flèches de transition $A_n\rightarrow A_m$ (resp. $A_m\rightarrow A_n$) soient des isomorphismes pour $m\ge n\ge n_0$.

$(ii)$ On dit d'un système normique $A\in {{_{\underline\G}\C}}$ qu'il se stabilise via la restriction (resp. corestriction) si le système inductif
(resp. projectif) sous-jacent se stabilise.

$(iii)$ On dit d'un objet $A\in D^b({_{\underline\G}\C})$ qu'il se stabilise (via la restriction ou la corestriction) s'il en est ainsi pour ses
objets de cohomologie.
\end{defn}

Commençons par quelques faits généraux.

\begin{lem} \label{addstab} Soit $A\rightarrow B\rightarrow C\rightarrow A[1]$ un triangle distingué de $D^b({_{\underline\G}\C})$.
Si deux des trois sommets se stabilisent via la restriction ou la corestriction, il en est de même du troisième.
\end{lem}
Preuve: On applique le lemme des cinq. \begin{flushright}$\square$\end{flushright}

\begin{lem} \label{stabadd2} Soient $A$, $B$, $C$ des systèmes inductifs de $\Zp$-modules.

$(i)$ On suppose que les $B_n$ sont de type fini et qu'il y a une suite exacte courte $0\rightarrow A\rightarrow B\rightarrow C\rightarrow 0$. Alors
$B$ se stabilise si et seulement si $A$ et $C$ se stabilisent.

$(ii)$ On suppose que les $A_n$ et les $C_n$ sont de type fini, et qu'il y a une suite exacte courte $A\rightarrow B\rightarrow C$. Si $A$ et $C$ se
stabilisent, alors $B$ aussi.
\end{lem}

\noindent

Preuve: $(i)$ Si $A$ et $C$ se stabilisent alors $B$ aussi, clairement. Si $B$ se stabilise, alors $A_m\hookrightarrow A_n$ pour $n\ge m>>0$. Mais
alors $\limind A_n$ est une union croissante à l'intérieur du $\Zp$-module noethérien $\limind B_n$, donc se stabilise. On en déduit que $C_n$ se
stabilise aussi, par le lemme du serpent.

$(ii)$ découle immédiatement de $(i)$.
\begin{flushright}$\square$\end{flushright}

\subsubsection{Critères de stabilisation}

\noindent Pour  $A\in D^b({_{\underline\G}\C})$, on cherche à établir des critères pour:

- la stabilisation éventuelle de $C(j_A)$ via la corestriction.

- la stabilisation éventuelle de $C(k_A)$ via la restriction.

\noindent Dans la pratique, ces deux problèmes sont équivalents. En effet:

\begin{prop} \label{stabjk} Soit $A\in D^b({_{\underline\G}\C})_{tf}$, alors

$(i)$ Si $C(k_A)$ se stabilise via la restriction, alors $C(j_A)$ se stabilise via la corestriction.

$(ii)$ Soit $A\in D^b({_{\underline\G}\C})_{tf}$, vérifiant $X_\infty^*\in D^b({_\L\C})_{tf}$. Si $C(j_A)$ se stabilise via la corestriction, alors
$C(k_A)$ se stabilise via la restriction.
\end{prop}

\noindent La preuve repose sur le lemme suivant:

\begin{lem} \label{stabstab} Soit $A\in D^b({_{\underline\G}\C})$, alors

$(i)$  Si $A$ stabilise via la restriction, alors $C(j_A)$ se stabilise via la corestriction.

$(ii)$ On suppose $A\in D^b({_{\underline\G}\C})_{tf}$. Si $A$ stabilise via la corestriction, alors $C(k_A)$ stabilise via la restriction.
\end{lem}

\noindent Preuve de \ref{stabstab}: $(i)$ Considérant le triangle distingué tautologique $A\rightarrow A\Ltens_\Zp\Qp\rightarrow
A\Ltens\Qp/\Zp\rightarrow A[1]$, on se ramène tout de suite au cas où $A$ est de torsion. Par troncature et récurrence sur la longueur de $A$, on
peut supposer que $A$ est concentré en degré $0$, et quitte à remplacer $\G$ par un sous-groupe ouvert, que tous les $j_{n,m}$ sont des
isomorphismes. Mais alors $C(j_A)$ se réduit à la collection $(H^1(\G_n,A_\infty))$ (munie de sa structure naturelle de système normique, via $res$
et $cor$) placée en degré $1$. L'action de $\G$ sur $A_\infty$ étant triviale, on voit tout de suite que $cor:H^1(\G_m,A_\infty)\rightarrow
H^1(\G_n,A_\infty)$ est un isomorphisme pour tout $m\ge n$, et cela termine la preuve.

$(ii)$ A l'aide de \ref{compjk} appliqué à $A^*$, on construit facilement un isomorphisme $RHom_\Zp(C(j_{A^*}),\Zp)\simeq C(k_{Bid_\Zp(A)})$; comme
$A\simeq Bid_\Zp(A)$, le résultat souhaité se déduit de $(i)$.
\begin{flushright}$\square$\end{flushright}

Preuve de \ref{stabjk}:  $(i)$ Comme $C(j_A)=C(j_{C(k_A)})$ (cf \ref{theoprincipal} $(i)$), il suffit d'appliquer \ref{stabstab} $(i)$ à l'objet
$C(k_A)\in D^b({_{\underline\G}\C})$.

$(ii)$ Puisque $C(j_A)$ se stabilise via la corestriction, alors $Bid_\Zp C(j_A)$ aussi. Comme $C(k_A)\simeq C(k_{Bid_\Zp C(j_A)})[-1]$
(\ref{theoprincipal} $(i)$), il suffit d'appliquer \ref{stabstab} $(ii)$ à l'objet $Bid_\Zp C(j_A)\in D^b({_{\underline\G}\C})$.
\begin{flushright}$\square$\end{flushright}

\begin{rem} \label{remstab} Soit $A\in D^b({_{\underline\G}\C})_{tf}$.

$(i)$ Comme $\limind C(j_A)=0$, la stabilisation de $C(j_A)$ via la restriction ne présente aucun intérêt. Aussi, on dira désormais simplement
``$C(j_A)$ se stabilise'' au lieu de ``$C(j_A)$ se stabilise via la corestriction''. Pour des raisons similaires, on dira ``$C(k_A)$ se stabilise''
au lieu de ``$C(k_A)$ se stabilise via la restriction''.

$(ii)$ Il ne faut pas croire que la stabilisation de $C(j_A)$ ou $C(k_A)$ se produise systématiquement. Pour s'en convaincre, observer le cas d'un
système normique dont la limite inductive et la limite projective soient toutes deux triviales.
\end{rem}

Au vu de \ref{stabjk}, on s'intéresse désormais uniquement à la stabilisation de $C(k_A)$. La proposition suivante est très utile:

\begin{prop} \label{Hstab} $ $ \\
\noindent 1. Soit $A\rightarrow B\rightarrow C\rightarrow A[1]$ un triangle distingué de $D^b({_{\underline\G}\C})$. Si $C(k_A)$ et $C(k_B)$ se
stabilisent, alors il en est de même de $C(k_C)$.

\noindent 2. Soit  $A\in D^b({_{\underline\G}\C})_{tf}$ tel que $X_\infty\in D^b({_\L\C})_{tf}$. Sont alors équivalents:

$(i)$ $C(k_A)$ se stabilise.

$(ii)$ $C(k_{H^q(A)})$ se stabilise pour tout $q$.
\end{prop}

\noindent Pour la preuve, on utilisera le lemme suivant:

\begin{lem}\label{stabnoeth} Soit $A\in {{_{\underline\G}\C}}_{tf}$. Si $X_\infty\in{{_\L\C}}_{tf}$, alors $H^{-2}(C(k_A))$ se stabilise automatiquement. \end{lem}
\noindent Preuve: Le système inductif sous-jacent à $(H^{-2}(C(k_A)))\in {{_{\underline\G}\C}}$, s'identifie à la réunion croissante des
$H_1(\G_n,X_\infty)\simeq X_\infty^{\G_n}$, et celle-ci se stabilise par noethérianité. \begin{flushright}$\square$\end{flushright}

\noindent Preuve de \ref{Hstab}: 1. On applique \ref{addstab} au triangle distingué $$C(k_A)\rightarrow C(k_B)\rightarrow C(k_C)\rightarrow
C(k_A)[1]$$ déduit de celui de l'énoncé par application du foncteur $C(k_{(-)})$.

2. Par troncature et récurrence sur la longueur de $A$, l'implication $(ii)\Rightarrow (i)$ découle du point 1.

Passons à l'implication $(i)\Rightarrow (ii)$. Soit $q_0$ tel que $H^q(A)=0$ pour $q> q_0$ et considérons le triangle distingué $$C(k_{\tau_{\le
q_0-1}A})\rightarrow C(k_A)\rightarrow C(k_{H^{q_0}(A)})\rightarrow C(k_{\tau_{\le q_0-1}A})[1]$$ déduit du triangle distingué tautologique par
application du foncteur $C(k_{(-)})$. Par récurrence sur la longueur de $A$, il nous suffit de montrer que si $C(k_A)$ se stabilise, alors il en est
de même de $C(k_{\tau_{\le q_0-1}A})$ et $C(k_{H^{q_0}(A)})$. On peut toujours supposer $q_0=0$. Soit donc $A$ tel que $H^q(A)=0$ pour $q>0$, tel que
$C(k_A)$ se stabilise. D'après la remarque \ref{amplitude}, $H^q(C(k_{\tau_{\le -1}A}))=0$ (resp. $H^q(C(k_A))=0$, $H^q(C(k_{H^{0}(A)}))=0$) si
$q>-1$ (resp. $q>0$, $q\ne -2,-1,0$). En prenant le cohomologie du triangle distingué ci-dessus, on obtient donc:

- $H^0(C(k_{H^0(A)}))=H^0(C(k_A))$,  se stabilise par hypothèse.

- Pour $q<-2$ $H^q(C(k_{\tau_{\le q_0-1}A}))=H^q(C(k_A))$,  se stabilise par hypothèse.

- Une suite exacte  $$H^{-2}(C(k_{\tau_{\le q_0-1}A}))\hookrightarrow H^{-2}(C(k_A))\rightarrow H^{-2}(C(k_{H^0(A)})) \hspace{4cm}$$
$$\hspace{4cm}\rightarrow H^{-1}(C(k_{\tau_{\le q_0-1}A}))\rightarrow H^{-1}(C(k_A))\twoheadrightarrow H^{-1}(C(k_{H^0(A)}))$$
dans laquelle le second et le cinquième terme se stabilisent par hypothèse, ainsi que troisième en vertu de \ref{stabnoeth}. Par \ref{stabadd2}, on
conclut à la stabilisation des 6 termes de la suite exacte, et cela termine la preuve.
\begin{flushright}$\square$\end{flushright}

\noindent Pour appliquer \ref{Hstab} 2., on a parfois recours au lemme suivant:

\begin{lem} \label{lemtd} Soit $q_0\in \mathbb N$, et  $Z_1\hookrightarrow Z_2\rightarrow\dots \rightarrow Z_q\rightarrow \dots\twoheadrightarrow Z_{q_0}$ une suite exacte
à $q_0$ termes dans une catégorie abélienne $\mathcal A$. Alors il existe un triangle distingué (non unique) de $D^b(\mathcal A)$ dont la suite
exacte longue de cohomologie s'identifie à la précédente. \end{lem}

\noindent Preuve (indication): Par récurrence sur $q_0$. \begin{flushright}$\square$\end{flushright}

\begin{cor} \label{stabtech2}
Soit $q_0\in \mathbb N$, et  $A^1\hookrightarrow A^2\rightarrow\dots \rightarrow A^q\rightarrow \dots\twoheadrightarrow A^{q_0}$ une suite exacte à
$q_0$ termes dans ${{_{\underline\G}\C}}_{tf}$. On suppose que chaque $\L$-module $\limproj_n A^q_n$ est de type fini. Soit $i\in \{1,2,3\}$, si
$C(k_{A^q})$ se stabilise pour $q$ non congru à $i$ modulo $3$, alors $C(k_{A^q})$ se stabilise aussi pour $q$ congru à $i$ modulo $3$.
\end{cor}
Preuve: On applique \ref{Hstab} 2., 1., puis 2. à un triangle distingué fourni par le lemme précédent. Il est naturellement possible d'établir ce
corollaire directement, en procédant par découpage (il faut alors considérer $q_0-3$ suites exactes à $9$ termes, et appliquer systématiquement
\ref{amplitude}, \ref{stabnoeth} et \ref{stabadd2}).
\begin{flushright}$\square$\end{flushright}





\begin{prop} \label{stabtech} Soit $A\hookrightarrow B\twoheadrightarrow C$ une suite exacte courte dans ${{_{\underline\G}\C}}_{tf}$. On suppose que
$\limproj B_n$ est de type fini sur $\L$, et que $C(k_B)$ se stabilise. Sont alors équivalents:

$(i)$ $H^{0}C(k_A)$ se stabilise.

$(ii)$ $H^{-1}C(k_C)$ se stabilise.

$(iii)$ $C(k_A)$ et $C(k_C)$ se stabilise.
\end{prop}

\noindent Preuve de \ref{stabtech}: Comme $(iii)$ $\Rightarrow $($(i)$ et $(ii)$), il suffit de montrer que ($(i)$ ou $(ii)$) $\Rightarrow$ $(iii)$.
On note d'abord que $(iii)$ équivaut à la stabilisation de $H^0(C(k_C))$ et $H^{-1}(C(k_C))$ (\ref{stabnoeth} et \ref{Hstab} 1.).

Ensuite, on observe la suite exacte suivante de systèmes inductifs de $\Zp$-modules de types finis (NB: $\limproj A_n$, $\limproj B_n$ et $\limproj
C_n$ sont de $\L$-type fini):
$$H^{-1}(C(k_B))\rightarrow H^{-1}(C(k_C))\rightarrow H^{0}(C(k_A)) \rightarrow H^{0}(C(k_B))\rightarrow H^{0}(C(k_C))\rightarrow 0$$
\noindent Par hypothèse, le premier et le quatrième terme se stabilisent. Si de plus le second ou le troisième se stabilisent aussi, alors il en est
de même des $2$ termes restants par \ref{stabadd2}.
\begin{flushright}$\square$\end{flushright}

\subsubsection{Conséquences de la stabilisation}

Signalons qu'on fait dans ce paragraphe un usage systématique de la dualité de Poincaré; les isomorphismes en jeu sont donc fonctoriels, mais
dépendent de l'isomorphisme  $\G\simeq\Zp$ choisi.

\begin{prop} \label{corstab} Soit $A\in D^b({_{\underline\G}\C})_{tf}$ tel que $X_\infty\in D^b({_\L\C})_{tf}$.

Si $C(k_A)$ se stabilise, alors:

\noindent 1. $(i)$ $X_\infty^*\in D^b({_\L\C})_{tf}$ et $RHom_\Zp(C(j_A),\Zp)\in D^b({_{\underline\G}\C})_{tf}$.

$(ii)$ $C(j_A)$ se stabilise aussi.

\noindent 2. Les modules de cohomologie de $\Delta_j(A)$ et $\Delta_k(A)$ sont fixés par $\G_n$ pour $n$ suffisamment grand (et sont donc de torsion
sur $\L$). De plus il y a des isomorphismes fonctoriels:

$(i)$ $\Delta_k(A)\simeq RHom_\Zp(\Delta_j(A),\Zp)\simeq RHom_\Zp(R\limproj C(j_A),\Zp)$.

$(ii)$ $\Delta_j(A)\simeq RHom_\Zp(\Delta_k(A),\Zp)\simeq \limind C(k_A)[-1]$.

$(iii)$ Il existe dans $D^b(_\L\C)$ un triangle distingué fonctoriel:\vspace{-.1cm} $$R\limproj C(j_A)\rightarrow \limind C(k_A)[-1]\rightarrow
RHom_\Zp(\Qp,R\limproj C(j_A))[1]\rightarrow R\limproj C(j_A)[1]$$\vspace{-.5cm}

En particulier, $R\limproj C(j_A)\simeq \limind C(k_A)\Ltens_\Zp \Qp/\Zp[-2]$.

\noindent  3. Si de plus les $\Zp$-modules $H^q(A_n)$ sont tous de torsion (et donc finis), alors

$(i)$ les $\L$-modules de cohomologie de $X_\infty$ et $X_\infty^*$ sont de torsion.

$(ii)$ $\chi(\Delta_j(A))=\chi(\Delta_k(A))=\L$.
\end{prop}

\noindent Preuve: 1. $(i)$ Comme $A$ et $\Zp\underline\Ltens X_\infty$ sont dans $D^b({_{\underline\G}\C})_{tf}$, on voit que $C(k_A)$ et
$RHom_\Zp(C(k_A),\Zp)$ le sont aussi. Mais alors l'hypothèse de stabilisation montre que les modules de cohomologie sont de type fini sur $\Zp$, et
donc sur $\L$. La première assertion provient alors de \ref{critnoethcod}, et la seconde de \ref{remnoethjk} $(ii)$.

$(ii)$ résulte de \ref{stabjk} $(i)$.

2. Que la cohomologie de $\Delta_j(A)=R\limproj Bid_\Zp(C(j_A))$ et celle de $\Delta_k(A)=RHom_\Zp(\limind C(k_A),\Zp)[1]$ soit fixée par $\G_n$ pour
$n>>0$ est une conséquence immédiate de la stabilisation des limites en jeu.

$(i)$ D'après \ref{theoprincipal} $(iv)$ et \ref{dupoin} $(iii)$, on a $$\Delta_j(A)\simeq RHom_\L(\Delta_k(A),\L)[1]\simeq
RHom_\Zp(\Delta_k(A),\Zp)$$ d'où le premier isomorphisme de l'énoncé. Pour obtenir le second, on observe que la flèche de bidualité
$C(j_A)\rightarrow Bid_\Zp(C(j_A))$ induit dans $D^b(_\L\C)$ un carré commutatif: $$\diagram{RHom_\Zp(R\limproj Bid_\Zp(C(j_A)),\Zp)&\hfl{}{}&
RHom_\Zp(R\limproj C(j_A),\Zp)\cr \vflupcourte{}{}&&\vflupcourte{}{}\cr RHom_\Zp(Bid_\Zp\pi_n(C(j_A)),\Zp)&\hfl{}{}& RHom_\Zp(\pi_nC(j_A),\Zp)}$$
dans lequel la flèche inférieure est un isomorphisme puisque la cohomologie de $RHom_\Zp(\pi_nC(j_A),\Zp)$ est de $\Zp$-type fini. Ici, $\pi_n$ est
le foncteur de projection de la remarque \ref{rempin}, et on regarde $\pi_n C(j_A)$ comme un objet de $D^b(_\L\C)$, par abus de notation. Par
stabilisation, les deux flèches verticales sont des isomorphismes pour $n$ suffisamment grand; d'où le second isomorphisme de l'énoncé.

$(ii)$ est similaire, mais plus facile.

$(iii)$ En combinant $(i)$ et $(ii)$, on obtient un isomorphisme $$Bid_\Zp(R\limproj C(j_A))\simeq \limind C(k_A)[-1]$$ Maintenant, l'analogue dans
$D^b(_\L\C)$ de \ref{bidettate} $(ii)$ montre que $Bid_\Zp(R\limproj C(j_A))$ s'identifie à $RHom_\Zp(\Qp/\Zp,R\limproj C(j_A))[1]$; le triangle
distingué recherché provient donc du triangle distingué tautologique $\Zp\rightarrow \Qp\rightarrow \Qp/\Zp\rightarrow \Zp[1]$. L'isomorphisme de
l'énoncé suit, puisque $RHom_\Zp(\Qp,R\limproj C(j_A))\Ltens_\Zp \Qp/\Zp$ $=0$ et $R\limproj C(j_A)\simeq R\limproj C(j_A)\Ltens_\Zp\Qp/\Zp[-1]$.

3. Quitte à s'y ramener par troncature et récurrence sur la longueur de $A$, on peut toujours supposer que $A$ est concentré en degré $0$. Dans ce
cas:

$(i)$ $X_\infty$ (resp. $X_\infty^*$) est concentré en degré $0$ (resp. $1$), et la cohomologie du triangle distingué $\alpha_k(A)$ donne une suite
exacte:
$$0\rightarrow E^0(X_\infty)\rightarrow H^{0}(\Delta_k(A))\rightarrow H^1(X_\infty^*)\rightarrow E^1(X_\infty)$$ sur laquelle il apparaît que:

- $E^0(X_\infty)$ est de torsion, puisque $H^0(\Delta_k(A))$ l'est (cf 2.). D'où en fait $E^0(X_\infty)=0$, si bien que $X_\infty$ est de
$\L$-torsion.

- $H^1(X_\infty^*)$ est de $\L$-torsion puisque $H^0(\Delta_k(A))$ et $E^1(X_\infty)$ le sont (cf \ref{lemmod} $(i)$).

$(ii)$ Il suffit de calculer $\chi(\Delta_k(A))$ (cf. \ref{charjk}) Regardant $\pi_n C(k_A)$ comme un objet de $D^b(_\L\C)$, on a pour $n$
suffisamment grand un isomorphisme dans $D^b(\C_\L)$:  $\Delta_k(A)\simeq RHom_\Zp(\pi_nC(k_A),\Zp)$. Observons alors le triangle distingué
$\pi_n(C(k_A))\in Tr(D^b(_\L\C))$:
$$\Zp\Ltens_{\Zp[[\G_n]]}X_\infty \rightarrow A_n\rightarrow \pi_nC(k_A)\rightarrow \Zp\Ltens_{\Zp[[\G_n]]} X_\infty[1]$$ Comme $A_n$ est fini, on voit
que $\chi(\pi_nC(k_A))=\chi(\Zp\Ltens_{\Zp[[\G_n]]}X_\infty)$. Ce dernier idéal caractéristique est trivial puisque $X_\infty$ est de $\L$-torsion.
Finalement, $\chi(\Delta_k(A))=\L$.
\begin{flushright}$\square$\end{flushright}

En vue des applications on se propose d'expliciter la proposition précé\-dente. Rappelons que si $M$ est un $\Zp$-module, $Ext^1_\Zp(M,\Zp)$
s'identifie naturellement au dual de Pontryagin du sous-module de $\Zp$-torsion de $M$. Avec les conventions expliquées dans la section \ref{not1},
ce dernier est muni de l'action de $\L$ (à droite) par transport de structure ($(f.\lambda)(m)=f(\lambda m)$), lorsque $M$ est un $\L$-module (à
gauche).

\begin{prop} \label{exhaust}

Soit $A=(A_n)$ un système normique dans lequel les $A_n$ sont de type fini sur $\Zp$. Comme toujours, on note $X_\infty=\limproj A_n$,
$A_\infty=\limind A_n$ et $X_\infty^*=RHom_\Zp(A_\infty,\Zp)$. Pour plus de clarté, notons $$j_n:A_n\rightarrow (A_\infty)^{\G_n}   \hspace{1.5cm}
\hbox{(resp.}\,\, \, k_n:(X_\infty)_{\G_n}\rightarrow A_n \hbox{)}$$ l'application de descente (resp. codescente). De sorte que les systèmes
normiques formés par la cohomologie de $C(j_A)$ et $C(k_A)$ valent respectivement:

- Pour $C(j_A)$: $H^{-1}=(Ker\,  j_n)$, $H^0=(Coker\, j_n)$, $H^1=(H^1(\G_n,A_\infty))$, $H^2=(H^2(\G_n,A_\infty))$ et $H^q=0$ si $q\ne -1,0,1,2$.

- Pour $C(k_A)$: $H^{-2}=(H_1(\G_n,X_\infty))$, $H^{-1}=(Ker\, k_n)$, $H^0=(Coker\, k_n)$ et $H^q=0$ si $q\ne -2,-1,0$.\vspace{0.2cm}

\noindent On suppose que $Ker\, k_n$ et $Coker\, k_n$ sont de type fini sur $\Zp$, et se stabilisent. Alors: \vspace{0.1cm}

\noindent 1. $X_\infty$ est de type fini sur $\L$, et $H^q(C(j_A))$ se stabilise aussi, $\forall q$.

\noindent 2. $(i)$ $\forall q$, on a $H^q(\Delta_j(A))\simeq \limind H^{q-1}(C(k_A))$ et il y a une suite exacte courte
$$\limproj_{n,k} H^q(C(j_A))/p^k\hookrightarrow H^q(\Delta_j(A)) \twoheadrightarrow \limproj_{n,k}H^{q+1}(C(j_A))[p^k]$$
\noindent où $(-)[p^k]$ (resp. $(-)/p^k$) désigne le plus grand sous-module (resp. quotient) tué par $p^k$.

$(ii)$ $\forall q$, on a $H^q(\Delta_k(A))\simeq Ext^1_\Zp(\limproj H^{1-q}(C(j_A)),\Zp)$ et il y a une suite exacte courte
$$Ext^1_\Zp(\limind H^{-q}(C(k_A)),\Zp)\hookrightarrow H^q(\Delta_k(A))\rightarrow Hom_\Zp(\limind H^{-1-q}(C(k_A)),\Zp)$$

$(iii)$ Pour tout $q$, $H^q(\Delta_j(A)))\approx E^1H^{-q}(\Delta_k(A)))$. De plus les $\L$-modules de cohomologie de $\Delta_j(A)$ sont tous
pseudo-isomorphes à un quotient convenable de $\Zp[G_n]^a$, pour $n$ et $a$ suffisamment grand.

$(iv)$ Si les $A_n$ sont finis, alors $\prod_q \mathcal Char(H^q(\Delta_k(A)))^{(-1)^q}=\L$. \vspace{0.1cm}

\noindent 3. Les $\L$-modules $H^0(X_\infty^*)=Hom_\Zp(A_\infty,\Zp)$ et $H^1(X_\infty^*)=Ext^1_\Zp(A_\infty,\Zp)$ sont de type fini et le second est
de torsion. Il sont de plus sujets à des isomorphismes et des suites exactes:

$(i)$ $E^{q+1}(H^1(X_\infty^*))\hookrightarrow E^q(X_\infty^*)\twoheadrightarrow E^q(H^0(X_\infty^*))$.

$(ii)$ $H^{-1}(\Delta_j(A))\hookrightarrow X_\infty\rightarrow E^0(X_\infty^*)\twoheadrightarrow H^0(\Delta_j(A))$.

\noindent et $E^1(X_\infty^*)\simeq H^1(\Delta_j(A))$.

$(iii)$ $H^{-1}(\Delta_k(A))\hookrightarrow H^0(X_\infty^*)\rightarrow E^0(X_\infty)\rightarrow H^0(\Delta_k(A))\rightarrow$ \vspace{-.25cm}
\begin{flushright}$\rightarrow H^1(X_\infty^*)\rightarrow E^1(X_\infty)\twoheadrightarrow H^1(\Delta_k(A))$\end{flushright}  \vspace{-.25cm}

\noindent et  $E^2(X_\infty)\simeq H^2(\Delta_k(A))$.
\end{prop}

\noindent Preuve: Expliquons seulement les points qui ne découlent pas immédiatement de  \ref{corstab}.

2. $(i)$ Pour obtenir la suite exacte, il faut expliciter l'isomorphisme $\Delta_j(A)\simeq R\limproj Bid_\Zp(C(j_A))$. On procède en deux étapes:
d'abord \ref{bidettate} $(ii)$ donne une suite spectrale à valeurs dans ${{_{\underline\G}\C}}_{tf}$:
$$E_2^{p,q}=\limproj_k Tor_{-p}^\Zp(H^q(C(j_A)),\Z/p^k)\Rightarrow H^{p+q}(Bid_\Zp(C(j_A)))=E^{p+q}$$
Celle-ci dégénère en des suites exactes courtes, dont celles de l'énoncé se déduisent par passage à la limite projective.

$(ii)$ La stabilisation de $C(j_A)$ assure que les objets de cohomologie de $R\limproj C(j_A)$ sont de $\Zp$-torsion, et \ref{corstab} 2. $(i)$ donne
l'isomorphisme annoncé.

$(iii)$ Comme les $\L$-modules de cohomologie de $\Delta_j(A)$ sont de type fini et de torsion, le bord de la suite spectrale d'hypercohomologie (cf
\ref{theoprincipal} $(iv)$)
$$E_2^{p,q}=E^{p+1}H^{-q}(\Delta_k(A))\Rightarrow H^{p+q}(\Delta_j(A))=E^{p+q}$$ réalise pour tout $q$ un pseudo-isomorphisme $$H^q(\Delta_j(A))
\mathop\rightarrow\limits^\approx E^1H^{-q}(\Delta_k(A))$$ Le reste se déduit du fait que $H^q(\Delta_j(A))$ est fixé par $\G_n$ pour $n$
suffisamment grand.

3. Que $H^1(X_\infty^*)$ soit de $\L$-torsion se lit sur la suite exacte $(iii)$, compte-tenu de 2. $(iii)$ et \ref{lemmod} $(i)$.

$(ii)$ La suite exacte de l'énoncé, obtenue en prenant la cohomologie du triangle distingué $\alpha_j(A)$, commence \textit{a priori} par
$$0\rightarrow E^{-1}(X_\infty^*)\rightarrow H^{-1}(\Delta_j(A))\rightarrow H^0(X_\infty)\rightarrow \dots$$
Maintenant $E^{-1}(X_\infty^*)=E^0(H^{1}(X_\infty^*))$ est trivial puisque $H^{1}(X_\infty^*)$ est de $\L$-torsion.
\begin{flushright}$\square$\end{flushright}

\begin{cor} \label{corexhaust} Soit $A=(A_n)\in {{_{\underline\G}\C}}_{tf}$. Si tous les $A_n$ sont de $\Zp$-torsion, et que $C(k_A)$ se stabilise,
alors:

$(i)$ les $\L$-modules $X_\infty=\limproj A_n$ et $H^1(X_\infty^*)=Hom_\Zp(\limind A_n,\Qp/\Zp)$ sont noethériens de torsion.

$(ii)$ Non canoniquement, les quatre modules suivants sont pseudo-iso\-morphes: $H^{-1}(\Delta_j(A)), H^0(\Delta_j(A)), H^0(\Delta_k(A))$ et
$H^1(\Delta_k(A))$.
\end{cor}

\noindent Preuve: $(i)$ Comme les $A_n$ sont de $\Zp$-torsion, c'est que $X_\infty^*$ est concentré en degré $1$. Maintenant, $H^1(X_\infty^*)$ est
de $\L$-torsion (\ref{exhaust} 3.) et les modules de cohomologie de $\Delta_j(A)$ aussi (2. $(iii)$); il en est donc de même de $X_\infty$.

$(ii)$ D'après \ref{exhaust} 2. $(iii)$, il suffit de vérifier que l'on a \break $\chi(H^1(\Delta_k(A)))=\chi(H^0(\Delta_k(A)))$. Mais cela découle
de \ref{exhaust} 2. $(iv)$, puisque d'après \ref{exhaust} 3. $(iii)$ $H^{-1}(\Delta_k(A))\subset H^0(X_\infty^*)=0$, $H^2(\Delta_k(A))\simeq
E^2(X_\infty)\approx 0$ (\ref{lemmod} $(i)$).
\begin{flushright}$\square$\end{flushright}

Mentionnons encore un corollaire de la description \ref{corstab} 2.

\begin{cor} \label{troncaturejk} Soit $A\in {_{\underline\G}\C}_{tf}$ et supposons, comme en \ref{corstab}, que
les systèmes inductifs $Ker\, k_n$ et $Coker\, k_n$ se stabilisent. Notons $$\phi:R\limproj C(j_A)\rightarrow \limind C(k_A)[-1]$$ le morphisme
fonctoriel figurant dans le triangle \ref{corstab} 2. $(iii)$. Alors:

\noindent 1.  $(i)$ $Ker\,  H^q(\phi)$ est $p$-divisible.

$(ii)$ $Coker\, H^q(\phi)$ est $\Zp$-libre.

\noindent 2. Considérons un triangle distingué $A\rightarrow A'\rightarrow A''\rightarrow A[1]$, où $A'$ et $A''$ vérifient les mêmes hypothèses que
$A$. Dans ces conditions, pour $q\in \Z$ fixé:

$(i)$ $\phi$ induit un isomorphisme de suites exactes courtes \vspace{-0.7cm}
{\small$$\diagram{\limproj_{n,k}H^q(C(j_A))/p^k&\hflcourte{}{}&\limproj_{n,k}H^q(C(j_{A'}))/p^k &\hflcourte{}{}&\limproj_{n,k}H^q(C(j_{A''}))/p^k\cr
\vflcourte{\wr}{}&&\vflcourte{\wr}{}&&\vflcourte{\wr}{}&& \cr t_\Zp \limind H^{q-1}(C(k_{A})) &\hflcourte{}{}&t_\Zp \limind H^{q-1}(
C(k_{A'}))&\hflcourte{}{}&t_\Zp \limind H^{q-1}(C(k_{A''}))}$$}\vspace{-0.7cm}

$(ii)$ Si $\limproj H^q(C(j_{A''}))$ est fini, alors le diagramme ci-dessus se prolonge en un isomorphisme de suites exactes longues:

- à gauche à l'aide de l'isomorphisme $\Delta_j(-)\simeq\limind C(k_{(-)})[1]$: \vspace{-0.7cm}{\small$$\diagram{\dots&\hflcourte{}{}&
H^{q-1}(\Delta_j(A'))&\hflcourte{}{}&H^{q-1}(\Delta_j(A''))&\hflcourte{}{}&\limproj_{n,k}H^q(C(j_A))/p^k&\hflcourte{}{}\cr
&&\vflcourte{\wr}{}&&\vflcourte{\wr}{}&&\vflcourte{\wr}{}&& \cr \dots&\hflcourte{}{}&\limind H^{q-2}(C(k_{A'}))&\hflcourte{}{}&\limind H^{q-2}(
C(k_{A''})) &\hflcourte{}{}&t_\Zp \limind H^{q-1}(C(k_{A}))&\hflcourte{}{}}$$}\vspace{-0.7cm}

- à droite à l'aide de $R\limproj C(j_{(-)})\simeq \limind C(k_{(-)})\Ltens_\Zp \Qp/\Zp[-2]$:\vspace{-0.7cm}
{\small$$\diagram{\limproj_{n,k}H^q(C(j_{A''}))/p^k&\hflcourte{}{}&\limproj H^{q+1} C(j_A)&\hflcourte{}{}&\limproj
H^{q+1}C(j_{A'})&\hflcourte{}{}&\dots \cr \vflcourte{\wr}{}&&\vflcourte{\wr}{}&&\vflcourte{\wr}{}&& \cr t_\Zp \limind
H^{q-1}(C(k_{A''}))&\hflcourte{}{}&Tor_{q-1}^\Zp(\limind C(k_A),\Qp/\Zp)&\hflcourte{}{}&Tor_{q-1}^\Zp(\limind
C(k_{A'}),\Qp/\Zp)&\hflcourte{}{}&\dots}$$}\vspace{-0.7cm}
\end{cor} \noindent Preuve: 1. est clair, puisque les modules de cohomologie du cône de $\phi$ sont uniquement $p$-divisibles.

2. $(i)$ Comme $H^q(C(j_A))$ est de $\Zp$-torsion et que $H^q(C(k_A))$ est de type fini, on peut préciser 1.: $Ker\, H^q(\phi)$ est le sous groupe
$p$-divisible maximal de $H^q(R\limproj C(j_A))$ (ie. le noyau de $\limproj H^q(C(j_A))\rightarrow \limproj_{n,k} H^q(C(j_A))/p^k$) et $Im\,
\phi=t_\Zp \limind H^{q-1}(C(k_A))$. D'où l'isomorphisme de suites exactes de l'énoncé.

$(ii)$ - Au vu de \ref{exhaust} 2. $(i)$, \ref{corstab} 2. $(ii)$ et de la définition de $\phi$, le diagramme commutatif suivant se produit en toute
généralité: \vspace{-1cm}

{\small$$\diagram{\dots&\hflcourte{}{}&
H^{q-1}(\Delta_j(A''))&\hflcourte{}{}&H^{q}(\Delta_j(A))&\hflrevcourte{\supset}{}&\limproj_{n,k}H^q(C(j_A))/p^k&\hflcourte{}{}\cr
&&\vflcourte{\wr}{}&&\vflcourte{\wr}{}&&\vflcourte{\wr}{}&& \cr \dots&\hflcourte{}{}&\limind H^{q-2}(C(k_{A''}))&\hflcourte{}{}&\limind H^{q-1}(
C(k_{A})) &\hflrevcourte{\supset}{}&t_\Zp \limind H^{q-1}(C(k_{A'}))&\hflcourte{}{}}$$} \vspace{-0.8cm}

Puisque par hypothèse $H^{q-1}(\Delta_j(A''))$  est fini (\ref{exhaust} 2. $(i)$), l'image des flèches de gauche est incluse dans celle des flèches
de droite. Le diagramme de l'énoncé est donc bien défini et ses lignes sont automatiquement exactes.

- De manière analogue, on peut prolonger le diagramme à droite comme indiqué dans l'énoncé, à condition que $\limproj H^q(C(j_{A''}))\rightarrow
\limproj H^{q+1}(C(j_A))$ se factorise par le quotient fini maximal $\limproj_{n,k}H^q(C(j_A))/p^k$. C'est ici le cas, puisque par hypothèse,
$\limproj H^q(C(j_{A''}))$ est lui-même fini.
\begin{flushright}$\square$\end{flushright}

\section{Applications}

\subsection{Cohomologie galoisienne continue}

\subsubsection{Généralités}

Soit $G$ un groupe profini, muni d'une surjection fixée $G\twoheadrightarrow \G$. Pour chaque sous-groupe ouvert $\G_n$ de $\G$, notons $H_n$ l'image
réciproque de $\G_n$ dans $G$, et $H_\infty=\cap H_n$. Ainsi, $G_n=\G/\G_n=G/H_n$, et $\G=G/H_\infty$. On note $_G\C$ la catégorie des $\Zp$-modules
munis d'une action discrète de $G$, et $Rep_\Zp(G)$ la catégorie des $\Zp$-représentations continues de $G$ dont le $\Zp$-module sous-jacent est
$\Zp$-libre de rang fini. La définition \cite{J1} pour la cohomologie continue se prête à la variante suivante:

\begin{defn} \label{defcohcont}Soit $*=b$ si $cd_p G<\infty$, et $*=+$ sinon.

$(i)$ Soit $\underline\G(G,-):{_G\C}\rightarrow {{_{\underline\G}\C}}$ le foncteur qui à $M\in {_G\C}$ associe le système normique $(M^{H_n})$. On
note $R\underline\G(G,-):D^*(_G\C)\rightarrow D^*({_{\underline\G}\C})$ le dérivé droit de $\underline\G(G,-)$.

$(ii)$ Soit $_G\C^{\N^\circ}$ la catégorie des systèmes projectifs d'objets de $_G\C$ indexés sur $\N$. On note
$R\underline\G(G,-):D^*(_G\C^{\N^\circ})\rightarrow D^*({_{\underline\G}\C})$ le foncteur composé \vspace{-0.8cm}
$$\diagram{{_G\C}^{\N^\circ}&\hfl{R\underline \G(G,-)}{}& {_{\underline\G}\C}^{\N^\circ}&\hfl{R{\mathop{\lim}\limits_{\leftarrow}}_\N}{}&
{_{\underline\G}\C}}$$ \vspace{-0.8cm}

$(iii)$ Si $T$ est une $\Zp$-représentation continue de $G$, on note $R\underline\G(G,T):=R\underline\G(G,(T/p^k))$, où $(T/p^k)$ désigne le système
projectif évident.
\end{defn}

\begin{rem} \label{remcohcont} $(i)$ Soit $A\in D^*(_G\C)$ et $(A)\in D^*({_G\C^{\N^{\circ}}})$ le système projectif de valeur constante $A$. On a alors
$R\underline\G(G,A)=R\underline\G(G,(A))$. Plus généralement, si $A\in {_G\C^{\N^{\circ}}}$ est un système projectif qui se stabilise, sa cohomologie
au sens $(ii)$ coïncide avec la cohomologie de sa limite, au sens $(i)$.

$(ii)$ Si l'image de $G$ dans $Aut_\Zp(T)$ est finie, on peut donner un sens à $R\underline\G(G,T)$ soit par $(i)$ soit $(iii)$. On convient
systématiquent du sens $(iii)$ lorsque $G$ est un groupe de Galois du type $G_F$ ou $G_{S,F}$ (cf \ref{ssjannsen}).

$(iii)$ Soit $K^*$ un foncteur résolvant, au sens de \cite{Se2}, Annexe 2. def 1.4, et notons $C_k^\bullet=K^*(T/p^k)$. Alors la collection des
résolutions $T/p^k\rightarrow C_k^\bullet$ est organisée en système projectif, et donne donc une résolution dans ${_G\C^{\N^{\circ}}}$
$(T/p^k)\rightarrow C^\bullet$, où $C^\bullet=(C_k^\bullet)\in Kom^+({_G\C})^{\N^{\circ}}=Kom^+({_G\C^{\N^{\circ}}})$. \\ Celle-ci vérifie:

$(*)$ Les $G$-modules $C_k^q$ sont $G$-cohomologiquement triviaux.

$(**)$ Les systèmes projectifs $(\underline\G (C_k^q))\in ({_{\underline\G}\C})^{\N^\circ}$ sont flasques.

\noindent Aussi, $R\underline\G(G,T)$ est représenté par $\limproj_k \underline\G(G,C_k^\bullet)\in Kom^+({_{\underline\G}\C^{\N^{\circ}}})$.

$(iv)$ Comme expliqué dans \cite{J1}, le propriétés $(*)$ et $(**)$ ci-dessus sont vérifiées par toute résolution injective $(T/p^k)\rightarrow
I^\bullet$, $I^\bullet\in Kom^+({_G\C})^{\N^{\circ}}$. On en déduit que $R\limproj_k\circ R\underline\G$ est aussi le dérivé droit du foncteur
composé $\limproj\circ \underline\G$.

$(v)$ Bien sûr, on retrouve la cohomologie continue de $H_n$ en appliquant le foncteur $\pi_n$ de \ref{rempin}: $\pi_nR\underline\G(G,T)\simeq
R\G(H_n,T)$.

$(vi)$ En considérant le cône de la multiplication par $p^t$: $(T/p^k)\rightarrow (T/p^k)$, on voit tout de suite qu'il y a isomorphisme naturel
$R\underline\G(G,T)\Ltens_\Zp\Z/p^t\simeq R\underline\G(G,T/p^t)$; mieux, la collection de ces isomorphismes pour $t$ variable se relève en un
isomorphisme de $D^+({_{\underline\G}\C^{\N}})$, d'où en passant à la limite inductive: $$R\underline\G(G,T)\Ltens_\Zp\Qp/\Zp\simeq
R\underline\G(G,T\otimes\Qp/\Zp)\hspace{1.5cm}\hbox{dans $D^+({_{\underline\G}\C})$}$$
\end{rem}

On renvoie à \cite{J1} pour les propriétés générales de la cohomologie continue, et la comparaison avec la définition de \cite{Ta} (par cochaînes
continues). On se contente de détailler le résultat suivant, qui reflète essentiellement la suite spectrale de Hochschild-Serre en cohomologie
continue pour la collection des extensions de groupes $H_n\hookrightarrow G\twoheadrightarrow G_n$:

\begin{prop} \label{desccohcont} Soit $T\in Rep_\Zp(G)$, et notons $A:=R\underline\G(G,T)\in D^*({_{\underline\G}\C})$. Alors $C(j_A)=0$. \end{prop}

\noindent Preuve: On note $k$ l'indice des systèmes projectifs sur $\N$ et $n$ ou $m$ celui des systèmes projectifs sur $N$, sous-jacents aux objets
de ${{_{\underline\G}\C}}$. Il s'agit de vérifier que la flèche de restriction $res_n:R\G(H_n,T)\rightarrow R\G(\G_n,\limind_m R\underline\G(G,T))$
est un isomorphisme (NB: avec les notations expliquées plus haut, la cohomologie de $\limind_m R\underline\G(G,T)$ vaut $\limind_m
H^\bullet(H_m,T)$).

Soit $C^\bullet=(C^\bullet_k)\in Kom^+(_G\C^{\N^{\circ}})$ un complexe qui possède les propriétés $(*)$ et $(**)$ (cf \ref{remcohcont}). Les objets
$R\G(H_n,T)$ et $\limind_m R\underline\G(G,T)$ sont alors représen\-tés respectivement par $\limproj_k{C_k^\bullet}^{H_n}$ et
$\limind_m\limproj_k{C_k^\bullet}^{H_m}$, et il nous suffit de montrer que les objets de ce dernier complexe sont $\G_n$-acycliques. En effet, si
c'est le cas la flèche $res_n$ est alors représentée (cf \ref{remdeftr}) par la flèche suivante: $$\limproj_k{C_k^\bullet}^{H_n}\rightarrow
(\limind_m\limproj_k{C_k^\bullet}^{H_m})^{\G_n}$$ qui est visiblement un isomorphisme.

Calculons donc, pour $p\ge 1$:
$$H^p(\G_n,\limind_m\limproj_k{C_k^q}^{H_m})\simeq\limind_m H^p(\G_n/\G_m,\limproj_k{C_k^q}^{H_m})$$
$$\hbox{par $(**)$}\hspace{3cm} \simeq\limind_m
H^p(\G_n/\G_m,R\limproj_k{C_k^q}^{H_m})$$ Maintenant, $$H^p(\G_n/\G_m,R\limproj_k{C_k^q}^{H_m})=Ext^p_{\Zp[\G_n/\G_m]}(\Zp,R\limproj_k{C_k^q}^{H_m})$$
$$\hbox{par \ref{remlimproj}}\hspace{3cm} \simeq R^p\limproj_kRHom_{\Zp[\G_n/\G_m]}(\Zp,{C_k^q}^{H_m})$$
Mais d'après $(*)$, ${C_k^q}^{H_m}$ est $\G_n/\G_m$-acyclique, si bien que
$$R^p\limproj_kRHom_{\Zp[\G_n/\G_m]}(\Zp,{C_k^q}^{H_m})\simeq R^p\limproj_k ({C_k^q}^{H_m})^{\G_n/\G_m}$$ $$\hspace{5cm}\simeq R^p\limproj_k
{C_k^q}^{H_n}$$ $$\hbox{par $(**)$}\hspace{1.5cm} \simeq 0$$ Ainsi, $H^p(\G_n,\limind_m\limproj_k{C_k^q}^{H_m})\simeq 0$ pour $p\ge 1$, et cela
termine la preuve.
\begin{flushright}$\square$\end{flushright}

\subsubsection{Un calcul d'adjoint à la Jannsen} \label{ssjannsen}

Soit $F$ un corps de nombres, et $F_\infty/F$ une extension galoisienne de groupe $\G$. Soit $S$ un ensemble de places de $F$ (fini ou non) contenant
celles qui se ramifient dans $F_\infty/F$. On note $G_S=G_{S,F}$ le groupe de Galois de la clôture $S$-ramifiée de $F$, et $G_{S,F_n}$ le sous-groupe
qui fixe $F_n$. Nous sommes donc dans la situation du paragraphe précédent avec $G=G_S\twoheadrightarrow \G$, $H_n=G_{S,F_n}$.

\begin{thm}\label{thmcohcont} Si $cd_p G_S<\infty$, alors:

$(i)$ $R\limproj Bid_\Zp(R\underline\G(G_{S,F_n},T))\simeq RHom_\L(RHom_\Zp(\limind  R\underline\G(G_{S,F},T),\Zp),\L)$

$(ii)$ Si $S$ est fini, $R\underline\G(G_{S,F_n},T)\simeq Bid_\Zp(R\underline\G(G_{S,F_n},T))$.

\end{thm}
Preuve: $(i)$ résulte de \ref{remtheoprincipal}, compte-tenu de \ref{desccohcont}.

$(ii)$ Lorsque $S$ est fini, on a $R\underline\G(G_S,T)\in D^b({_{\underline\G}\C})_{tf}$ (cf \cite{Ta} cor. prop. 2.1).
\begin{flushright}$\square$\end{flushright}

\begin{cor} \label{ssjansen} Si $S$ est fini et $cd_p G_S<\infty$, il
y a une suite spectrale à valeurs dans ${_\L\C}$: $$Ext_\L^p(Hom_\Zp(H^q(G_{S,F_\infty},T\otimes \Qp/\Zp),\Qp/\Zp),\L)\Rightarrow \limproj
H^{p+q}(G_{S,F_n},T)$$
\end{cor}
Preuve: Notons qu'il y a un isomorphisme naturel $$RHom_\Zp(\limind  R\underline\G(G_{S,F},T),\Zp))\simeq RHom_\Zp(\limind
R\underline\G(G_{S,F},T\otimes\Qp/\Zp),\Qp/\Zp))$$ La suite spectrale d'hypercohomologie du foncteur $RHom_\L(-,\L)$, à coefficients dans
$Hom_\Zp(H^q(G_{S,F_\infty},T\otimes \Qp/\Zp),\Qp/\Zp)$ donne donc le résultat, compte-tenu de \ref{thmcohcont}.

\begin{rem} Cette suite spectrale fait l'objet d'une note non publiée de Jannsen \cite{J2}. Elle y est obtenue ``à la main''
indépendamment de nos hypothèses (selon lesquelles $\L$ est noethérien, $cd_p\G<\infty$ et $cd_p G_S<\infty$).
\end{rem}

\subsection{Le groupe des $(p)$-classes dans les $\Zp^d$-extensions} \label{ClZpd}

On conserve les notations du paragraphe précédent, mais on suppose désormais $\G\simeq \Zp^d$, $S=S_p\cup S_\infty$ et $p\ne 2$ si $F$ possède une
place réelle (ces conditions assurent $cd_p G_S\le 2$). Comme première illustration des méthodes de cet article, on se propose d'établir le théorème
\ref{propcl} ci-dessous et ses corollaires. On peut voir cette étude comme un analogue partiel de \cite{LFMN} pour les $\Zp^d$-extensions.

Les résultats qui suivent (\ref{propcl}, \ref{corNV1} et \ref{corNV2}) ont été obtenus en collaboration avec J. \nekovar. Leur preuve occupe toute la
présente section. \vspace{0.1cm}

La taille des groupes de décomposition de $\G$ joue un rôle important dans notre étude. Notons $S^{dec}=\{v\in S_p(F),\G_v=0\}$ et considérons
l'hypothèse suivante:

$(Dec_a)$ Pour tout $v\in S_p(F)$, $rg_\Zp \G_v\ge a$. \vspace{0.1cm}

\noindent (ainsi $(Dec_1)\Leftrightarrow S^{dec}=\emptyset$). Soit encore $Z:=Ker\left(\oplus_{v\in S_p(F),\G_v\ne 0}\Zp[[\G/\G_v]]\rightarrow
\Zp\right)$.

\begin{thm}\label{propcl} Notons $\mathcal Cl'_n:=\mathcal Cl(\mathcal O_{F_n}[{1\over p}])\otimes\Zp$.

$(i)$ Les $\L$-modules $\limproj \C l_n'$ et $Hom_\Zp(\limind \mathcal Cl'_n,\Qp/\Zp)$ sont noethériens de torsion. De plus, il existe une morphisme
naturel naturel (défini en \ref{theoprincipal}):
$$\alpha_k:Hom_\Zp(\limind \mathcal Cl'_n,\Qp/\Zp)\rightarrow E^1(\limproj \C l_n')$$ \nopagebreak
 et un pseudo-isomorphisme (non canonique) $Ker\alpha_k\approx Coker\alpha_k$.

$(ii)$ $Ker\, \alpha_k\hookrightarrow E^1(Z)$. En particulier, $\alpha_k$ est un pseudo-isomorphisme injectif si l'hypothèse $(Dec_2)$ est vérifiée.
\end{thm}

\begin{cor}\label{corNV1} ~~

$(i)$ $\mathcal Char(\limproj\C l'_n)=\C har(Hom_\Zp(\limind \mathcal Cl'_n,\Qp/\Zp))$. De plus:
$$\limproj \C l'_n\approx 0\Leftrightarrow \limind \C l'_n=0$$

$(ii)$ Désignons par $(\limproj \C l'_n)^0$ le sous-module pseudo-nul maximal de $\limproj \C l'_n$ et considérons l'application de capitulation
$j_n:\mathcal Cl'_n\rightarrow (\limind_m \mathcal Cl'_m)^{\G_n}$. Sous $(Dec_2)$, il y a une suite exacte canonique: $$\limproj_n Ker\,
j_n\hookrightarrow (\limproj \C l'_n)^0\twoheadrightarrow \limproj_{n,k} Coker\, j_n[p^k]$$
\end{cor}

\begin{cor}\label{corNV2} (comp. \cite{Ne} 9.4.11). Notons $E_n':=\mathcal O_{F_n}[{1\over p}]^\times\otimes\Zp$. Alors: \vspace{.1cm}

$(i)$ $t_\L Hom_\Zp(\limind E'_n\otimes \Qp/\Zp,\Qp/\Zp)\approx Z$, non canoniquement.

$(ii)$ Sous $(Dec_3)$, il y a un isomorphisme naturel: $$Coker\alpha_k\simeq t_\L Hom_\Zp(\limind E'_n\otimes \Qp/\Zp,\Qp/\Zp)$$
 \end{cor}

Principe de la preuve de \ref{propcl}: Considérons les systèmes normiques naturels $\mathcal Cl'=(\C l'_n)\in {_{\underline\G}\C}_{tf}$ et
$E'=(E_n')\in {_{\underline\G}\C}_{tf}$. Prouver \ref{propcl} revient à expliciter la suite exacte longue de cohomologie de $\alpha_k(\C l')$. Pour
cela, nous utiliserons les systèmes normiques auxiliaires suivants:

- $R\G_S:=R\underline\G(G_S,\Zp(1))\in D^b({{_{\underline\G}\C}})_{tf}$ (cf \ref{defcohcont}).

- $H^q_S:=H^q(R\underline\G_S)=(H^q(G_{S,F_n},\Zp(1)))\in {{_{\underline\G}\C_{tf}}}$, $q=1,2$.

- $Br':=(Br_n')\in {{_{\underline\G}\C_{tf}}}$, où  $Br_n'=\limproj_k H^2(G_{S,F_n},\Qp/\Zp(1))[p^k]$ le module de Tate du groupe de Brauer de
$\mathcal O_{F_n}[{1\over p}]$.

- $\underline \L_v=\Zp[[\G/\G_v]]\underline\otimes \Zp=(\Zp[[\G/\G_v]]_{\G_n})\in {{_{\underline\G}\C_{tf}}}$, où $\G_v$ est le groupe de
décomposition de $\G$ en $v$. Concrètement, $k_{m,n}:\Zp[[\G/\G_v]]_{\G_m}\rightarrow \Zp[[\G/\G_v]]_{\G_n}$ (resp.
$j_{n,m}:\Zp[[\G/\G_v]]_{\G_n}\rightarrow \Zp[[\G/\G_v]]_{\G_m}$) est la flèche induite par l'identité (resp. la multiplication par
$\sum_{g\in\G_n/\G_m}g$).

- $\underline \Z_p=\Zp\underline\otimes\Zp=(\Zp)\in {{_{\underline\G}\C_{tf}}}$. Concrètement, $k_{m,n}=id$ et $j_{n,m}=[\G_n:\G_m]$.

 \vspace{0.1cm}

\noindent Résumons en un lemme les relations qui existent entre ces différents objets:

\begin{lem} \label{tdtriv} Il y a dans $D^b({_{\underline\G}\C})_{tf}$ des triangles distingués:

$(i)$ $H^1_S[-1]\rightarrow R\G_S\rightarrow H^2_S[-2]\rightarrow H^1_S$.

$(ii)$ $\C l'\rightarrow H^2_S\rightarrow Br'\rightarrow \C l'[1]$. De plus, $H^1_S\simeq E'$.

$(iii)$ $Br'\rightarrow \oplus_{v\mid p} \underline \L_v\rightarrow \underline\Z_p\rightarrow Br'[1]$. \end{lem} \noindent Preuve (esquisse): La
suite exacte de Kummer (resp. le principe de Hasse) donne lieu à un isomorphisme et une suite exacte (resp. une suite exacte):
$$E_n'\simeq H^1(G_{S,F_n},\Zp(1))\hspace{1.5cm} \mathcal Cl'_n\hookrightarrow H^2(G_{S,F_n},\Zp(1))\twoheadrightarrow Br_n'$$ $$Br_n'\hookrightarrow
\oplus_{v\mid p} \Zp[[G_n/G_{n,v}]]\twoheadrightarrow \Zp$$ compatibles à la restriction et à la corestriction.
\begin{flushright}$\square$\end{flushright}

Pour $A\in {_{\underline\G}\C}$ ou $D^b({_{\underline\G}\C})$, on note $A_\infty(A),$ $X_\infty(A),$ $A_\infty^*(A),$ $X_\infty^*(A)$ les objets
notés simplement $A_\infty,$ $X_\infty,$ $A_\infty^*,$ $X_\infty^*$ dans le paragraphe \ref{parcod}. Commençons par deux lemmes.

\begin{lem} \label{cledesc} $(i)$ $C(j_{R\G_S})=0$, et $X_\infty^*(R\G_S)\in D^b(\C_\L)_{tf}$.

$(ii)$ $C(k_{R\G_S})=0$, et $X_\infty(R\G_S)\in D^b({_\L\C})_{tf}$.

$(iii)$ $\Delta_k(R\G_S)=\Delta_j(R\G_S)=0$. \end{lem} Preuve: Tout provient de la proposition \ref{desccohcont}, moyennant \ref{remnoethjk} et
\ref{corcjk}. \begin{flushright}$\square$\end{flushright}

\begin{lem} \label{cleBr} Soit $Br'_{non-dec}:=Ker(\mathop\oplus\limits_{v\notin S^{dec}} \underline\L_v\rightarrow \underline\Z_p)$,
de sorte que le module $Z$ de $\ref{corNV2}$ s'identifie naturellement à $X_\infty(Br'_{non-dec})$. Alors:

$(i)$ Il y a dans $D^b({_{\underline\G}\C})_{tf}$ un triangle distingué naturel: $$Br'_{non-dec}\rightarrow Br'\rightarrow \oplus_{v\in S^{dec}}
\underline\L\rightarrow Br'_{non-dec}[1]$$ lequel induit $Z=X_\infty(Br'_{non-dec})\simeq t_\L X_\infty(Br')$ et $f_\L(X_\infty(Br'))\simeq
\oplus_{v\in S^{dec}}\L$.

$(ii)$ On a $X_\infty^*(Br'_{non-dec})=0$, $\Delta_j(Br'_{non-dec})\simeq \Delta_j(Br')\simeq Z[1]$ et $\Delta_k(Br')\simeq
\Delta_k(Br'_{non-dec})\simeq RHom_\L(Z,\L)\simeq\tau_{\ge 1}RHom_\L(X_\infty(Br'),\L)$.

$(iii)$ Sous $(Dec_a)$, on a $Ext^q(X_\infty(Br'),\L)=0$ pour $q\le a-1$. En particulier, $X_\infty(Br')\approx 0 \Leftrightarrow (Dec_2)$.
\end{lem}
Preuve: $(i)$ Comme $S^{dec}\varsubsetneq S$, la projection naturelle $Br'\rightarrow\oplus_{v\in S^{dec}} \underline\L$ est surjective, d'où le
triangle distingué de l'énoncé. Le reste suit.

$(ii)$ Ayant remarqué que:

- si $v\notin S^{dec}$, alors $A_\infty(\underline\L_v)$ est uniquement $p$-divisible, et donc
$X_\infty^*(\underline\L_v)=RHom_\Zp(A_\infty(\underline\L_v),\Zp)=0$.

- $X_\infty^*(\underline\Z_p)=0$.

\noindent on voit que $X_\infty(Br'_{non-dec})=0$. Comme par ailleurs $\Delta_k(\underline\L)=\Delta_j(\underline\L)=0$, le résultats de l'énoncé se
déduisent facilement de $(i)$, en appliquant les foncteurs $\alpha_j$ et $\alpha_k$ aux objets et aux flèches du triangle distingué de $(i)$.


$(iii)$ Comme $\Zp[[\G/\G_v]]=\L\otimes_{\Zp[[\G_v]]}\Zp\in {_\L\C}$, on a dans $D^b(\C_\L)$: $$RHom_\L(\Zp[[\G/\G_v]],\L)\simeq
RHom_{\Zp[[\G_v]]}(\Zp,\L)\simeq RHom_{\Zp[[\G_v]]}(\Zp,\Zp[[\G_v]])\Ltens_{\Zp[[\G_v]]}\L$$ Maintenant, $RHom_{\Zp[[\G_v]]}(\Zp,\Zp[[\G_v]])\simeq
\Zp[-d_v]$ et $RHom_\L(\Zp,\L)\simeq \Zp[-d]$ où $d_v$ (resp. $d$) désigne le rang de $\G_v$ (resp. $\G$) (cf \ref{dupoin} $(iii)$). Aussi
\ref{tdtriv} $(iii)$ donne-t-il un triangle distingué de $D^b(\C_\L)$:
$$\Zp[-d]\rightarrow \oplus_{v\mid p} \Zp[[\G/\G_v]][-d_v]\rightarrow RHom_\L(X_\infty(Br'),\L)\rightarrow \Zp[1-d]$$ D'où la première assertion de
l'énoncé. L'assertion concernant la pseudo-nullité résulte de \ref{lemmod} $(iv)$.
\begin{flushright}$\square$\end{flushright}

Nous sommes maintenant en mesure d'expliciter le théorème \ref{theoprincipal} appliqué au système normique $H^2_S$. Pour plus de lisibilité, et en
vue de références futures, on adopte les notations suivantes:

- $H^q_{Iw}(F_\infty,\Zp(1)):=\limproj H^q(G_{S,F_n},\Zp(1))=H^q(X_\infty(R\G_S))\simeq X_\infty(H^q_S)$.

- $\mathcal H^q(F_\infty,\Zp(1)):= \limind H^q(G_{S,F_n},\Zp(1))=H^q(A_\infty(R\G_S))\simeq A_\infty(H^q_S)$.

- $j_n^q:H^q(G_{S,F_n},\Zp(1))\rightarrow \mathcal H^q(F_\infty,\Zp(1))^{\G_n}$, de sorte que $H^0(j_{H^q_S})=(j_n^q)$.

- $k_n^q:H^q_{Iw}(F_\infty,\Zp(1))_{\G_n}\rightarrow H^q(G_{S,F_n},\Zp(1))$, si bien que $H^0(k_{H^q_S})=(k_n^q)$.

\begin{prop}\label{propZpi} $ $ \\
\noindent 1. Les $\L$-modules $H^q_{Iw}(F_\infty,\Zp(1))$, $Ext_\Zp^p(\mathcal H^q(F_\infty,\Zp(1)),\Zp)$, $p=0,1$, $q=1,2$, sont tous de type fini.

\noindent 2. $(i)$  $C(j_{H^2_S})\simeq C(j_{H^1_S})[2]$ (resp.  $C(k_{H^2_S})\simeq C(k_{H^1_S})[2]$) est acyclique hors de $[-1,d-1]$ (resp.
$[-d-1,-2]$).

\noindent Explicitement:

$(ii)$ $Ker\, j_n^1=Coker\, j_n^1=0$, $Ker\,  k_n^2=Coker\, k_n^2=0$.

$(iii)$ $Ker\, j_n^2=H^1(\G_n,\mathcal H^1(F_\infty,\Zp(1)))$, $Coker\, j_n^2=H^2(\G_n,\mathcal  H^1(F_\infty,\Zp(1)))$,

$H_2(\G_n,H^2_{Iw}(F_\infty,\Zp(1)))=Ker\, k_n^1$, $H_1(\G_n,H^2_{Iw}(F_\infty,\Zp(1)))=Coker\, k_n^1$.

$(iv)$ Pour $q\ge 1$ on a, $$H^q(\G_n,\mathcal H^2(F_\infty,\Zp(1)))=H^{q+2}(\G_n,\mathcal H^1(F_\infty,\Zp(1)))$$
$$H_{q+2}(\G_n,H^2_{Iw}(F_\infty,\Zp(1)))=H_{q}(\G_n,H^1_{Iw}(F_\infty,\Zp(1)))$$

\noindent 3. $(i)$ $\Delta_k(H^2_S)=\Delta_k(H^1_S)[-2]$ est acyclique hors de $[1,d+1]$.

$(ii)$ Les modules de cohomologie de $\Delta_k(H^2_S)$ sont noethériens de torsion. De plus:

- $H^q(\Delta_k(H^2_S))$ est pseudo-nul si $q\ne 1$.

- $H^1(\Delta_k(H^2_S))\simeq Hom_\Zp(\limind H_1(\G_n,H^2_{Iw}(F_\infty,\Zp(1))),\Zp)$.

$(iii)$ $Hom_\Zp(\mathcal H^2(F_\infty,\Zp(1)),\Zp)=E^0(H^2_{Iw}(F_\infty,\Zp(1)))$. En particulier, $H^2_{Iw}(F_\infty,\Zp(1))$ et $\limproj Br'_n$
possède le même $\L$-rang. Explicitement, il y a un isomorphisme naturel $Hom_\Zp(\mathcal H^2(F_\infty,\Zp(1)),\Zp)\simeq \oplus_{v\in S^{dec}}\L$.

$(iv)$  Il y a une suite exacte naturelle $$Ext^1_\Zp(\mathcal H^2(F_\infty,\Zp(1)),\Zp)\hookrightarrow
E^1(H^2_{Iw}(F_\infty,\Zp(1)))\twoheadrightarrow H^1(\Delta_k(H^2_S))$$

$(v)$ Pour $q\ge 2$, $H^q(\Delta_k(H^2_S))=E^q(H^2_{Iw}(F_\infty,\Zp(1)))$.

\noindent 4.  $(i)$ $\Delta_j(H^2_S)=\Delta_j(H^1_S)[2]$ est acyclique hors de $[-1,d-1]$.

$(ii)$ Les modules de cohomologie de $\Delta_j(H^2_S)$ sont noethériens de torsion. De plus:

- $H^q(\Delta_j(H^2_S))$ est pseudo-nul si $q\ne -1$.

- $H^{-1}(\Delta_j(H^2_S))\approx E^1(H^1(\Delta_k(H^2_S)))$.

 $(iii)$ Il y a des suites exactes naturelles $$H^{-1}(\Delta_j(H^2_S))\hookrightarrow H^2_{Iw}(F_\infty,\Zp(1))\rightarrow
E^0(X_\infty^*(H^2_S))\twoheadrightarrow H^{0}(\Delta_j(H^2_S))$$
$$E^1(Ext^1_\Zp(\mathcal H^2(F_\infty,\Zp(1)),\Zp))\hookrightarrow E^0(X_\infty^*(H^2_S))\twoheadrightarrow E^0(Hom_\Zp(\mathcal
H^2(F_\infty,\Zp(1)),\Zp))$$
$$\limproj Ker\, j_n^2\hookrightarrow H^{-1}(\Delta_j(H^2_S))\twoheadrightarrow \limproj_{n,k}  Coker j_n^2[p^k]$$
$$\limproj_{n,k} (Coker j_n^2)/p^k\hookrightarrow H^{0}(\Delta_j(H^2_S))\twoheadrightarrow \limproj_{n,k} H^1(\G_n,\mathcal  H^2(F_\infty,\Zp(1)))[p^k]$$

$(iv)$ Pour $q\ge 1$, il y a un isomorphisme $$E^{q+1}(Ext^1_\Zp(\mathcal H^2(F_\infty,\Zp(1)),\Zp))\simeq H^q(\Delta_j(H^2_S))$$
\end{prop}

\noindent Preuve: 1. est une conséquence directe de \ref{cledesc} $(i)$ et $(ii)$.

\noindent Le point $(i)$ de 2., 3., 4. résulte de \ref{cledesc} $(iii)$, \ref{tdtriv} $(i)$ et \ref{amplitude}.

2. $(ii)$ - $(iv)$ découlent immédiatement de 2. $(i)$. Ces résultats sont en fait bien connus, et peuvent aussi se lire sur des suites spectrales
convenables (e.g. \cite{Kat} 2.3 pour la codescente).

3. $(ii)$ D'après 2. $(ii)$, la flèche de descente $\Zp\underline\otimes X_\infty(H^2_S)\rightarrow H^2_S$ est un isomorphisme, et l'on peut donc
appliquer \ref{deltak}.

4. $(ii)$, on raisonne comme en \ref{exhaust} 2. $(iii)$.

Le reste de l'énoncé exprime essentiellement la cohomologie  des triangles distingués $\alpha_k(H^2_S)$ et $\alpha_j(H^2_S)$, aux exceptions
suivantes près:

3. $(iii)$ Comme $\C l'_n$ est fini, la suite exacte de Kummer (cf \ref{tdtriv}) montre que $Hom_\Zp(Br',\Zp)=Hom_\Zp(H^2_S,\Zp)$. La description
explicite de l'énoncé résulte donc de \ref{cleBr} $(i)$.

4. $(iii)$ Pour obtenir les deux dernières suites exactes, on raisonne comme en \ref{exhaust} 2. $(i)$.
\begin{flushright}$\square$\end{flushright}

Dans la proposition \ref{propZpi}, on s'est attaché à expliciter la cohomologie des triangles distingués $\alpha_j(A)$ et $\alpha_k(A)$ avec
$A=H^2_S$. Pour avoir une vision complète de la situation, il conviendrait de faire de même avec $A=H^1_S$, puis de recoller les deux avec $R\G_S$. A
défaut d'écrire les détails d'une étude aussi fastidieuse, on se contente d'évoquer deux points significatifs.
\begin{cor}\label{H1} $(i)$ $H^1_{Iw}(F_\infty,\Zp(1))\simeq E^0(RHom_\Zp(\mathcal H^1(F_\infty,\Zp(1)),\Zp))$.

$(ii)$ Il y a une suite exacte naturelle $$H^{1}(\Delta_k(H^2_S))\hookrightarrow Hom_\Zp(\mathcal H^1(F_\infty,\Zp(1)),\Zp)\rightarrow
E^0(H^1_{Iw}(F_\infty,\Zp(1)))\hspace{2cm}$$
$$\hspace{8.3cm}\rightarrow H^2(\Delta_k(H^2_S))\rightarrow \Zp(1)$$
En particulier, $H^{1}(\Delta_k(H^2_S))\simeq t_\L Hom_\Zp(\limind E'_n,\Zp)$.
\end{cor}
\noindent Preuve:  Pour $(i)$ (resp. $(ii)$), on écrit un morceau de la suite exacte longue de cohomologie de $\alpha_j(H^1_S)$ (resp.
$\alpha_k(H^1_S)$), compte-tenu de \ref{propZpi} 4. $(i)$. La dernière assertion résulte de la suite exacte, puisque $t_\L
E^0(H^1_{Iw}(F_\infty,\Zp(1)))=0$, et que $H^1_S\simeq E'$ (\ref{tdtriv} $(ii)$).
\begin{flushright}$\square$\end{flushright}

Le résultat suivant est l'ingrédient principal pour la preuve de \ref{propcl}.

\begin{prop} \label{propcle} Il existe des pseudo-isomorphismes non canoniques
$$Z\approx H^1(\Delta_k(H^2_S))\approx H^{-1}(\Delta_j(H^2_S))$$ En particulier $(Dec_2)\Leftrightarrow H^1(\Delta_k(H^2_S))\approx
0\Leftrightarrow H^{-1}(\Delta_j(H^2_S))\approx 0$.\end{prop} \noindent Preuve: Pour ne pas exclure le cas où $(Dec_1)$ n'est pas vérifiée,
introduisons $H^2_{S,dec}$, le noyau de la flèche composée $H^2_S\rightarrow Br'\rightarrow \oplus_{v\in S^{dec}}\underline\L$. On a alors deux
triangles distingués naturels
$$H^2_{S,non-dec}\rightarrow H^2_S\rightarrow \oplus_{v\in S^{dec}}\underline\L\rightarrow H^2_{S,non-dec}[1]$$
$$\C l'\rightarrow H^2_{S,non-dec}\rightarrow Br'_{non-dec}\rightarrow \C l'[1]$$
sur lesquels on lit les faits suivants:

- $C(k_{H^2_{S,non-dec}})\simeq C(k_{H^2_S})$. En particulier, $\Zp\underline\otimes X_\infty(H^2_{S,non-dec})\simeq H^2_{S,non-dec}$.

- $\Delta_k(H^2_S)\simeq \Delta_k(H^2_{S,non-dec})$.

- $f_\Zp(X_\infty(H^2_{S,non-dec})_{\G_n})\simeq f_\Zp((H^2_{S,non-dec})_n)\simeq (Br'_{non-dec})_n$.

\noindent D'après \ref{propZpi} 3. $(iii)$, le $\L$-module $M:=X_\infty(H^2_{S,non-dec})$ est de $\L$-torsion. On peut donc lui appliquer
\ref{deltak} 2. $(i)$ et \ref{critDeltak}, et cela donne une pseudo-surjection non-canonique $Z\simeq \limproj M_{\G_n}\rightarrow
H^1(\Delta_k(H^2_{S,non-dec}))\simeq H^1(\Delta_k(H^2_S))$.

Mais $H^1(\Delta_k(H^2_S))\simeq t_\L Hom_\Zp(\limind E'_n,\Zp)$ (\ref{H1} $(ii)$), or d'après \cite{Ne} 9.4.11 $(iii)$, $$\C har(Z)\mid \C har (t_\L
Hom_\Zp(\limind E'_n,\Zp))$$ (noter que $Hom_\Zp(\limind E_n',\Zp)\simeq Hom_\Zp(\limind E_n'\otimes\Qp/\Zp,\Qp/\Zp)$). Le noyau de la
pseudo-surjection précédente est donc nécessairement pseudo-nul. D'où le premier pseudo-isomorphisme de l'énoncé. Le reste suit, grâce à
\ref{propZpi} 4. $(ii)$ et \ref{cleBr} $(iii)$.
\begin{flushright}$\square$\end{flushright}

\begin{cor} \label{corZpi} $ $

$(i)$ $Ext^1_\Zp(\mathcal H^2(F_\infty,\Zp(1)),\Zp)$ ne possède aucun sous-module pseudo-nul.

$(ii)$ Sous $(Dec_2)$, $H^{-1}(\Delta_j(H^2_S))$ est le sous-module pseudo-nul maximal de $H^2_{Iw}(F_\infty,\Zp(1))$. \end{cor} \noindent Preuve:
$(i)$ Au vu de \ref{propZpi} 3. $(iv)$, il suffit de montrer que $E^1(X_\infty(H^2_S))$ ne possède aucun sous-module pseudo-nul. Mais d'après la
preuve de \ref{propcle}, $X_\infty(H^2_S)$ est somme directe d'un $\L$-module libre et du $\L$-module de torsion $X_\infty(H^2_{S,non-dec})$. D'où
$E^1(X_\infty(H^2_S))\simeq E^1(X_\infty(H^2_{S,non-dec}))$, et le résultat par \ref{lemmod} $(iii)$.

$(ii)$ Le $\L$-module $Ext^1_\Zp(\mathcal H^2(F_\infty,\Zp(1)),\Zp)$ étant de torsion (\ref{propZpi} 3. $(iv)$ et \ref{lemmod} $(i)$), les deux
premières suites exactes de \ref{propZpi} 4. $(iii)$, jointes à \ref{lemmod} $(iii)$, montrent que le sous-module pseudo-nul maximal de
$H^2_{Iw}(F_\infty,\Zp(1))$ est contenu dans $H^{-1}(\Delta_j(H^2_S))$. D'où le résultat, par \ref{propcle}.
\begin{flushright}$\square$\end{flushright}

Preuve de \ref{propcl}: $(i)$ Comme $\C l'=Ker(H^2_S\rightarrow Br')$, \ref{propZpi} 1. et  3. $(iii)$ montrent que $\limproj \C l'_n$ est noethérien
de torsion. Que le module $$Hom_\Zp(\limind \C l'_n,\Qp/\Zp)=H^1(X_\infty^*(\C l'))\simeq Ext_\Zp^1(\mathcal H^2_S(F_\infty,\Zp(1)),\Zp)$$ soit
noethérien de torsion résulte de \ref{propZpi} 1.,  3. $(iv)$ et \ref{lemmod} $(i)$.

Comme $X_\infty(\C l')$ est de $\L$-torsion, la suite exacte longue de cohomologie du triangle distingué $\alpha_k(\C l')$ s'écrit $$0\rightarrow
H^{0}(\Delta_k(\C l'))\rightarrow H^1(X_\infty^*(\C l'))\rightarrow E^1(X_\infty(\C l'))\rightarrow H^1(\Delta_k(\C l'))\rightarrow 0$$ La flèche
centrale est celle notée $\alpha_k$ dans l'énoncé, et il reste à montrer que les deux termes extrèmes sont pseudo-isomorphes.

Appliquer $\Delta_k$ à \ref{tdtriv} $(ii)$ donne un triangle distingué dont la suite exacte longue de cohomologie s'écrit:
$$H^0(\Delta_k(H^2_S))\rightarrow H^0(\Delta_k(\mathcal Cl'))\rightarrow
H^1(\Delta_k(Br'))\hspace{3cm} (***)$$ $$\hspace{4cm}\rightarrow H^1(\Delta_k(H^2_S))\rightarrow H^1(\Delta_k(\mathcal Cl'))\rightarrow
H^2(\Delta_k(Br'))$$ \noindent Dans celle-ci, le premier terme est trivial (\ref{propZpi} 3. $(i)$), et le dernier est pseudo-nul (\ref{cleBr} $(ii)$
et \ref{lemmod} $(i)$). Maintenant, \ref{propcle}, \ref{lemmod} $(iv)$ et \ref{cleBr} $(ii)$ donnent, non-canoniquement:
$$H^1(\Delta_k(H^2_S))\approx t_\L X_\infty(Br')\approx E^1(X_\infty(Br'))\simeq H^1(\Delta_k(Br'))$$ Comme $t_\L X_\infty(Br')$ est
pseudo-semi-simple, la suite exacte précédente montre que $H^0(\Delta_k(\mathcal Cl'))$ et $H^1(\Delta_k(\mathcal Cl'))$ le sont aussi, et qu'ils
possèdent le même idéal caractéristique. D'où $H^0(\Delta_k(\mathcal Cl'))\approx H^1(\Delta_k(\mathcal Cl'))$, non canoniquement.

$(ii)$ Puisque $X_\infty(\C l')$ est de torsion (cf $(i)$), c'est que $H^0(\Delta_k(\C l'))=Ker\, \alpha_k$. Mais d'après \ref{cleBr} $(ii)$
$$H^1(\Delta_k(Br'))\simeq E^1(X_\infty(Br'))\simeq E^1(Z)$$ et la suite exacte $(***)$ ci-dessus donne donc l'injection de l'énoncé:
$Ker\, \alpha_k\hookrightarrow E^1(Z)$. Si l'hypothèse $(Dec_2)$ est vérifiée, \ref{cleBr} $(iii)$ montre qu'en fait, $E^1(Z)=0$, et cela termine la
preuve.
\begin{flushright}$\square$\end{flushright}

Preuve de \ref{corNV1}: $(i)$ résulte directement de \ref{propcl}. \ref{corZpi} $(i)$.

$(ii)$ En toute généralité, la suite exacte longue de cohomologie du triangle distingué $\alpha_j(\C l')$ s'écrit $$0\rightarrow H^{-1}(\Delta_j(\C
l'))\rightarrow X_\infty(\C l')\rightarrow E^1(H^1(X_\infty^*(\C l')))\rightarrow H^0(\Delta_j(\C l'))\rightarrow 0$$ Dans celle-ci:

- Le troisième terme ne possède aucun sous-module pseudo-nul non nul (\ref{lemmod} $(iii)$).

- La flèche centrale coïncide avec $E^1(\alpha_k)$; sous $(Dec_2)$, c'est un pseudo-isomorphisme et le premier terme est donc pseudo-nul.

\noindent Aussi, $(\limproj \C l'_n)^0\simeq H^{-1}(\Delta_j(\C l'))$.

La suite exacte de l'énoncé est semblable à \ref{exhaust} 2. $(i)$.
\begin{flushright}$\square$\end{flushright}

Preuve de \ref{corNV2}: $(i)$ a été établi au cours de la preuve de \ref{propcle}.

$(ii)$ Sous $(Dec_3)$, on a $H^1(\Delta_k(Br'))=H^2(\Delta_k(Br'))=0$ (\ref{cleBr} $(iii)$). Le raisonnement de \ref{propcl} $(i)$ donne alors
$$Coker\, \alpha_k\simeq H^1(\Delta_k(\C l'))\simeq H^1(\Delta_k(H^2_S))\simeq t_\L Hom_\Zp(\limind E_n'\otimes\Qp/\Zp,\Qp/\Zp)$$
\begin{flushright}$\square$\end{flushright}

\begin{rem} \label{remrem}
1. $(i)$ L'équivalence de \ref{corNV1} $(i)$, longtemps suspectée, avait déjà été établie dans \cite{NV}, conditionnellement à la conjecture de Gross
pour chaque $F_n$; la preuve reposait alors de façon cruciale sur les résultats de \cite{LFMN}. Des résultats partiels figurent aussi dans \cite{LN}
et \cite{MC}, voir \cite{NV} à ce sujet.

$(ii)$ La divisibilité $\C har(Hom_\Zp(\limind \mathcal Cl'_n,\Qp/\Zp))\mid \mathcal Char(\limproj\C l'_n)$ était déjà connue (\cite{Ne} (prop.
9.4.1. $(iii)$). La preuve de \ref{propcle} présentée ici repose d'ailleurs sur les résultats de \cite{Ne} 9.4.

$(iii)$ Sous $(Dec_2)$, \ref{propcle} est en fait indépendant de \cite{Ne} 9.4. Dans le cas général, une réponse positive à la question
\ref{questionIF} engloberait le résultat de \emph{loc. cit.}. On peut se demander si une adaptation des méthodes de \emph{loc. cit} permettrait de
répondre à \ref{questionIF}; nous espérons y revenir dans un travail ultérieur. \vspace{0.1cm}\\
\noindent 2. Pour $T=\Zp(0)$, \ref{corZpi} $(i)$ possède un analogue bien connu, relié à la conjecture de Leopoldt faible. Celui-ci est
originellement dû à \cite{G1}, et \cite{N1} en a donné une autre preuve. On pourra aussi consulter \cite{Ne} 9.3.1 pour un résultat général.
\vspace{0.1cm}

\noindent 3. La proposition \ref{H1} $(i)$ montre que $f_\L H^1_{Iw}(F_\infty,\Zp(1))$ est réflexif, et $(ii)$ pourrait être le point de départ d'une
discussion sur la réflexivité du module $f_\L Hom_\Zp(\mathcal H^1(F_\infty,\Zp(1)),\Zp)$. Il serait sans doute intéressant de comparer cette
approche avec celle de \cite{G3}.
\end{rem}

Toute la présente étude concerne les modules d'Iwasawa attachés à la $\Zp$-représentation continue $T=\Zp(1)$. Pour $T=\Zp(i)$, l'étude analogue
présente quelques différences dont la plus marquante est peut-être la description du sous-module pseudo-nul maximal de $H^2_{Iw}(F_\infty,\Zp(i))$:
pour $i\ge 2$, celui-ci vaut tout simplement $\limproj Ker\, j_n^2$. Cette remarque est à rapprocher du fait suivant: les méthodes de \cite{V2},
\cite{V3} permettent de produire des éléments de $Ker\, j_n\subset H^2(G_{S,F_n},\Zp(i))$ (pour $n$ fixé) seulement pour $i\ge 2$, et non pour $i=1$.

\begin{q} Est-il possible d'adapter les méthodes de \cite{V2}, \cite{V3} pour produire des éléments du sous-quotient $\limproj_{n,k} Coker\, j_n^2[p^k]$
de $H^2_{Iw}(F_\infty,\Zp(1))$? \end{q}

\noindent Une réponse positive à cette question demanderait certainement la construction préalable de systèmes projectifs d'unités le long de
sous-extensions bien choisies de $F_\infty$.

\subsection{Descente et codescente dans la tour cyclotomique}

On propose ici une étude systématique de la (co-)descente pour les systè\-mes normiques usuels de la théorie d'Iwasawa cyclotomique. On retrouve
ainsi de façon concise et unifiée de nombreux résultats connus (\cite{Kuz}, \cite{Gr}, \cite{LFMN}, \cite{Be}, \cite{KN}, \cite{NL}, \cite{G4} et
\cite{I2}). Pour simplifier, on suppose $p\ne 2$.

L'approche axiomatique présente l'avantage de mettre en lumière les hypothèses dont dépendent les résultats. On ne fait ici aucun usage de la
conjecture principale. Par ailleurs, l'étude est essentiellement indépendante de la conjecture de Leopoldt (th. de Baker-Brumer, dans notre
contexte), à l'exception du corollaire \ref{corqNQD}. Les propriétés ``arithmétiques'' qu'on utilise sys\-té\-matiquement sont rassemblées dans le
lemme \ref{hyppart}.

Le phénomène de stabilisation joue un rôle important. Il est automatique pour les systèmes normiques qui proviennent de la cohomologie de $\Zp(1)$;
Pour les unités cyclotomiques, l'étude repose entièrement sur un résultat préalable de \cite{Be} (cf \ref{stabC}).\vspace{0.2cm}

Fixons $F$ un corps de nombres abélien, et considérons sa $\Zp$-extension cyclotomique $F_\infty/F$. Comme toujours, $\G_n=Gal(F_\infty/F_n)$,
$\G/\G_n=G_n$. Munies de la norme et de l'extension, les collections suivantes constituent les systèmes normiques que l'on souhaite étudier:

- $E_n:=\mathcal O_{F_n}^\times\otimes \Zp$ le $\Zp$-module des unités, $E=(E_n)\in {{_{\underline\G}\C}}_{tf}$.

- $\mathcal Cl_n:=\mathcal Cl(\mathcal O_{F_n})\otimes\Zp$ le $p$-groupe de classes, $\mathcal Cl=(\mathcal Cl_n)\in {{_{\underline\G}\C}}_{tf}$.

- $C_n\subset E_n$ le $\Zp$-module engendré par les unités cyclotomiques (cf \cite{Si} pour une définition précise), $C=(C_n) \in
{{_{\underline\G}\C}}_{tf}$.

- $B_n:=E_n/C_n$, $B=(B_n)\in {{_{\underline\G}\C}}_{tf}$. \vspace{0.1cm}

\noindent  On utilisera aussi les systèmes normiques auxiliaires suivants:

- $E'=(E_n')$, $\C l'=(\C l_n')$, $Br'=(Br_n')$, $R\G_S=R\underline\G(G_S,\Zp(1))$, $H^1_S=(H^1(G_{S,F_n},\Zp(1))\simeq E'$,
$H^2_S=(H^2(G_{S,F_n},\Zp(1)))$, $\underline\L_v=(\Zp[G_n/G_{n,v}])$ (voir section \ref{ClZpd} pour les notations).

- $I=(I_n)$ où $I_n$ est le $\Zp$-module libre sur les $p$-places de $F_n$, et les application de transitions sont induites par l'extension et la
norme des idéaux.\vspace{0.2cm}

Si $A$ est un système normique, on note  $$A_\infty=\limind A_n,\, X_\infty=\limproj A_n,\, j_n:A_n\rightarrow (A_\infty)^{\G_n},\, \hbox{et}\,
k_n:(X_\infty)_{\G_n}\rightarrow A_n$$ Si plusieurs systèmes normiques sont en jeu, on écrira parfois $A_\infty(A)$, $X_\infty(A)$, $j_{n,A}$ ou
encore $k_{n,A}$ pour indiquer auquel on se rapporte.



\begin{prop} \label{stabcyclo} Si $A=E,\mathcal Cl,C$ ou $B$, alors:

\noindent 1. $X_\infty$ est un $\L$-module de type fini.

\noindent 2. $(i)$ Les systèmes inductifs $H_1(\G_n,X_\infty)$, $Ker\, k_n$ et $Coker\, k_n$ se stabilisent; en particulier, leur limite est de type
fini sur $\Zp$.

$(ii)$ Les systèmes projectifs $Ker\, j_n$ , $Coker\, j_n$, $H^1(\G_n,A_\infty)$ et $H^2(\G_n,A_\infty)$ se stabilisent; en particulier leur limite
est de torsion sur $\Zp$.

\noindent 3. Les $\L$-modules $Hom_\Zp(A_\infty,\Zp)$ et $Ext^1_\Zp(A_\infty,\Zp)=Hom_\Zp(t_\Zp A_\infty,\Qp/\Zp)$ sont de type fini.
\end{prop}

\noindent Preuve: Commençons par un lemme:

\begin{lem} \label{lemloc} A lieu dans ${_{\underline\G}\C}$ une suite exacte ``de localisation hors des $p$-places'' qui se matérialise par une famille
de suites exactes: {\small$$E_n\hookrightarrow H^1(G_{S,F_n},\Zp(1))\rightarrow I_n\rightarrow \mathcal Cl_n\rightarrow
H^2(G_{S,F_n},\Zp(1))\rightarrow \oplus_{v\mid p}\Zp[G_n/G_{n,v}]\twoheadrightarrow \Zp$$}
\end{lem}

\noindent Preuve: La construction directe (par localisation) d'un triangle distingué adéquat ne semblant pas relever des méthodes de cet article, on
se contentera d'une construction à la main. D'abord, la théorie de Kummer donne comme en \ref{tdtriv} $(ii)$ $E_n'\simeq H^1(G_{S,F_n},\Zp(1))$ et
$\mathcal Cl_n'\hookrightarrow H^2(G_{S,F_n},\Zp(1))\rightarrow Br'_n$. La suite exacte de l'énoncé est alors obtenue en recollant cette dernière à
la suite exacte de valuation: $E_n\hookrightarrow E_n'\rightarrow \Zp[S(F_n)]\rightarrow \mathcal Cl_n$, via la flèche naturelle $\mathcal
Cl_n\rightarrow \mathcal Cl_n'$. La famille de suites exactes ainsi construite est $G_n$-équivariante, compatible à la norme et à l'extension, et
donne donc bien une suite exacte dans ${_{\underline\G}\C}$.
\begin{flushright}$\square$\end{flushright}

\noindent Preuve de \ref{stabcyclo}:  1. Dans la suite exacte obtenue en appliquant $\limproj_n$ à celle du lemme précédent, on connaît déjà la
noethérianité du second et du cinquième terme grâce à \ref{propZpi} 1., et celle des troisième et sixième terme est évidente, compte-tenu de leur
description explicite. Le résultat suit.

  2. $(i)$ D'après (\ref{cledesc} $(ii)$), $C(k_{R\G_S})=0$. Par \ref{Hstab} 2., on en déduit que $C(k_{H^q_S})$ se stabilise pour $q=1,2$. Par
ailleurs, les systèmes normiques $I$ et $\underline\L_v$ se stabilisent via la corestriction, puisque pour $n$ suffisamment grand, l'extension
$F_\infty/F_n$ est totalement ramifiée en chaque $p$-place. Il s'ensuit que $C(k_{I})$ et $C(k_{\underline\L_v})$ se stabilisent automatiquement
(\ref{stabstab}), puis $C(k_E)$ et $C(k_{\C l})$ aussi, par \ref{stabtech2}, appliqué à la suite exacte \ref{lemloc}.

Reste la stabilisation de $C(k_C)$ et $C(k_B)$. Par définition de $B$, il y a une suite exacte $C\hookrightarrow E\twoheadrightarrow B$, et il nous
suffit donc par \ref{stabtech}  d'établir la stabilisation de $Coker\, k_n$ pour $A=C$ (condition $(i)$, \emph{loc. cit.}). Or celle-ci est assurée
par la

\begin{prop} \label{stabC} (\cite{Be}, lemme 2.5) $Coker\, k_n$ se stabilise dès que toutes les $p$-places sont totalement ramifiées
dans $F_\infty/F_n$. En particulier, $Coker\, k_n$ se stabilise. \end{prop} La preuve de ce résultat repose sur un examen attentif de la définition
des unités cyclotomiques.
\begin{flushright}$\square$\end{flushright}

Terminons la preuve de \ref{stabcyclo}: 2. $(ii)$ procède de \ref{stabjk} $(i)$, et 3. découle de 1. et 2. $(i)$, compte-tenu de \ref{critnoethcod}.
\begin{flushright}$\square$\end{flushright}

Pour aller plus loin, nous devons prendre en compte la spécificité des quatre systèmes normiques qui nous occupent. Rassemblons en un lemme les
propriétés qui nous seront utiles:
\begin{lem} \label{hyppart} $ $

$(i)$ Pour $A=E$: $t_\Zp A_n=\mu_{p^\infty}(F_n)$, $Ker\, j_n=0$ et $Coker\, j_n=0$.

$(ii)$ Pour $A=C$: $t_\Zp A_n=\mu_{p^\infty}(F_n)$, $Ker\, j_n=0$ et $Coker\, j_n$ est fini.

$(iii)$ Pour $A=\mathcal Cl$: $A_n$ est fini et $Coker\, k_n=0$ pour $n>>0$.

$(iv)$ Pour $A=B$: $B_n$ est fini.

\end{lem}
Preuve: $(i)$ est clair.

$(ii)$ La trivialité de $Ker\, j_n$ se déduit de la propriété analogue pour $E$. Quant à $Coker j_n$, c'est d'une part un sous-quotient de $E_n$ donc
un $\Zp$-module de type fini, et d'autre part un $\Zp$-module de torsion, comme chaque objet de cohomologie de $C(j_A)$. Il est donc nécessairement
fini.

$(iii)$  Comme $k_n$ est un isomorphisme pour $A=(H^2(G_{S,F_n},\Zp(1)))$ et $A=(Br'_n)$, appliquer le foncteur homologique $H^\bullet(C(k_{(-)}))$
aux deux suites exactes suivantes (découpées dans \ref{lemloc}): $I\rightarrow \C l\rightarrow\C l'\rightarrow 0$, et $0\rightarrow \C l'\rightarrow
H^2_S\rightarrow Br'_n\rightarrow 0$ en donne une troisième: $Coker\, k_{n,I}\rightarrow Coker\, k_{n,\C l}\rightarrow Ker\, k_{n,Br'}$, dans
laquelle les deux termes extrèmes s'annulent dès que $F_\infty/F_n$ est totalement ramifiée en chaque $p$-place. On a donc bien $Coker\, k_n=0$ pour
$A=\C l$ et $n>>0$.

$(iv)$ La finitude de $B_n$ est bien connue (\cite{Si} th. 4. 1 donne même une formule pour l'ordre de $B_n$).
\begin{flushright}$\square$\end{flushright}


\begin{rem} \label{remleop}

$(i)$ (Leopoldt) Pour $A=B$, $Ker\, k_n$ est fini $\forall n$. En effet, la finitude de $Ker\, k_n$ (ie. de $(\limproj B_m)_{\G_n}$) résulte de la
conjecture de Leopoldt sous sa forme ``fonctions $L$ $p$-adiques'', modulo la généralisation par \cite{Ts} aux corps abéliens quelconques du théorème
d'Iwasawa reliant unités cyclotomiques et fonctions $L$ $p$-adiques. On renvoie à \cite{Be2} th. 1.1 pour plus d'explications.

$(ii)$ Considérons la condition suivante: $$(S_p=S_p^+): \hbox{$F_\infty$ et $F_\infty^+$ possèdent le même nombre de $p$-places}$$ \noindent  Pour
$A=\C l$, cette condition assure la finitude de $Ker\, k_n$ $\forall n$. En effet, pour la partie $-$, cela résulte de la conjecture de Gross (cf
\cite{Ko}), et de celle de Leopoldt pour la partie $+$ (cf \cite{G0}, th. 1).

En fait, $(S_p=S_p^+)\Leftrightarrow Ker\, k_n$ est fini $\forall n$, comme nous le verrons plus loin (\ref{remfinale}) 2.
\end{rem}

Reste maintenant à écrire les résultats de la proposition \ref{exhaust} pour les quatre systèmes normiques qui nous occupent. On procède en trois
étapes:

- \ref{resAB} énoncé des résultats pour $A=\C l$ ou $B$.

- \ref{resEC} énoncé des résultats pour $A=E$ ou $C$.

- \ref{compfinale} lien entre les deux premières étapes.

\noindent à chaque étape, on indique une liste de corollaires (connus ou non), dont les notations suivantes simplifieront les énoncés.
\begin{defn}

\noindent 1. Si $M$ est un $\L$-module de type fini, on note:

$(i)$ $M^{\delta}=\cup M^{\G_n}\simeq \limind H_1(\G_n,M)$.

$(ii)$ $M^0=t_\Zp M^\delta$ son sous-module fini maximal.

\noindent 2.  $(ii)$ Dualement, si $M$ est un $\L$-module topologique discret ($\Rightarrow$ de $\Zp$-torsion) tel que $Hom_\Zp(M,\Qp/\Zp)$ soit un
$\L$-module de type fini, on note:

$(i)$ $M_\delta=Hom_\Zp(Hom_\Zp(M,\Qp/\Zp)^\delta,\Qp/\Zp)\simeq \limproj M_{\G_n}\simeq \limproj H^1(\G_n,M)$.

$(ii)$ $M_0=Hom_\Zp(Hom_\Zp(M,\Qp/\Zp)^0,\Qp/\Zp)\simeq \limproj_{n,k}M_{\G_n}/p^k$. \end{defn}

\begin{prop} \label{resAB}
Soit $A=\C l$ ou $B$, alors:

\noindent 1. Les $\L$-modules $X_\infty$ et $Hom_\Zp(A_\infty,\Qp/\Zp)$ sont noethériens de torsion.

\noindent 2. Il y a deux flèches adjointes naturelles:
$$\alpha_j: X_\infty\rightarrow E^1(Hom_\Zp(A_\infty,\Qp/\Zp))\hspace{0.5cm}
\alpha_k:Hom_\Zp(A_\infty,\Qp/\Zp)\rightarrow E^1(X_\infty)$$ et celles-ci vérifient:

$(i)$ Non canoniquement, $Ker\, \alpha_j$, $Coker\, \alpha_j$, $Ker\, \alpha_k$, $Coker\, \alpha_k$ sont tous pseudo-isomorphes.

$(ii)$ Noyau et conoyau de $\alpha_j$ sont décrits par:
$$\limproj Ker\, j_n \hookrightarrow  Ker\, \alpha_j \twoheadrightarrow \limproj_{n,k} Coker\, j_n[p^k]$$
$$Ker\alpha_j\simeq \limind H_1(\G_n,X_\infty)$$
$$\limproj_{n,k} Coker\, j_n/p^k\hookrightarrow Coker\,\alpha_j \twoheadrightarrow \limproj_{n,k} H^1(\G_n,A_\infty)[p^k]$$
$$Coker\, \alpha_j\simeq \limind Ker k_n$$ De plus, on a $$E^2(Hom_\Zp(A_\infty,\Qp/\Zp))\simeq \limproj_{n,k} H^1(\G_n,A_\infty)/p^k\simeq \limind Coker\, k_n$$ Pour $A=\C l$, ce dernier
module vaut $0$. \vspace{0.1cm}

$(iii)$ Noyau et conoyau de $\alpha_k$ sont décrits par:
$$Ker\alpha_k\simeq Hom_\Zp(\limproj
H^1(\G_n,A_\infty),\Qp/\Zp)$$ $$Hom_\Zp(\limind Coker\, k_n,\Qp/\Zp)\hookrightarrow Ker\alpha_k\twoheadrightarrow  Hom_\Zp(\limind Ker\, k_n,\Zp)$$
\noindent Pour $A=\C l$, le premier terme de la suite exacte précédente est trivial.

Pour $A=\C l$ ou $B$, on a de plus
$$Coker\, \alpha_k\simeq Hom_\Zp(\limproj Coker\, j_n,\Qp/\Zp)$$
$$Hom_\Zp(t_\Zp \limind Ker\, k_n,\Qp/\Zp)\hookrightarrow Coker\, \alpha_k\twoheadrightarrow Hom_\Zp(\limind H_1(\G_n,X_\infty),\Zp)$$
$$E^2(X_\infty)\simeq Hom_\Zp(\limproj Ker\, j_n,\Qp/\Zp)\simeq Hom_\Zp(t_\Zp \limind H_1(\G_n,X_\infty),\Qp/\Zp)$$
\end{prop}

Preuve: Le corollaire \ref{corexhaust} donne  1. et 2. $(i)$. Le point 2. $(ii)$ (resp. 2. $(iii)$) est donné par \ref{exhaust} 2. $(i)$ et 3. $(ii)$
(resp. \ref{exhaust} 2. $(ii)$ et 3. $(iii)$).
\begin{flushright}$\square$\end{flushright}

\begin{rem} \label{remAB} Dans la proposition précédente,

$(i)$ Si $Ker\, \alpha_j$ est fini, alors le troisième terme de chacune des quatre suites exactes de \ref{resAB} disparaît, puisque le terme central
doit être fini (\ref{resAB} 2. $(i)$). De plus, les suites exactes en question mettent en évidence l'équivalence entre la finitude des modules
suivants: $Ker\, \alpha_j$, $\limproj Coker\, j_n$, $\limproj H^1(\G_n,A_\infty)$, $\limind Ker\, k_n$, $\limind H_1(\G_n,X_\infty)$. Comme les
systèmes projectifs et inductifs en jeu se stabilisent, on voit facilement que la finitude de la limite équivaut à celle de tous les arguments; par
exemple, $\limind Ker\, k_n$ est fini si et seulement si chacun des $Ker\, k_n$ l'est. On renvoie à \ref{remleop} pour l'interprétation arithmétique
de cette condition.

$(ii)$ Pour $A=B$, la conjecture de Leopoldt assure la finitude de $Ker\, \alpha_j$ (cf. $(i)$ et \ref{remleop} $(i)$). Dans cette situation l'expert
appréciera de voir la finitude des groupes de cohomologie $H^q(\G_n,\limind B_m)$ apparaître sans plus d'efforts. Dans ce contexte, on note que
\ref{resAB} 2. $(iii)$ donne deux isomorphismes: $$Hom_\Zp(\limproj Coker\, j_n,\Qp/\Zp)\simeq Hom_\Zp(\limind Ker\, k_n,\Qp/\Zp)\simeq
Coker\,\alpha_k$$ dont le premier répond au attentes de \cite{LN} Cor. 4. 6, en toute généralité.

$(iii)$ Pour $A=\C l$, les résultats sont indépendants de toute hypothèse sur $F$. Par ailleurs, il tiennent \textit{verbatim} si l'on remplace $\C
l$ par $\C l'$. Dans ce contexte, le rapport entre $Coker\, \alpha_k$, $\limproj Coker\, j_n$ et $\limind\, Ker\, k_n$ est à rapprocher de
l'interprétation de \cite{LFMN} de ces deux derniers objets en termes du module de Bertrandias-Payan (cf \textit{loc. cit}).


\end{rem}

\begin{cor} \label{remreg} $ $ Avec les notations de la définition précédente, on a:

\noindent 1. $(i)$  $(\limproj \C l_n)^0\simeq \limproj Ker\, j_{n,\C l}$.

$(ii)$ $(\limind \C l_n)_0\simeq\limind Coker\, k_{n,\C l}=0$.

\noindent 2. $(i)$  $(\limproj B_n)^0\simeq \limproj Ker\, j_{n,B}$.

$(ii)$ $(\limind B_n)_0\simeq \limind Coker\, k_{n,B}$.
\end{cor}
Preuve: D'après \ref{resAB}, 2. $(ii)$, on a $Ker\alpha_j=(X_\infty)^\delta$. Les points 1. $(i)$ et 2. $(i)$ découlent aussitôt de \ref{resAB} 2.
$(ii)$. De manière analogue, les points 1. $(ii)$ et 2. $(ii)$ s'obtiennent en dualisant l'isomorphisme
$Ker\alpha_k=Hom_\Zp((A_\infty)_\delta,\Qp/\Zp)$, à l'aide de \ref{resAB} 2. $(iii)$.
\begin{flushright}$\square$\end{flushright}

\begin{prop} \label{resEC} Soit $A=E$ ou $C$, $X_\infty=\limproj A_n$, $A_\infty=\limind A_n$. Alors:

\noindent 1. Les $\L$-modules $X_\infty$ et $Hom_\Zp(A_\infty,\Zp)$ sont noethériens.

\noindent 2. Il y a deux flèches duales naturelles: $$\tilde\alpha_j:X_\infty\rightarrow E^0(Hom_\Zp(A_\infty,\Zp))\hspace{0.5cm}
\alpha_k:Hom_\Zp(A_\infty,\Zp)\rightarrow E^0(X_\infty)$$ et celles-ci vérifient:

$(i)$ Noyau et conoyau de $\tilde \alpha_j$ sont décrits par:
$$Ker\, \tilde\alpha_j=\limproj \mu_{p^\infty}(F_n)\hspace{1.5cm} Coker\, \tilde\alpha_j\simeq \limind Ker\, k_n$$
$$\limproj Coker\, j_n\hookrightarrow Coker\, \tilde\alpha_j\twoheadrightarrow \limproj_{n,k}
H^1(\G_n,A_\infty)[p^k]$$ Dans la suite exacte précédente, le premier terme est fini si $A=C$, trivial si $A=E$.

De plus,  $E^1(Hom_\Zp(A_\infty,\Zp))\simeq \limind Coker\, k_n$ et il y a une suite exacte $$\limproj_{n,k} H^1(\G_n,A_\infty)/p^k\hookrightarrow
E^1(Hom_\Zp(A_\infty,\Zp))\twoheadrightarrow \limproj H^2(\G_n,A_\infty)[p^k]$$

$(ii)$ Noyau et conoyau de $\alpha_k$ sont décrits par:
$$Ker\, \alpha_k\simeq Hom_\Zp(\limproj H^2(\G_n,A_\infty),\Qp/\Zp)\simeq Hom_\Zp(\limind Coker\, k_n,\Zp)$$ $$Coker\,\alpha_k\simeq
Hom_\Zp(\limproj H^1(\G_n,A_\infty),\Qp/\Zp)$$ $$Hom_\Zp(t_\Zp \limind Coker\, k_n,\Qp/\Zp)\hookrightarrow Coker\, \alpha_k\twoheadrightarrow
Hom_\Zp(\limind Ker \, k_n,\Zp)$$ De plus, il y a une suite exacte $$\limproj\mu_{p^\infty}(F_n)\hookrightarrow E^1(X_\infty)\twoheadrightarrow
Hom_\Zp(t_\Zp \limind Ker\, k_n,\Qp/\Zp)$$ et $Hom_\Zp(t_\Zp \limind Ker\, k_n,\Qp/\Zp)\simeq Hom_\Zp(\limproj Coker\, j_n,\Qp/\Zp)$; ce dernier
module est trivial si $A=E$.
\end{prop}

Preuve: 1. vient de \ref{exhaust} 1. et 3. Si $\mu_p(F)=0$, alors 2. $(i)$ (resp. 2. $(ii)$) est donné par \ref{exhaust} 2. $(i)$ et 3. $(ii)$ (resp.
\ref{exhaust} 2. $(ii)$ et 3. $(iii)$). Pour traiter le cas $\mu_p(F)\ne 0$, on introduit le système normique des racines de l'unité:
$\underline\Z_p(1)=(\mu_{p^\infty}(F_n))$. Comme $\mu_{p^\infty}(F_\infty)$ est $\G$-cohomologiquement trivial, on voit que
$C(j_{\underline\Z_p(1)})$ et $C(k_{\underline\Z_p(1)})$ sont triviaux, si bien qu'il n'y a aucune difficulté à se ramener à l'étude du système
normique $A/\underline\Z_p(1)$; laquelle est analogue au cas $\mu_p(F)=0$. \begin{flushright}$\square$\end{flushright}

Si $M\hookrightarrow \L^r$, il est bien connu (et facile à démontrer) que $M$ est libre si et seulement si le conoyau est sans $\Zp$-torsion. De
plus, la $\Zp$-torsion du conoyau en question ne dépend pas du plongement choisi, à unique isomorphisme près; on définit donc ainsi un invariant
structurel de $M$, appelé suggestivement le défaut de liberté du $\L$-module $M$.

\begin{cor} $ $

\noindent 1. $(i)$ On a $t_\L \limproj E_n=\Zp(1)$ ou $0$ selon que $F$ contient ou non $\mu_p$;\vspace{-.2cm} {\small $$t_\L Hom_\Zp(\limind
E_n,\Zp)\simeq Hom(\limproj_n H^2(\G_n,\limind_m E_m),\Qp/\Zp)\simeq Hom_\Zp(\limind Coker\, k_{n,E},\Zp)$$}\vspace{-.2cm}

$(ii)$ $f_\L \limproj E_n$ est libre, et le défaut de liberté de $f_\L Hom_\Zp(\limind E_n,\Zp)$ est dual à $\limproj_{n,k} H^1(\G_n,\limind_m
E_m)/p^k\simeq t_\Zp \limind Coker\, k_{n,E}$.

\noindent 2. $(i)$ On a $t_\L \limproj C_n=\Zp(1)$ ou $0$ selon que $F$ contient ou non $\mu_p$;\vspace{-.2cm} {\small $$t_\L Hom_\Zp(\limind
C_n,\Zp)\simeq Hom(\limproj_n H^2(\G_n,\limind_m C_m),\Qp/\Zp)\simeq Hom_\Zp(\limind Coker\, k_{n,C},\Zp)$$}\vspace{-.2cm}

$(ii)$ Le défaut de liberté du module $f_\L \limproj C_n$ est isomorphe à $\limproj Coker\, j_{n,C}$ $\simeq t_\Zp\limind Ker\, k_{n,C}$, et celui de
$f_\L Hom_\Zp(\limind C_n,\Zp)$, dual à $\limproj_{n,k} H^1(\G_n,\limind_m C_m)/p^k$ $\simeq t_\Zp \limind Coker\, k_{n,C}$.
\end{cor}
Preuve: C'est une conséquence immédiate de la proposition qui précéde, compte-tenu du fait que le but de $\tilde\alpha_j$, comme celui de $\alpha_k$,
est un $\L$-module libre, puisqu'obtenu en appliquant $E^0$ à un $\L$-module de type fini.
\begin{flushright}$\square$\end{flushright}

\begin{rem}
1. $(i)$ apparaît en filigrane dans \cite{I2}.

1. $(ii)$ est généralement attribué à \cite{Kuz}. Voir aussi \cite{Gr}, pour une preuve algébrique de la première partie, dans le contexte de l'étude
axiomatique de systèmes normiques particuliers.

2. $(iii)$ La première partie de l'énoncé est dûe à \cite{Be}. Le module $\limproj\,  Coker j_{n,C}$ a fait l'objet de nombreuses études, cf
\cite{BN} pour quelques références.
\end{rem}

\begin{prop} \label{compfinale} $ $

$(i)$ Il y a un diagramme commutatif naturel, dans lequel la colonne centrale est exacte: \vspace{-0.5cm}{\small$$\diagram{&&0&&\cr
&&\vflcourte{}{}&&\cr &&\limind Ker\, k_{n,E}&\hflcourte{\sim}{}&\limproj_{n,k} H^1(\G_n,\limind_m E_m)[p^k] \cr &&\vflcourte{}{}&&\cr &&\limproj
I_n&& \cr &&\vflcourte{r}{}&&\cr \limproj Ker\, j_{n,\C l}& \injfl{}{}&Ker\, \alpha_{j,\C l}&\surjfl{}{}&\limproj_{n,k} Coker\, j_{n,\C l}[p^k] \cr
\vflcourte{}{}&&\vflcourte{}{}&&\vflcourte{}{}\cr\limproj_{n,k} H^1(\G_n,\limind_m E_m)/p^k&\injfl{}{}& \limind Coker\,
k_{n,E}&\surjfl{}{}&\limproj_{n,k} H^2(\G_n,\limind_m E_m)[p^k] \cr &&\vflcourte{}{}&&\cr
 && \limproj Br_n'&&
\cr &&\vflcourte{}{}&& \cr \limproj_{n,k} Coker\, j_{n,\C l}/p^k&\injfl{}{}&\limind Ker\, k_{n,\C l}&\surjfl{}{}&\limproj_{n,k} H^1(\G_n,\limind \C
l_m)[p^k]\cr &&\vflcourte{}{}&&\cr &&0&&}$$}

Dans ce diagramme, $r:\limproj I_n\rightarrow Ker\,\alpha_{j,\C l}=(\limproj \C l_m)^{\delta}$ est induite par la flèche naturelle $I_n\rightarrow \C
l_n$.

$(ii)$ Il y a une suite exacte naturelle $$\limproj Ker\, j_{n,\C l}\hookrightarrow \limproj_n H^1(\G_n,\limind_m E_m)\rightarrow (\limproj
I_n)\otimes \Qp/\Zp \rightarrow \limproj Coker\, j_{n,\C l} \hspace{3cm}$$
$$\hspace{2cm} \rightarrow \limproj_n H^2(\G_n,\limind_m E_m)\rightarrow (\limproj
Br'_n)\otimes\Qp/\Zp\twoheadrightarrow \limproj H^1(\G_n,\limind \C l_m)$$

\end{prop}

\noindent Preuve: Le moyen le plus rapide serait l'utilisation du triangle distingué sous-jacent à la suite exacte de localisation \ref{lemloc}.
Comme il n'est pas disponible directement, on en fabrique un à la main:

\begin{lem} \label{lemtdloc} On peut trouver $Z,Z'\in D^b({{_{\underline\G}\C}})$, et un triangle distingué de $D^b({_{\underline\G}\C})_{tf}$:
$T=\left(Z\rightarrow R\G_S\rightarrow Z'\rightarrow  Z[1]\right)$ vérifiant le propriétés suivantes: \vspace{0.1cm}

$(i)$ $Z$ et $Z'$ sont acycliques en degrés $\ne 1,2$, et la suite exacte longue de cohomologie de $T$ s'identifie à \ref{lemloc}, tronquée à
l'avant-dernier terme:  \vspace{-0.7cm}{\small
$$\diagram{H^1(Z)&\injfl{}{}& H^1(R\G_S)&\hflcourte{}{}& H^1(Z')&\hflcourte{}{}&
H^2(Z)&\hflcourte{}{}&H^2(R\G_S)&\surjfl{}{}& H^2(Z')\cr
\vflcourte{\wr}{}&&\vflcourte{\wr}{}&&\vflcourte{\wr}{}&&\vflcourte{\wr}{}&&\vflcourte{\wr}{}&&\vflcourte{\wr}{}\cr
E&\injfl{}{}&H^1_S&\hflcourte{}{}&I&\hflcourte{}{}& \C l&\hflcourte{}{}& H^2_S&\surjfl{}{}& Br'}$$}   \vspace{-0.8cm}

$(ii)$ $\Delta_j(Z')\simeq X_\infty(Z')[1]$.

$(iii)$ $\Delta_j(Z)\simeq \Delta_j(Z')[-1]$.
\end{lem}

\noindent Preuve: $(i)$ L'existence de $Z,Z'$ et $T$ vérifiant $(i)$ relève du lemme \ref{lemtd}. Les propriétés $(ii)$ et $(iii)$ se déduisent de
$(i)$. En effet:

$(ii)$ Puisque $\tau_{\le 1}Z'=I$ et $\tau_{\ge 2}Z'=Br'$ se stabilisent via la corestriction, on voit tout de suite que la cohomologie de
$A_\infty(Z')$ est uniquement $p$-divisible. D'où, automatiquement, la trivialité de $X_\infty^*(Z')\simeq RHom_\Zp(A_\infty(Z'),\Zp)$, si bien que
le triangle distingué $\alpha_j(Z')$ se réduit à l'isomorphisme de l'énoncé.

$(iii)$ est obtenu en appliquant $\Delta_j$ à $T$, puisque $\Delta_j(R\G_S)=0$ (\ref{cledesc} $(iii)$). \begin{flushright}$\square$\end{flushright}

Retour à la preuve de \ref{compfinale}: Tronquer $Z$ et appliquer le foncteur $\Delta_j$ donne un triangle distingué
$$\Delta_j(E[-1])\rightarrow \Delta_j(Z)\rightarrow \Delta_j(\C l[-2])\rightarrow \Delta_j(E)$$ dans lequel le premier (resp. troisième) sommet
est \emph{a priori} (\ref{amplitude}) acyclique hors de $[0,2]$ (resp. $[1,3]$), et le second, hors de $[1,2]$ (\ref{lemtdloc} $(ii)$). D'où une
suite exacte longue de cohomologie \begin{flushleft}$H^{0}(\Delta_j(E))\hookrightarrow H^1(\Delta_j(Z))\rightarrow H^{-1}(\Delta_j(\C l))\rightarrow
H^{1}(\Delta_j(E))$\end{flushleft} \begin{flushright}$\rightarrow H^2(\Delta_j(Z))\twoheadrightarrow H^{0}(\Delta_j(\C l))$\end{flushright} dont il
reste à vérifier qu'elle donne bien la colonne du diagramme de l'énoncé. La description explicite des second et cinquième termes relève de
\ref{lemtdloc}. Pour les autres, et pour les suites exactes horizontales, on applique \ref{resAB} 2. $(ii)$ et \ref{resEC} 2. $(i)$, compte-tenu de
la description suivante de la cohomologie de $\Delta_j(A)$, valable aussi bien pour $A=E$ que pour $A=\C l$ (\ref{corstab} 2. $(ii)$):

- $H^{-1}(\Delta_j(A))\simeq Ker\alpha_j$ et $H^0(\Delta_j(A))\simeq Coker\alpha_j\simeq \limind Ker\, k_n$.

- $H^1(\Delta_k)\simeq \limind Coker\, k_n$.

Vérifier que la flèche $r$ est bien induite par $I_n\rightarrow \C l_n$ ne pose pas de diffculté particulière, si l'on observe l'effet du morphisme
fonctoriel $\Delta_j[-1]\rightarrow X_\infty$ (début du triangle $\alpha_j$) sur la flèche composée suivante (cf. \ref{lemtdloc}): $$I\simeq
(\tau_{\le 1}Z')[1]\rightarrow Z'[1]\rightarrow Z[2]\rightarrow (\tau_{\ge 2} Z)[2]\simeq \C l$$

$(ii)$ La suite exacte de l'énoncé est obtenue en appliquant le foncteur $\limind C(j_{(-)})$ au triangle $T$ de \ref{lemtdloc}, compte-tenu de
\ref{corstab} 2. $(iii)$, et du fait que $C(k_A)=(H_1(\G_n,X_\infty))[2]\simeq (X_\infty^{\G_n})[2]$, pour $A=I$ ou $Br'$.
\begin{flushright}$\square$\end{flushright}

\begin{rem} \label{remfinale}
Une suite exacte à sept termes, similaire à celle de \ref{compfinale} 2. apparaît dans la littérature en plusieurs endroits (e.g. \cite{I2} et ref.).
Cette suite exacte présente plusieurs intérêts:

1. Comme on sait que les limites en question se stabilisent, on obtient une information sur la structure des groupes $H^q(\G_n,\limind E_m)$, $q=1,2$
pour $n>>0$. C'est l'objet de \cite{I2}.

2. On voit se produire le pendant algébrique du phénomène des zéros trivaux: comme $H^q(\G_n,\limind E_m)$ est fixé par la conjugaison complexe, on
obtient dans la partie imaginaire: $(\limproj I_n\otimes\Qp/\Zp)^-\simeq \limproj (Coker\, j_{n,\C l})^-$. La conjecture de Leopoldt prédisant par
ailleurs la finitude de $(Coker\, j_{n,\C l})^+$ (cf \ref{remleop}), on voit que $rg_\Zp (\limproj \C l_m)^\delta=corg_\Zp \limproj (Coker\, j_{n,\C
l})$ (cf \ref{resAB} 1. $(ii)$) coïncide avec le nombre de $p$-places qui se décomposent dans $F_\infty/F_\infty^+$.
\end{rem}

\begin{cor} $(i)$ $t_\Zp \limind  Coker\, k_{n,E}\simeq (\limproj \C l'_n)^0$. \end{cor}
Preuve: Cela découle de la description explicite de la flèche $r$ dans \ref{compfinale} $(i)$. \begin{flushright}$\square$\end{flushright}

\begin{cor} \label{corqNQD} (utilise la conjecture de Leopoldt, cf \ref{remleop} $(ii)$)

Il y a un diagramme commutatif à lignes exactes: {\small
$$\diagram{&&\limproj_{n,k} H^1(\G_n,C_\infty)[p^k]&\injfl{}{} &\limproj_{n,k} H^1(\G_n,E_\infty)[p^k]& \hflcourte{}{}&\limproj Cok\, j_{n,B}\cr
                        &&           \vflupcourte{}{}&           &\vflupcourte{\wr}{}&   &                                       \vflcourte{\wr}{} \cr
(\limproj B_m)^0&\injfl{}{}&  \limind Ker\,  k_{n,C}&\hflcourte{}{}& \limind Ker\, k_{n,E}&\hflcourte{}{}                        &\limind Ker\, k_{n,B}
}$$ \vspace{-1cm}$$
\diagram{&\hflcourte{}{}&(\limproj_n H^1(\G_n,C_\infty))_0&\hflcourte{}{}& (\limproj_n H^1(\G_n,E_\infty))_0 &\hflcourte{}{}& (\limproj_n H^1(\G_n,\limind B_m))_0\cr
     &           &\vflcourte{\wr}{}&                        &\vflcourte{\wr}{}&                                &\vflcourte{\wr}{}&&\cr
 &\hflcourte{}{}&t_\Zp \limind Coker\, k_{n,C}&\hflcourte{}{}&  t_\Zp\limind Coker\, k_{n,E}&\hflcourte{}{}&         \limind Coker\, k_n}$$
 \vspace{-1cm}
$$\diagram{ &&&&&&&&\hflcourte{}{}&         \limproj H^2(\G_n,C_\infty)&\surjfl{}{}&\limproj H^2(\G_n,E_\infty)\cr
            &&&&&&                 && &\vflcourte{\wr}{}&                                &\vflcourte{\wr}{}\cr
             &&&&&& && \hflcourte{}{}&               \limind Cok\, k_{n,C}\otimes \Qp/\Zp&\surjfl{}{}&\limind Cok\, k_{n,E}\otimes
             \Qp/\Zp\cr}$$}
\noindent Le première flèche verticale est surjective, et son noyau vaut $\limproj Coker\, j_{n,C}$.
\end{cor}
\noindent Preuve: On applique \ref{troncaturejk} 2. au triangle distingué tautologique $C\rightarrow E\rightarrow B\rightarrow C[1]$, avec $q=1$.
L'hypothèse de 2. $(i)$ est satisfaite grâce à Leopoldt: en effet $\limproj H^1(C_{j_B})=\limproj Coker\, j_n$ se stabilise et $(B_\infty)^{\G_n}$
est fini, par \ref{remAB}.

Dans l'énoncé, on a remplacé $H^0(\Delta_j(C))=Coker\, \alpha_j=Coker\, \tilde \alpha_j$ par son quotient sans torsion (cf \ref{resEC} 2. $(i)$),
c'est pourquoi la première flèche verticale n'est pas injective a priori.
\begin{flushright}$\square$\end{flushright}

\begin{rem} \label{remqln} Ce résultat recouvre à la fois les th. 3.7 (descente) et 4.5 (codescente) de \cite{NL}, et montre que l'hypothèse
concernant la nullité de $\limproj Coker\, j_{n,C}$ ($(DG)$ dans la terminologie de \emph{loc. cit.}) ne joue aucun rôle dans cette étude.
\end{rem}

\section{Appendice: Construction explicite du complexe dualisant}
Dans le cas d'un groupe $\G$ de dimension cohomologique finie pour lequel $\L:=\Zp[[\G]]$ est noethérien, on construit le complexe dualisant pour la
cohomologie galoisienne de la catégorie ${_\G\C_{tor}}$ des $\G$-modules discrets de $p$-torsion. La construction, basée sur les propriétés des
foncteurs d'induction (cf \ref{induction}) d'une part, et le rapport entre cohomologie galoisienne et foncteurs $Ext$ d'autre part, présente les
avantages suivants:

- On obtient une formule explicite pour le complexe dualisant, et pas seulement pour sa cohomologie.

- Au lieu d'une collection d'isomorphismes fonctoriels dans $D^b(Ab)$, exprimant que les foncteurs
$Hom_\Zp(H^q(\G_n,-),\Qp/\Zp):{_\G\C}_{tor}\rightarrow Ab$ sont tous représentables par le même objet de $D(\G)\in D^b(_\G\C)$, on établit un
isomorphisme fonctoriel dans $D^b({_{\underline\G}\C})$ (ce qui est plus précis), exprimant la ``représenta\-bilité'' du foncteur
$Hom_\Zp(R\underline\G(-),\Qp/\Zp)$ (cf. \ref{compdu} et \ref{compdurem} pour l'énoncé exact). Cette précision est essentielle en vue des passages à
la limites. \vspace{0.2cm}

Comme toujours, on regarde ${_\G\C}$ comme une sous-catégorie pleine de ${_\L\C}$. Lorsqu'on écrit $RHom_\L$ ou $R\underline{Hom}$, il s'agit
toujours des foncteurs obtenus par dérivation à droite sur ${_\L\C}$ (et non sur ${_\G\C}$).

\begin{lem}\label{discdep} $ $

$(i)$ Soit $B$ un objet injectif de $_\G\C$. Si $B$ est de $\Zp$-torsion, alors le foncteur $Hom_\L(-,B):{_\L\C}_{tf}\rightarrow {_\Zp\C}$ est exact.

$(ii)$ Soit $B\in Kom^+(_\L\C)$ un complexe dont les objets vérifient $(i)$, alors les flèches canoniques de $D^+({_\Zp\C})$ (resp.
$D^+({_{\underline\G}\C})$ ou $D^+(\C_{\underline\G})$):
$$Hom_\L(A,B)\rightarrow RHom_\L(A,B)\hspace{1cm}\hbox{(resp.}\hspace{0.2cm} \underline{Hom}(A,B)\rightarrow R\underline{Hom}(A,B)
\hspace{0.2cm}\hbox{)}$$ \noindent
sont des isomorphismes pour tout $A\in Kom^-({_\L\C})_{tf}$.
\end{lem}

\noindent Preuve: $(i)$ La sous-catégorie de ${_\G\C}$ formée des objets de $\Zp$-torsion s'identifie à celle des $\L$-modules topologiques discrets,
et ${_\L\C}_{tf}$ à une sous-catégorie pleine des $\L$-modules topologiques compacts. Le résultat provient alors de \cite{Br} lemma 2.2.

$(ii)$ découle facilement de $(i)$.
\begin{flushright}$\square$\end{flushright}

On utilise maintenant les foncteurs d'induction (normique et discrète) définis dans les préliminaires (cf. \ref{not1}).

\begin{lem} \label{inddu} $ $

Pour $(A,B,C)\in {_\G\C}\times {_\Zp\C}\times {_\L\C}$, il y a dans ${_{\underline\G}\C}$ un morphisme trifonctoriel
$$Hom_\Zp(A,B)\underline{\otimes}C\longrightarrow \underline{Hom}(A,Coind_\G(B)\otimes_\L C)$$
c'est un isomorphisme si $C$ est projectif de type fini.
 \end{lem}
Preuve: On a une série d'isomorphismes, expliqués ci-dessous:

$$\begin{array}{rcl}Hom_\Zp(A,B)\underline{\otimes}C&\mathop\simeq\limits^{(1)}& \underline{Ind}(Hom_\Zp(A,B))\otimes_\L C\\
 &\mathop\simeq\limits^{(2)}& Hom_\Zp(\underline{Coind}(A),B)\otimes_\L C\\
  &\mathop\simeq\limits^{(3)}& Hom_\Zp(\underline{Ind}(A),B)\otimes_\L C\\
 & \mathop\simeq\limits^{(4)} & Hom_\L(\underline{Ind}(A),Coind_\G(B))\otimes_\L C\end{array}$$

\noindent $(1)$ est simplement l'isomorphisme d'induction, pour $\underline{Ind}:{\C_\L}\rightarrow {_{\underline\G}\C_\L}$ (noter que $Hom_\Zp(A,B)$
est un $\L$-module à droite, via $A$). \vspace{0.2cm}

\noindent $(2)$ provient de l'isomorphisme d'évaluation suivant, dans ${_{\underline\G}\C_\L}$:
$$Hom_\Zp(A,B)\underline\otimes \L\simeq Hom_\Zp(\underline{Hom}(\L,A),B)$$

\noindent $(3)$ résulte de \ref{indvscoind}, qui identifie les foncteurs $\underline{Ind}$ et $\underline{Coind}:{_\L\C}\rightarrow
{_\L\C_{\underline\G}}$. \vspace{0.2cm}

\noindent $(4)$ Par hypothèse, $\G$ agit discrètement sur $A$. Comme $\G_n$ est distingué, on voit tout de suite que l'action à gauche de $\G\subset
\L$ sur $\L\otimes_{\Zp[[\G_n]]} A$ via le facteur $\L$ est discrète aussi. Munissant $B$ d'une action triviale de $\G$, on obtient alors un
isomorphisme d'induction discrète:
$$Hom_\Zp(\L\otimes_{\Zp[[\G_n]]} A,B)\simeq Hom_\L(\L\otimes_{\Zp[[\G_n]]} A,Coind_\G B)$$
Celui-ci a lieu dans ${_{G_n}\C_\L}$, si l'on fait agir $G_n$ à gauche via la conjugaison à droite sur $\L\otimes_{\Zp[[\G_n]]} A$, et $\L$ à droite
via le facteur $\L$ (resp, $Coind_\G B$) sur le premier (resp. second) terme (dans le contexte du paragraphe intitulé ``induction discrète'' dans les
préliminaires, cette action de $\L$ doit être interprétée comme une action de $\G$ par conjugaison). D'où un isomorphisme dans
${_{\underline\G}\C_\L}$, en faisant varier $n$; $(4)$ suit. \vspace{0.2cm}

Par ailleurs, il y a dans ${_{\underline\G}\C}$ un isomorphisme d'induction, avec cette fois $\underline{Ind}:{_\L\C}\rightarrow
{_\L\C_{\underline\G}}$:
$$\underline{Hom}(A,Coind_\G(B)\otimes_\L C)\simeq Hom_\L(\underline{Ind}(A),Coind_\G(B)\otimes_\L C)$$

La flèche de l'énoncé provient alors finalement de la flèche naturelle de multiplication:
$$Hom_\L(\underline{Ind}(A),Coind_\G(B))\otimes_\L C\rightarrow Hom_\L(\underline{Ind}(A),Coind_\G(B)\otimes_\L C)$$
C'est visiblement un isomorphisme lorsque $A$ est projectif de type fini.
\begin{flushright}$\square$\end{flushright}

Nous sommes maintenant en mesure d'établir le théorème de dualité.

\begin{thm} \label{compdu} \emph{(Complexe dualisant)}

Soit $D(\G):=\limind \Qp/\Zp\underline\Ltens\Zp \in D^b({_\G\C})$. \vspace{0.1cm}

\noindent Si $A\in D^b(_\G\C)$, il y a dans $D^b({_{\underline\G}\C})$ un isomorphisme fonctoriel:
$$RHom_\Zp(A,\Qp/\Zp)\underline\Ltens\Zp\mathop\rightarrow\limits^\sim R\underline{Hom}(A,D(\G))$$
\end{thm}

\noindent Preuve de \ref{compdu}: Soit $P\rightarrow \Zp$ une résolution parfaite de $\Zp$ dans ${_\L\C}$. D'après \ref{inddu}, on a dans
$Kom^b({_{\underline\G}\C})$ un isomorphisme fonctoriel, pour $A\in Kom^-({_\G\C})$ (avec les conventions de signes adéquates pour la formation des
complexes totaux...):
$$Hom_\Zp(A,\Qp/\Zp)\underline\otimes P\mathop\rightarrow\limits^\sim\underline{Hom}(A,Coind_\G(\Qp/\Zp)\otimes_\L P)$$ Fixant un quasi-isomorphisme
$Coind_\G(\Qp/\Zp)\otimes_\L P\rightarrow I$ où $I\in Kom^+(_\L\C)$ est à objets injectifs, on obtient alors un morphisme composé, fonctoriel dans
$Kom^+({_{\underline\G}\C})$: \begin{eqnarray}\label{eq11}Hom_\Zp(A,\Qp/\Zp)\underline\otimes P\mathop\rightarrow\limits^\sim
\underline{Hom}(A,Coind_\G(\Qp/\Zp)\otimes_\L P)\mathop\rightarrow\limits^{?} \underline{Hom}(A,I)\end{eqnarray} Comme $Coind_\G(\Qp/\Zp)\otimes_\L
P=\limind \underline{Coind}(\Qp/\Zp)\otimes_\L P\simeq \limind \Qp/\Zp\underline\otimes P$ représente $D(\G)$, et que $I$ en est une résolution
injective dans ${_\L\C}$, on voit que le but de $(\ref{eq11})$ ci-dessus représente $R\underline{Hom}(A,D(\G))$. Reste donc à montrer que la flèche
marquée ``$?$'' est un quasi-isomorphisme.

Puisque $Coind_\G(\Qp/\Zp)$ est un injectif de ${_\G\C}$, de $\Zp$-torsion, c'est que les objets de $B:=Coind_\G(\Qp/\Zp)\otimes_\L P$ vérifient les
conditions de \ref{discdep} $(i)$. Si de plus on suppose $A$ de type fini, \ref{discdep} $(ii)$ affirme que ``$?$'' est un quasi-isomorphisme et cela
termine la preuve dans ce cas. Le cas général se ramène au cas où $A$ est de type fini, par un argument de passage à la limite. Détaillons. Soit
$(A_i)\in Kom^b({{_\G\C}_{tf}}^{I})=Kom^b({{_\G\C}_{tf}})^{I}$ un système inductif de complexes de type fini tel que $A=\limind A_i$, et
$Cofl^*((A_i))\rightarrow (A_i)$ la résolution coflasque envisagée dans la preuve de \ref{proplimproj} 4. ($Cofl^*((A_i))\in Kom^-({_\G\C}^{I})$).
Par fonctorialité la flèche $(\ref{eq11})$ donne un diagramme commutatif de $Kom^+({_{\underline\G}\C}^{I^\circ})$: \vspace{-0.8cm}
$$\diagram{Hom_\Zp((A_i),\Qp/\Zp)\underline\otimes P&\hfl{1}{}&\underline{Hom}((A_i),I)\cr \vflcourte{2}{}&&\vflcourte{3}{}&&\cr
Hom_\Zp(Cofl^*((A_i)),\Qp/\Zp)\underline\otimes P&\hfl{}{}& \underline{Hom}(Cofl^*((A_i)),I)\cr  \vflcourte{4}{}&&\vflcourte{5}{}&&\cr
Fl^*((Hom_\Zp(A_i),\Qp/\Zp)\underline\otimes P)&\hfl{6}{}& Fl^*(\underline{Hom}(A_i,I))}$$ \vspace{-1.1cm}

\noindent Dans celui-ci:

- La flèche $1$ est un qis d'après ce qu'on a déjà démontré, puisque les $A_i$ sont de type fini.

- Les flèches $2$ et $3$ sont des qis, puisque les foncteurs $Hom_\Zp(-,\Qp/\Zp)\underline\otimes P$ et $\underline{Hom}(-,I)$ sont exacts.

- Les flèches $4$ et $5$ sont des isomorphismes, puisque les deux foncteurs précédents transforment sommes en produits.

En conséquence, la flèche $6$ est elle aussi un qis. Ne reste plus qu'à appliquer le foncteur $\limproj_I$: on obtient alors un diagramme semblable
dans lequel:

- les flèches $2,3$ sont des qis, car les foncteurs $\limproj\circ Hom_\Zp(-,\Qp/\Zp)\underline\otimes P\simeq Hom_\Zp(-,\Qp/\Zp)\underline\otimes
P\circ \limind$ et $\limproj\circ \underline{Hom}(-,I)\simeq \underline{Hom}(-,I)\circ \limind$ sont exacts.

- Les flèches $4$ et $5$ sont toujours des isomorphismes.

- La flèche $6$ est un qis, car les objets des deux complexes inférieurs sont acycliques pour $\limproj$.

\noindent En conséquence, la flèche $Hom_\Zp(A,\Qp/\Zp)\underline\otimes P\rightarrow \underline{Hom}(A,I)$, déduite de $1$ par $\limproj_I$, est
elle aussi un qis, et cela termine la preuve dans le cas général. \begin{flushright}$\square$\end{flushright}

\begin{rem} \label{compdurem} $(i)$ A l'aide du lemme \ref{lem1}, on peut réécrire la formule défi\-nis\-sant $D(\G)$ de la manière suivante:
$$D(\G)\simeq \limind_{n,k} RHom_\Zp(R\underline\G(\Z/p^k),\Qp/\Zp)$$ On retrouve la formule de \cite{Se2} A.2 prop. 3.1, 5. c) en prenant la
cohomologie.

$(ii)$ Encore à l'aide du lemme \ref{lem1}, on voit que $D(\G)$ est effectivement le complexe dualisant pour la catégorie des $\G$-modules de
discrets de $p$-torsion, au sens de \cite{Se2}, annexe de Verdier. Mieux: pour $A\in {_\G\C_{tor}}$ il y a un isomorphisme fonctoriel dans
$D^b(_{\underline\G}\C)$ (ou dans $D^b(\C_{\underline\G})$, en inversant l'action, ce qu'on n'a pas fait jusqu'ici pour des raisons d'écritures liées
au foncteurs d'induction): \vspace{-0.2cm}
$$RHom_\Zp(R\underline\G(A),\Qp/\Zp)\simeq R\underline{Hom}(A,D(\G))$$
\vspace{-0.5cm}

$(iii)$ Si $A\in D^b(_\G\C)$, alors \ref{compdu} appliqué à $A\Ltens_\Zp\Qp/\Zp$ donne, dans $D^b({_{\underline\G}\C})$ (ou
$D^b(\C_{\underline\G})$):\vspace{-0.2cm}$$\begin{array}{rcl}RHom_\Zp(A,\Zp)\underline\Ltens \Zp & \simeq & RHom_\Zp(A\Ltens_\Zp
\Qp/\Zp,\Qp/\Zp)\underline\Ltens \Zp\\ & \mathop\simeq\limits^{7.3}& R\underline{Hom}(A\Ltens_\Zp\Qp/\Zp,D(\G))\\ &\simeq &
R\underline{Hom}(A,RHom_\Zp(\Qp/\Zp,D(\G)))\end{array}$$ le premier et troisième isomorphismes étant obtenus par adjonction.
\end{rem}

\end{document}